\documentclass[11pt]{amsart} 
\usepackage{amsmath} \usepackage{amsxtra}
\usepackage{amscd} \usepackage{amsthm}  
\usepackage{amsfonts}\usepackage{amssymb} 
\usepackage{eucal}\usepackage{bm}
\usepackage{graphicx} 
\newcommand{\cA}{{\mathcal A}}
\newcommand{\cB}{{\mathcal B}}
\newcommand{\cC}{{\mathcal C}}

\newcommand{\cH}{{\mathcal H}}
\newcommand{\cE}{{\mathcal E}}
\newcommand{\cG}{{\mathcal G}}

\newcommand{\cJ}{{\mathcal J}}
\newcommand{\cO}{{\mathcal O}}
\newcommand{\cL}{{\mathcal L}}
\newcommand{\cM}{{\mathcal M}}

\newcommand{\cF}{{\mathcal F}}
\newcommand{\cK}{{\mathcal K}}
\newcommand{\cP}{{\mathcal P}}

\newcommand{\cR}{{\mathcal R}}
\newcommand{\cS}{{\mathcal S}}
\newcommand{\cT}{{\mathcal T}}
\newcommand{\cU}{{\mathcal U}}
\newcommand{\cV}{{\mathcal V}}
\newcommand{\cW}{{\mathcal W}}
\newcommand{\cX}{{\mathcal X}}

\newcommand{\cXL}{{\mathcal XL}}

\renewcommand{\AA}{{\mathbb A}}
\newcommand{\BB}{{\mathbb B}}
\newcommand{\HH}{{\mathbb H}}
\newcommand{\GG}{{\mathbb G}}
\newcommand{\NN}{{\mathbb N}}
\newcommand{\ZZ}{{\mathbb Z}}
\newcommand{\QQ}{{\mathbb Q}}

\newcommand{\PP}{{\mathbb P}}
\newcommand{\OO}{{\mathbb O}}
\newcommand{\TT}{{\mathbb T}}

\newcommand{\UU}{{\mathbb U}}

\newcommand{\gp}{\mathfrak{p}}
\newcommand{\gq}{\mathfrak{q}}

\newcommand{\gt}{\mathfrak{t}}
\newcommand{\gr}{\mathfrak{r}}

\newcommand{\bH}{\mathbf{H}}

\newcommand{\on}{\operatorname}

\newcommand{\Rep}{{\on{Rep}}}
\newcommand{\Sch}{{\on{Sch}}}

\newcommand{\Qlb}{\mathbb{\bar Q}_\ell}
\newcommand{\Gm}{\mathbb{G}_m}
\newcommand{\Ga}{\mathbb{G}_a}
\newcommand{\A}{\mathbb{A}}

\newcommand{\toup}[1]{\stackrel{#1}{\to}}
\newcommand{\hook}[1]{\stackrel{#1}{\hookrightarrow}}

\newcommand{\getsup}[1]{\stackrel{#1}{\gets}}
\newcommand{\Sp}{\on{\mathbb{S}p}}
\newcommand{\Spin}{\on{\mathbb{S}pin}}
\newcommand{\GSpin}{\on{G\mathbb{S}pin}}
\newcommand{\GSp}{\on{G\mathbb{S}p}}
\newcommand{\IC}{\on{IC}}

\newcommand{\Hom}{\on{Hom}}

\newcommand{\Sym}{\on{Sym}}
\newcommand{\SO}{\on{S\mathbb{O}}}
\newcommand{\GO}{\on{G\mathbb{O}}}

\newcommand{\Aut}{\on{Aut}}

\newcommand{\RG}{\on{R\Gamma}}

\newcommand{\Pic}{\on{Pic}}

\newcommand{\Bun}{\on{Bun}}

\newcommand{\Bunt}{\on{\widetilde\Bun}}

\newcommand{\rk}{\on{rk}}

\newcommand{\Spec}{\on{Spec}}

\newcommand{\supp}{\on{supp}}

\newcommand{\Gr}{\on{Gr}}
\newcommand{\Grb}{\overline{\Gr}}

\newcommand{\GL}{\on{GL}}

\newcommand{\Fun}{{\on{Fun}}}

\newcommand{\pr}{\on{pr}}
\newcommand{\id}{\on{id}}
\newcommand{\tr}{\on{tr}}
\newcommand{\QED}{$\square$} 
\newcommand{\Fq}{\mathbb{F}_q}  
\newcommand{\Fp}{\mathbb{F}_p}  % what for??????    
\newcommand{\iso}{{\widetilde\to}}

\newcommand{\comp}{\circ}
\newcommand{\Four}{\on{Four}}
\renewcommand{\H}{{\on{H}}}   %cohomologies
\newcommand{\DD}{\mathbb{D}}  %for duality
\newcommand{\D}{\on{D}}       %for derived categories     
\newcommand{\DP}{\on{DP}}
\newcommand{\wt}{\widetilde}

\newcommand{\select}[1]{{\it{#1}}}

\renewcommand{\P}{{\on{P}}}

\newcommand{\<}{\langle}
\renewcommand{\>}{\rangle}

\newcommand{\Loc}{\on{Loc}}
\newcommand{\Lie}{\on{Lie}}
\newcommand{\Sph}{\on{Sph}}
\newcommand{\Res}{\on{Res}}
\newcommand{\gRes}{\on{gRes}}

\newcommand{\ttimes}{\tilde\times}

\newcommand{\act}{\on{act}}
\newcommand{\dimrel}{\on{dim.rel}}
\newcommand{\Funct}{\on{Funct}}

\newcommand{\SL}{\on{SL}}

\newcommand{\tboxtimes}{\,\tilde\boxtimes\,}

\newcommand{\Ind}{\on{Ind}}

\newcommand{\Levi}{\on{Levi}}
\newcommand{\ra}{\rightarrow}
\newcommand{\la}{\leftarrow}

\newcommand{\Weil}{\on{{\mathcal{W}}eil}}
\newcommand{\Char}{\on{Char}}
\newcommand{\Cr}{\on{Cr}}

\newcommand{\Out}{\on{{\mathbb O}ut}}
\newcommand{\cind}{\on{c-ind}}
\newcommand{\GSO}{\on{GS\mathbb{O}}}
\newcommand{\colim}{\on{colim}}
\newcommand{\DGCat}{\on{DGCat}}
\newcommand{\DG}{\on{DG}}
\newcommand{\Spc}{\on{Spc}}
\newcommand{\PreStk}{\on{PreStk}}
\newcommand{\Cat}{\on{Cat}}

\newtheorem{Cor}[subsubsection]{Corollary}
\newtheorem{Lm}[subsubsection]{Lemma}

\newtheorem{Pp}[subsubsection]{Proposition}
\newtheorem{Con}[subsubsection]{Conjecture}
\newtheorem{Th}[subsubsection]{Theorem}
\newtheorem{Def}[subsubsection]{Definition}
\newtheorem{Rem}[subsubsection]{Remark}

\newtheorem{PpA}{Proposition}    %%check this!

\newtheorem{PpAA}{Proposition}
\newenvironment{Prf}{\par\noindent {\it Proof }}{\QED}
\newcommand{\Step}[1]{\par\noindent{\bf Step {#1}}.}

\newcommand{\nc}{\newcommand}
\nc{\renc}{\renewcommand}
\nc{\ssec}{\subsection}
\nc{\sssec}{\subsubsection}

\begin{document}

\author{Sergey Lysenko}
\title{Geometric theta-lifting for the dual pair $\GSp_{2n}$, $\GSO_{2m}$}
\address{Institut Elie Cartan Lorraine, Universit\'e de Lorraine, 
 B.P. 239, F-54506 Vandoeuvre-l\`es-Nancy Cedex, France}
\email{Sergey.Lysenko@univ-lorraine.fr}
\thanks{\select{Acknowledgements.} I am grateful to V. Lafforgue for regular and stimulating discussions.}
\begin{abstract}
% \noindent{\scshape Abstract}\hskip 0.8 em 
Let $X$ be a smooth projective curve over an algebraically closed field of characteristic $>2$. Consider the dual pair $H=\GSO_{2m}$, $G=\GSp_{2n}$ over $X$, where $H$ splits over an \'etale two-sheeted covering $\pi:\tilde X\to X$. Write $\Bun_G$ and $\Bun_H$ for the stacks of $G$-torsors and $H$-torsors on $X$. We show that for $m\le n$ (respectively, for $m>n$) the theta-lifting functor $F_G:\D(\Bun_H)\to \D(\Bun_G)$ (respectively, $F_H:\D(\Bun_G)\to \D(\Bun_H)$) commutes with Hecke functors with respect to a morphism of the corresponding 
$L$-groups involving the $\SL_2$ of Arthur. So, this functor realizes the (nonramified) geometric Langlands functoriality for the corresponding morphism of dual groups. 

 As an application, we obtain a particular case of the geometric Langlands conjectures. Namely, we construct the automorphic Hecke eigensheaves on $\Bun_{\GSp_4}$ corresponding to certain endoscopic local systems on $X$. 
\end{abstract} 
\maketitle
%\tableofcontents 

\section{Introduction}

\sssec{} The classical theta correspondence for the dual reductive pair $(\GSp_{2n}, \GSO_{2m})$ is known to satisfy a version of strong Howe duality (cf. \cite{R}). In this paper, which is a continuation of \cite{L2}, we develop the geometric theory of theta-lifting for this dual pair in the everywhere unramified case. 

 The classical theta-lifting operators for this dual pair are as follows. Let $X$ be a smooth projective geometrically connected curve over $\Fq$ (with $q$ odd). Let $F=\Fq(X)$, $\AA$ be the ad\`eles ring of $X$, $\cO$ the integer ad\`eles. Write $\Omega$ for the canonical line bundle on $X$.
Pick a rank $2n$-vector bundle $M$ with symplectic form $\wedge^2 M\to\cA$ with values in a line bundle $\cA$ on $X$. Let $G$ be the group scheme over $X$ of automorphisms of the $\GSp_{2n}$-torsor $(M,\cA)$. 

 Let $\pi:\tilde X\to X$ be an \'etale two-sheeted covering with Galois group $\Sigma=\{1,\sigma\}$. Let $\cE$ be the $\sigma$-anti-invariants in $\pi_*\cO_{\tilde X}$. Fix a rank $2m$-vector bundle $V$ on $X$ with symmetric form $\Sym^2 V\to \cC$ with values in a line bundle $\cC$ on $X$ together with a compatible trivialization $\gamma: \cC^{-m}\otimes\det V\,\iso\, \cE$. This means that $\gamma^2: \cC^{-2m}\otimes (\det V)^2\,\iso\,\cO$ is the trivialization induced by the symmetric form. Let $H$ be the group scheme over $X$ of automorphisms of $V$ preserving the symmetric form up to a multiple and fixing $\gamma$. This is a form of $\GSO_{2m}$, which splits over $\tilde X$. Assume given an isomorphism $\cA\otimes\cC\,\iso\,\Omega$. 
 
 Let $G_{2nm}$  the group scheme of automorphisms of $M\otimes V$ preserving the symplectic form $\wedge^2(M\otimes V)\to\Omega$. Write $GH\subset G\times H$ for the group subscheme over $X$ of pairs $(g,h)$ such that $g\otimes h$ acts trivially on $\cA\otimes\cC$. The metaplectic cover $\wt G_{2nm}(\AA)\to G_{2nm}(\AA)$ splits naturally after restriction under $(GH)(\AA)\to G_{2nm}(\AA)$.  Let $S$ be the corresponding Weil representation of $GH(\AA)$. The space $S^{(GH)(\cO)}$ has a distinguished nonramified vector $v_0$. If $\theta: S\to\Qlb$ is a theta-functional then $\phi_0: (GH)(F)\backslash (GH)(\AA)/ (GH)(\cO)\to\Qlb$ given by $\phi_0(g,h)=\theta((g,h)v_0)$ is the classical theta-function. The theta-lifting operators 
$$
F_G: \Funct(H(F)\backslash H(\AA)/H(\cO))\to  \Funct(G(F)\backslash G(\AA)/G(\cO))
$$
and
$$
F_H:  \Funct(G(F)\backslash G(\AA)/G(\cO))\to \Funct(H(F)\backslash H(\AA)/H(\cO))
$$
are the integral operators with kernel $\phi_0$ for the diagram of projections
$$
\begin{array}{ccccc}
&& (GH)(F)\backslash (GH)(\AA)/ (GH)(\cO)\\
&  \swarrow\lefteqn{\scriptstyle \gq} && \searrow\lefteqn{\scriptstyle \gp}\\
H(F)\backslash H(\AA)/H(\cO) &&&&  G(F)\backslash G(\AA)/G(\cO)
\end{array}
$$

\sssec{Main result} The following claim would be an analog of a theorem of Rallis \cite{Ral} for similitude groups (the author has not found its proof in the litterature). If $m\le n$ (resp., $m>n$) then  $F_G$ (resp., $F_H$) commutes with the actions of global Hecke algebras $\cH_G$, $\cH_H$ with respect to certain homomorphism $\cH_G\to\cH_H$ (resp., $\cH_H\to \cH_G$). 

Our main result is a geometric version of this claim (cf. Theorem~\ref{Th_theta-lifting_functors-general}). Its precise formulation in the geometric setting involves the $\SL_2$ of Arthur (or rather its maximal torus). Write $\Bun_G$ for the stack of $G$-torsors on $X$, similarly for $H$. We define the theta-lifting functors $F_G: \D(\Bun_H)\to\D(\Bun_G)$ and $F_H: \D(\Bun_G)\to\D(\Bun_H)$ between the $\DG$-categories of $\Qlb$-sheaves on these stacks. Theorem~\ref{Th_theta-lifting_functors-general} claims that for $m\le n$ (resp., $m>n$) the functor $F_G$ (resp., $F_H$) commutes with the actions of Hecke functors with respect to a suitable morphism of $L$-groups $H^L\times\Gm\to G^L$ (resp., $G^L\times\Gm\to H^L$). 
In the particular case $n=m$ (resp., $m=n+1$) the above homomorphism is trivial on the $\Gm$-factor. This establishes a special case of the geometric Langlands functoriality.

 This extends a similar result for the pair $(\Sp_{2n}, \SO_{2n})$ from \cite{L2} in two directions. On one hand, we consider the similitude groups, and on the other hand, we allow our groups to split on an \'etale degree two cover of $X$. To take into account the case of non split groups, we propose in Section~\ref{Sect_Twisted_setting_Hecke_functors} a general setting for the geometric Langlands program for groups, which split on an \'etale Galois cover of the curve $X$. 
 
\sssec{Applications} There are two striking applications of our Theorem~\ref{Th_theta-lifting_functors-general}. Both provide proofs of some particular cases of (some version of) the geometric Langlands conjecture for $G=\GSp_4$ (cf. Conjecture~\ref{Con_2.2.3}).

 For the first application, consider an irreducible rank two smooth $\Qlb$-sheaf $\tilde E$ on $\tilde X$ equipped with an isomorphism $\pi^*\chi\,\iso\,\det \tilde E$, where $\chi$ is a smooth $\Qlb$-sheaf on $X$ of rank one. Then $\pi_* (\tilde E^*)$ is equipped with a natural symplectic form $\wedge^2(\pi_* \tilde E^*)\to\chi^{-1}$, so can be viewed as a $\check{G}$-local system $E_{\check{G}}$ on $X$, where $\check{G}$ is the Langlands dual group over $\Qlb$. We construct the automorphic sheaf $K$ on $\Bun_G$, which is a $E_{\check{G}}$-Hecke eigensheaf (cf. Corollary~\ref{Cor_1}). This partially establishes (\cite{L}, Conjecture~2).

The second application, which is one of our main motivations, is a construction of automorphic sheaves on $\Bun_G$ in \cite{L3} attached to $\check{G}$-local systems on $X$, whose standard representations are irreducible local systems of rank four on $X$. It owes its existence to the main result of this paper. 

\ssec{Informal comments on proofs}  

\sssec{} Our methods extend those of \cite{L2}, the global results are derived from the corresponding local ones. Let $F_x$ be the completion of $F$ at $x\in X$, $\cO_x\subset F_x$ the ring of integers. Remind that $S\,\iso\, \otimes'_{x\in X} S_x$ is the restricted tensor product of local Weil representations. The geometric analog of the 
$(GH)(F_x)$-representation $S_x$ is the Weil category $W(\wt\cL_d(W_0(F_x)))$ (cf. Sections~\ref{Sect_Background_Weil}, \ref{Sect_3.2}). It was originally introduced in \cite{LL}. However, the knowledge of the action of the nonramified Hecke algebra of $GH$ on $S_x^{(GH)(\cO_x)}$ is not sufficient to establish Theorem~\ref{Th_theta-lifting_functors-general}. 

To get the actions of the whole local nonramified Hecke algebras $\cH_{H, x}$, $\cH_{G, x}$, we essentially have to consider the compactly induced representation 
\begin{equation}
\label{repres_bar_S_x}
\bar S_x=\cind_{(GH)(F_x)}^{(G\times H)(F_x)} S_x
\end{equation} 
and the space of invariants $(\bar S_x)^{G(\cO_x)\times H(\cO_x)}$ (compare with \cite{R}). We introduce its geometric analog as a family of $\DG$-categories $\D_{\cT_a}(\wt\cL_d(W_a(F_x)))$ indexed by $a\in\ZZ$ (cf. Section~\ref{Sect_5.1_local_G_H}, \ref{Sect_5.2_stack_a_cXL}). In the local setting our group schemes $G, H$ over the formal disk around $x$ are constant. Let $\GG=\GSp_{2n}, \HH=\GSO_{2m}$ be split, write $\check{\GG}, \check{\HH}$ for their Langlands dual groups. In Section~\ref{Sect_5.3_Hecke_functors} we define the actions of $\Rep(\check{\HH}), \Rep(\check{\GG})$ on the above collection of categories. 

 The category $\D_{\cT_0}(\wt\cL_d(W_0(F_x)))$ contains a distinguished object, the perverse sheaf $S_{W_0(F)}$ introduced in (\cite{LL}, Section~6.5). This is an analog of the (unique up to a multiple) $\Sp(M\otimes V)(\cO)$-invariant vector in $S_x$. We reduce our Theorem~\ref{Th_theta-lifting_functors-general} to Theorem~\ref{Th_local_main}, which is our main local result. It says that for $m\le n$ (resp., $m>n$) the actions of $\Rep(\check{\GG})$ and $\Rep(\check{\HH})$ on $S_{W_0(F)}$ are compatible via a suitable homomorphism $\kappa: \check{\HH}\times\Gm\to \check{\GG}$ (resp., $\kappa: \check{\GG}\times\Gm\to \check{\HH}$). Major part of the paper (Section~\ref{Sect_Dual_pair_GSp_GO}) is devoted to a proof of Theorem~\ref{Th_local_main}.
  
\sssec{}  Our pattern of the proof of Theorem~\ref{Th_local_main} follows that of \cite{L2}.  However, we can not simply apply the general results of (\cite{L2}, Section~4), because our Hecke functors change the index $a\in\ZZ$ of the corresponding categories $\D_{\cT_a}(\wt\cL_d(W_a(F_x)))$. Though the categories $\D_{\cT_a}(\wt\cL_d(W_a(F_x)))$ are well adapted for relations with global applications, the action of Hecke functors on them is not sufficiently explicit for our purposes.
 
  For this reason we introduce suitable Levi subgroups $Q(\GG)\subset\GG, Q(\HH)\subset\HH$, and the corresponding Schr\"odinger models of the Weil category and the induced representation (\ref{repres_bar_S_x}). Their advantage over $\D_{\cT_a}(\wt\cL_d(W_a(F_x)))$ is that the group $\GG(F_x)$ (resp., $\HH(F_x)$) acts not just on the category, but on the spaces itself. 
  
  To achieve this, we consider two different Schr\"odinger models for $Q(\GG)$ and $Q(\HH)$ related by a canonical intertwining functor between them. We introduce for $a\in\ZZ$ the $F_x$-vector spaces $\Upsilon_a(F_x)$, $\Pi_a(F_x)$. For $a=0$ these are lagrangian subspaces in the standard representation of $(GH)(F_x)$, and for other $a\in\ZZ$ they are some twisted version. The corresponding Schr\"odinger models are the $\DG$-categories $\D(\Upsilon_a(F_x))$, $\D(\Pi_a(F_x))$ of sheaves on them. The canonical intertwining operator between them is given by the Fourier transform 
$$
\zeta_a: \D(\Upsilon_a(F_x))\,\iso\, \D(\Pi_a(F_x))
$$
for any $a\in\ZZ$ (cf. Section~\ref{Sect_4.5.3}). They are used to reformulate Theorem~\ref{Th_local_main} in a more convenient form as Theorem~\ref{Th2_local}. 

 Namely, in Definition~\ref{Def_Weil_a} we introduce a collection of abelian categories $\Weil_a$, $a\in\ZZ$, which model the $\GG(\cO_x)\times\HH(\cO_x)$-invariants in (\ref{repres_bar_S_x}) via the above Schr\"odinger models. The perverse sheaf $S_{W_0(F)}$ corresponds to a distinguished object $I_0\in\Weil_0$ defined in Section~\ref{Sect_5.6.3_about_I_0}. Theorem~\ref{Th2_local} says that the $\Rep(\check{\GG})$ and $\Rep(\check{\HH})$-actions on $I_0$ are compatible in the same sense as in Theorem~\ref{Th_local_main}.
  
\sssec{} The main tools in the proof of Theorem~\ref{Th2_local} are, on one hand, the weak Jacquet functors that we introduce in Section~\ref{Sect_Weak_Jacquet} and, on the other hand, a result on the local geometric theta-lifting for the dual pair $(\GL_n,\GL_m)$ (\cite{L2}, Proposition~4 and Corollary~5). Our proof also uses a result of Laumon (\cite{La}, Theorem~1.1.2).    
  
\sssec{} As a byproduct, we obtain some new results at the classical level of functions (Propositions~A.1 and A.2). For $a$ even they reduce to a result from \cite{MVW}, but for $a$ odd they are new and amount to a calculation of $K\times H(\cO_x)$-invariants in the Weil representation of $(GH)(F_x)$, where $K$ is the nonstandard maximal compact subgroup of $G(F_x)$. 

\ssec{Organization} Our main results are formulated in Section~\ref{Sect_Main_resuts}, and the proofs are given in the remaining sections. In Section~\ref{Sect_Root data} we describe explicitly the root data of the groups involved, identify their dual groups, and introduce some related objects used later. In Section~\ref{Sect_Local_theory} we remind the definition of the Weil category from \cite{LL} and explain our geometric approach to the necessaty induction along $\Gm$. In Section~\ref{Sect_Dual_pair_GSp_GO} we specialize to the case of the dual pair $(\GSp_{2n},\GSO_{2m})$ and prove Theorem~\ref{Th_local_main}, which is 
our main local result. In Section~\ref{Sect_Global_theory} we derive our 
global results from the local ones. 

 In Appendix~\ref{Sect_appendix_A} we prove Lemma~\ref{Lm_injectivity_appendixA} at the classical level of functions, it is used in Proposition~\ref{Pp_about_functor_cP}. We also establish Proposition~\ref{PpA.1}, which is not used in the rest of the paper. 
 
\section{Main results}
\label{Sect_Main_resuts}

\ssec{Notation} 
\label{Sect_Notations}
\sssec{} 
\label{Sect_2.2.1_revised}
From now on $k$ denotes an algebraically closed field of characteristic $p>2$, all the schemes (or prestacks) we consider are defined
over $k$ (except in Section~\ref{Sect_finite_field_k_0}). Fix a prime $\ell\ne p$. 

 % Let $X$ be a smooth projective connected curve. 
% Set $F=k(X)$. Write $\Omega$ for the canonical line bundle on $X$. For a closed point $x\in X$ write $F_x$ for the completion of $F$ at $x$, let $\cO_x\subset F_x$ be the ring of integers. Let $D_x=\Spec\cO_x$ denotes the disc around $x$. 
 
 We use the following notations and conventions from (\cite{L2}, Section~2.1). Let $\Sch^{aff}$ (resp., $\Sch^{aff}_{ft}$) be the category of affine schemes (resp., affine schemes of finite type) over $k$. Let $\Spc$ be the $\infty$-category of spaces as introduced in (\cite{GaRo}, ch. I.1, 1.1.2).
 
 We adapt the conventions regarding sheaf theory from (\cite{G5}, Section~1.2) taking as sheaf theory the \'etale $\Qlb$-sheaves on affine schemes of finite type. So, $\DGCat_{cont}$ denotes the $\infty$-category of $\Qlb$-linear cocomplete presentable $\DG$-categories with continuous exact $\Qlb$-linear functors (defined in \cite{GaRo}, ch. I.1, 10.3.3). The functor
\begin{equation}
\label{Shv_sheaf_theory_initial}
Shv: (\Sch^{aff}_{ft})^{op}\to\DGCat_{cont}
\end{equation}
attaches to $S$ the ind-completion of the $\DG$-category of constructible \'etale $\Qlb$-sheaves on $S$. 

 Recall that $\PreStk$ is the category of accessible functors
$$
(\Sch^{aff})^{op}\to\Spc
$$
The full subcategory $\PreStk_{lft}\subset \PreStk$ of prestacks locally of finite type consists of functors preserving filtered colimits. It identifies via restriction with the category of all functors $(\Sch^{aff}_{ft})^{op}\to\Spc$. We define
$$
Shv: (\PreStk_{lft})^{op}\to \DGCat_{cont}
$$
as the right Kan extension of (\ref{Shv_sheaf_theory_initial}). For a map $f: Y\to Y'$ in $\PreStk_{lft}$ it gives a functor $f^!: Shv(Y')\to Shv(Y)$, it admits a left adjoint $f_!: Shv(Y)\to Shv(Y')$. For $f$ ind-schematic of ind-finite type we have the functor $f_*: Shv(Y)\to Shv(Y')$. 

Recall the category $\DGCat^{non-cocmpl}$ defined in (\cite{GaRo}, ch. I.1, 10.3.1), it admits filtered colimits\footnote{The projection $\DGCat^{non-cocmpl}\to 1-\Cat$ preserves filtered colimits, where $1-\Cat$ is the $\infty$-category of $\infty$-categories.}. For $C\in\DGCat_{cont}$ we write $C^c\subset C$ for the full subcategory of compact objects, one has $C^c\in\DGCat^{non-cocmpl}$ naturally. 

 For an algebraic $k$-stack $S$ locally of finite type or an ind-scheme of ind-finite type we also write $\D(S)=Shv(S)$, it is equipped with a perverse t-structure. For an ind-scheme of ind-finite type $S$ we also set $\D(S)^{constr}=\D(S)^c$. 
 
 For an algebraic stack $S$ locally of finite type with an affine diagonal we define $\D(S)^{constr}\subset \D(S)$ as the full subcategory consisting of objects that pull back to an object of $Shv(S')^c$ for any $S'\to S$, where $S'\in\Sch^{aff}_{ft}$. 
 
 Recall that $\DGCat^{non-cocmpl}$ has a canonical involution $\cC\mapsto\cC^{op}$ sending $\cC$ to the opposite category (cf. \cite{GaRo}, ch. I.1, 10.3.2). By (\cite{AGKRRV}, Appendix C), if $S$ is an ind-scheme of ind-finite type or an algebraic stack  locally of finite type with an affine diagonal then the Verdier duality gives an equivalence $\DD: (\D(S)^{constr})^{op}\,\iso\, \D(S)^{constr}$ in $\DGCat^{non-cocmpl}$. 

 Write $\P(S)\subset \D(S)$ for the full subcategory of perverse sheaves. Set 
$\DP(S)=\oplus_{i\in\ZZ}\,  \P(S)[i]\subset \D(S)$. By definition, we let for $K,K'\in \P(S), i,j\in\ZZ$ 
$$
\Hom_{\DP(S)}(K[i], K'[j])=\left\{
\begin{array}{ll}
\Hom_{\P(S)}(K,K'), & \mbox{for}\;\; i=j\\
0, & \mbox{for}\;\; i\ne j
\end{array}
\right.
$$
Write $\P^{ss}(S)\subset\P(S)$ for the full subcategory of semi-simple perverse sheaves. Let $\DP^{ss}(S)\subset\DP(S)$ be the full subcategory of objects of the form $\oplus_{i\in\ZZ} \, K_i[i]$ with $K_i\in\P^{ss}(S)$ for all $i$.  

\sssec{} We freely use the extension of the above sheaf theory to placid ind-schemes as defined in (\cite{G5}, Appendix C). So, for a placid ind-scheme $Y\in\PreStk$ we have $Shv(Y)\in\DGCat_{cont}$ defined in (\cite{G5}, C.2), and for any morphism $f: Y\to Y'$ of placid ind-schemes we have the functor $f_*: Shv(Y)\to Shv(Y')$ in $\DGCat_{cont}$.

 Let $Z$ be a placid scheme written as $\lim_{i\in I^{op}} Z_i$, where $I$ is a filtered small $\infty$-category, $Z_i$ is a scheme of finite type and for $\alpha: i\to j$ in $I$ the map $f_{\alpha}: Z_j\to Z_i$ is affine, surjective and smooth of some relative dimension $d_{\alpha}$. We define $\D(Z)$ as $\colim_{i\in I} \D(Z_i)$ with respect to the transition functors $f_{\alpha}^*[\dimrel(f_{\alpha})]$. We have canonically $\D(Z)\,\iso\, Shv(Z)$, but our notation $\D(Z)$ is supposed to recall that $\D(Z)$ is naturally equipped with a perverse t-structure (while for an arbitrary placid ind-scheme $S$ we need more data to introduce a t-structure on $Shv(S)$, see \cite{G5}, Remark C.2.2). 
 
 We also define $\D(Z)^{constr}\in \DGCat^{non-cocmpl}$ as follows. Recall that for each $i\in I$, $\D(Z_i)$ is compactly generated so that $\Ind(\D(Z_i)^{constr})\,\iso\,\D(Z_i)$ canonically. Moreover, for $\alpha: i\to j$ in $I$ the functor $f_{\alpha}^*:\D(Z_i)\to \D(Z_j)$ preserves compactness, as its right adjoint is continuous. So, we get a diagram $I\to \DGCat^{non-cocmpl}$, $i\mapsto \D(Z_i)^{constr}$.
with the transition functors $f_{\alpha}^*[\dimrel(f_{\alpha})]$. We let 
$$
\D(Z)^{constr}=\mathop{\colim}\limits_{i\in I} \D(Z_i)
$$ 
taken in $\DGCat^{non-cocmpl}$. Then $\D(Z)^{constr}\subset \D(Z)$ is a full subcategory, and as in (\cite{GaRo}, ch. I.1, 7.2.7) one gets $\Ind(\D(Z)^{constr})\,\iso\, \D(Z)$ canonically.  
 
 Let $Y$ be a placid ind-scheme written as $Y\,\iso\, \colim_{j\in J} Y_j$ in $\PreStk$, where $J$ is a small filtered $\infty$-category, $Y_j$ is a placid scheme as above and for $\alpha: i\to j$ in $J$, the transition map $f_{\alpha}: Y_i\to Y_j$ is a placid closed embedding. Then we let $\D(Y)=\lim_{j\in J^{op}} \D(Y_j)$ with respect to the transition functors $f_{\alpha}^!$. Recall also that 
$$
\D(Y)\,\iso\, \mathop{\colim}\limits_{j\in J} \D(Y_j)
$$ 
in $\DGCat_{cont}$ with respect to the transition functors $(f_{\alpha})_!$. 
 
  We have $Shv(Y)\,\iso\, \D(Y)$, but our notation $\D(Y)$ is supposed to recall that $\D(Y)$ is naturally equipped with a perverse t-structure. By definition, the category of connective objects $\D(Y)^{\le 0}$ in $\D(Y)$ for the perverse t-structure is the smallest full subcategory containing $\D(Y_j)^{\le 0}$ for each $j\in J$ and stable under small colimits. 
  
 Define $\D(Y)^{constr}\in \DGCat^{non-cocmpl}$ as follows. For any $\alpha: i\to j$ in $J$, the functor $(f_{\alpha})_!: \D(Y_i)\to \D(Y_j)$ sends the full subcategory $\D(Y_i)^{constr}$ of $\D(Y_i)$ to $\D(Y_j)^{constr}$. We get a diagram $J\to \DGCat^{non-cocmpl}$, $j\mapsto \D(Y_j)^{constr}$. We let 
$$
\D(Y)^{constr}=\mathop{\colim}\limits_{j\in J} \D(Y_j)^{constr}
$$ 
taken in $\DGCat^{non-cocmpl}$. As above, $\D(Y)\,\iso\, \Ind(\D(Y)^{constr})$ canonically in $\DGCat_{cont}$. 
 
 Since the transition functors in the definition of $\D(Y)^{constr}$ commute with the Verdier duality, we get an equivalence $\DD: (\D(Y)^{constr})^{op}\,\iso\, \D(Y)^{constr}$ in $\DGCat^{non-cocmpl}$. 
 
\sssec{} Since we are working over an algebraically closed field, we systematically ignore Tate twists (except in Section~\ref{Sect_finite_field_k_0}, where we work over a finite subfield $k_0\subset k$. In this case we also fix a square root $\Qlb(\frac{1}{2})$ of the sheaf $\Qlb(1)$ over $\Spec k_0$). Fix a nontrivial character $\psi: \Fp\to\Qlb^*$ and denote by $\cL_{\psi}$ the corresponding Artin-Shreier sheaf on $\A^1$.
 
 If $V\to S$ and $V^*\to S$ are dual rank $n$ vector bundles over a stack $S$, we normalize 
the Fourier transform $\Four_{\psi}: \D(V)\to\D(V^*)$ by 
$\Four_{\psi}(K)=(p_{V^*})_!(\xi^*\cL_{\psi}\otimes p_V^*K)[n](\frac{n}{2})$,  
where $p_V, p_{V^*}$ are the projections, and $\xi: V\times_S V^*\to \A^1$ is the pairing.

 For a sheaf of groups $G$ on a scheme $S$, $\cF^0_G$ denotes the trivial $G$-torsor on $S$. For a representation $V$ of $G$ and a $G$-torsor $\cF_G$ on $S$ write $V_{\cF_G}=V\times^G\cF_G$ for the induced vector bundle on $S$. For a morphism of stacks $f:Y\to Z$ denote by $\dimrel(f)$ the function of  connected component $C$ of $Y$ given by $\dim C-\dim C'$, where $C'$ is the connected component of $Z$ containing $f(C)$. 
 
 Let $\cO$ be a complete discrete valuation $k$-algebra, $F$ its fraction field. For a local system $E$ on $X$ we denote by $E^*$ its dual. 

\ssec{Hecke functors} 

\sssec{} 
\label{Sect_2.2.1}
Let $X$ be a smooth connected projective curve. For $r\ge 1$ write $\Bun_r$ for the stack of rank $r$ vector bundles on $X$. The Picard stack $\Bun_1$ is also denoted $\Pic X$. For a connected reductive group $\GG$ over $k$, let $\Bun_{\GG}$ denote the stack of $\GG$-torsors on $X$. 

 Given a maximal torus and a Borel subgroup $\TT\subset \BB\subset \GG$, we write $\Lambda_{\GG}$ (resp., $\check{\Lambda}_{\GG}$) for the coweights (resp., weights) lattice of $\GG$. Let $\Lambda^+_{\GG}$ (resp., $\check{\Lambda}^+_{\GG}$) denote the set of dominant coweights (resp., dominant weights) of $\GG$. Write $\check{\rho}_{\GG}$ (resp., $\rho_{\GG}$) for the half sum of the positive roots (resp., coroots) of $\GG$, $w_0$ for the longuest element of the Weyl group of $\GG$. 
 
 Set $K=k(X)$. For a closed point $x\in X$ let $K_x$ be the completion of $K$ at $x$, $\cO_x\subset K_x$ be its ring of integers. Set $D_x=\Spec \cO_x$, $D_x^*=\Spec K_x$. 
 
 The following notations are borrowed from \cite{L2}. The affine grassmanian is denoted $\Gr_{\GG}=\GG(F)/\GG(\cO)$ and $\Gr_{\GG,x}=\GG(K_x)/\GG(\cO_x)$. The latter is an ind-scheme classifying a $\GG$-torsor $\cF_{\GG}$ on $X$ together with a trivialization $\beta: \cF_{\GG}\mid_{X-x}\,\iso\, \cF^0_{\GG}\mid_{X-x}$. For $\lambda\in\Lambda^+_{\GG}$ we write $\Gr_{\GG, x}^{\lambda}$ for the $\GG(\cO_x)$-orbit through $t^{\lambda}\GG(\cO_x)\in \Gr_{\GG, x}$ and $\Grb^{\lambda}_{\GG,x}\subset\Gr_{\GG,x}$ for its reduced closure.  
 
 For $\theta\in\pi_1(\GG)$ denote by $\Gr_{\GG}^{\theta}$ the connected component of $\Gr_{\GG}$ containing $\Gr_G^{\lambda}$ for any $\lambda\in\Lambda^+_{\GG}$ lying over $\theta$. 
  
 Denote by $\cA^{\lambda}_{\GG}$ the intersection cohomology sheaf of $\Grb^{\lambda}_{\GG}$. Write $\check{\GG}$ for the Langlands dual group to $\GG$, this is a reductive group over $\Qlb$ equipped with the dual maximal torus and Borel subgroup $\check{\TT}\subset\check{\BB}\subset\check{\GG}$. Write $\Sph_{\GG}$ for the category of $\GG(\cO_x)$-equivariant perverse sheaves on $\Gr_{\GG,x}$. This is a tensor category, and one has a canonical equivalence of tensor categories $\Loc: \Rep(\check{\GG})\,\iso\,\Sph_{\GG}$, where $\Rep(\check{\GG})$ is the category of finite-dimensional representations of $\check{\GG}$ over $\Qlb$ (cf. \cite{MV}). The category $\Sph_{\GG}$ is independent of $x\in X$ up to a canonical equivalence by (\cite{MV}, Proposition~2.2).   
 
 For the definition of the Hecke functors 
\begin{equation}
\label{Hecke_functors_basic}
\H^{\la}_{\GG}, \H^{\ra}_{\GG}: \Sph_{\GG}\times \D(\Bun_{\GG})\to \D(X\times\Bun_{\GG})
\end{equation}
we refer the reader to (\cite{L2}, Section~2.2.1). Write $\ast: \Sph_{\GG}\,\iso\, \Sph_{\GG}$ for the covariant equivalence induced by the map $\GG(K_x)\to \GG(K_x)$, $g\mapsto g^{-1}$. In view of $\Loc$, the corresponding functor $\ast: \Rep(\check{\GG})\,\iso\, \Rep(\check{\GG})$ sends an irreducible $\check{\GG}$-module with highest weight $\lambda$ to the irreducible $\check{\GG}$-module with highest weight $-w_0(\lambda)$. For $\lambda\in\Lambda^+_{\GG}$ we also write $\H^{\lambda}_{\GG}(\cdot)=\H^{\la}_{\GG}(\cA^{\lambda}_{\GG},\cdot)$. For $\cS\in\Sph_{\GG}$ one has functorially $\H^{\ra}_{\GG}(\ast\cS, \cdot)\,\iso\, \H^{\la}_{\GG}(\cS,\cdot)$.
 
 Set 
$$
\D\Sph_{\GG}=\oplus_{r\in\ZZ} \; \Sph_{\GG}[r] \subset \D(\Gr_{\GG})
$$ 
As in (\cite{L2}, Section~2.1.2), we equip it with a structure of  a tensor category in such a way that the Satake equivalence extends to an equivalence of tensor categories 
\begin{equation}
\label{def_Loc_gr}
\Loc^{\gr}: \Rep(\check{\GG}\times\Gm)\,\iso\, \D\Sph_{\GG}
\end{equation}
Our convention is that $\Gm$ acts on $\Sph_{\GG}[r]$ by the character $x\mapsto x^{-r}$. Extend $\ast$ to an involution $\ast: \D\Sph_{\GG}\,\iso\,\D\Sph_{\GG}$ by $\ast (K[r])\,\iso\, (\ast K)[r]$ for $K\in\Sph_{\GG}$. 

We extend (\ref{Hecke_functors_basic}) to
$$ 
\H^{\la}_{\GG}, \H^{\ra}_{\GG}: \D\Sph_{\GG}\times \D(\Bun_{\GG})\to \D(X\times\Bun_{\GG}),
$$
where for $\cS\in\Sph_{\GG}$, 
$$
\H^{\la}_{\GG}(\cS[m], \cdot)=\H^{\la}_{\GG}(\cS, \cdot)[m],
$$ 
and similarly for $\H^{\ra}_{\GG}$. In view of $\Loc^{\gr}$, we sometimes replace $\D\Sph_{\GG}$ in the definition of Hecke functors by $\Rep(\check{\GG}\times\Gm)$. 

\sssec{} For the convenience of the reader, we formulate the version of the geometric Langlands conjecture addressed in Section~\ref{Section_autom_sheaves_Bun_GSp4}. For the notion of a Hecke eigensheaf we refer to (\cite{L3}, Definition~2.1.1) or (\cite{GL}, Section~9.5.3). 

\begin{Con} 
\label{Con_2.2.3}
Let $E$ be a $\check{\GG}$-local system on $X$. There is a nonzero $K\in \D(\Bun_{\GG})$, which is a $E$-Hecke eigensheaf.
\end{Con}
\begin{Rem}
\label{Rem_2.2.4}
In the situation of Conjecture~\ref{Con_2.2.3} it is expected that, under some genericity assumptions on $E$, one may find a nonzero $E$-Hecke eigensheaf $K\in\D(\Bun_{\GG})$, which is moreover perverse.
\end{Rem}

\sssec{Geometric restriction functors} 
\label{Sect_Geometric restriction functors} 
Let $\PP\subset\GG$ be a standard parabolic subgroup with Levi quotient $\QQ$. Let $\check{\QQ}\subset\check{\GG}$ be the corresponding Langlands dual group, $Z(\check{\QQ})$ be the center of $\check{\QQ}$. Let $\kappa: \check{\QQ}\times\Gm\to \check{\GG}$ be a homomorphism, whose first component is the natural inclusion, and the second one factors as $\Gm\toup{i_{\kappa}}Z(\check{\QQ})\hook{} \check{\GG}$ for some cocharacter $i_{\kappa}: \Gm\to Z(\check{\QQ})$. As in (\cite{L2}, Section~4.7.2), define the geometric restriction functor $\gRes^{\kappa}: \Sph_{\GG}\to \D\Sph_{\QQ}$ as follows. Note that $i_{\kappa}\in \check{\Lambda}_{\GG}$ is orthogonal to the coroots of $\QQ$. 
For $\theta\in\pi_1(\QQ)$ consider the diagram 
$$
\Gr_{\QQ}^{\theta}\getsup{\gt_{\PP}^{\theta}}\Gr_{\PP}^{\theta}\toup{\gt_{\GG}^{\theta}}\Gr_{\GG}
$$ 
obtained by functoriality from $\QQ\gets\PP\to\GG$. Set 
$$
\gRes^{\kappa}(\cS)=\mathop{\oplus}\limits_{\theta\in\pi_1(\QQ)} (\gt_{\PP}^{\theta})_!(\gt_{\GG}^{\theta})^*\cS[\<\theta, 2(\check{\rho}_{\GG}-\check{\rho}_{\QQ})-i_{\kappa}\>]
$$
The diagram commutes 
$$
\begin{array}{ccc}
\Sph_{\GG} & \toup{\gRes^{\kappa}} & \D\Sph_{\QQ}\\
\uparrow\lefteqn{\scriptstyle \Loc} && \uparrow\lefteqn{\scriptstyle \Loc^{\gr}} \\
\Rep(\check{\GG}) &\toup{\Res^{\kappa}} & \Rep(\check{\QQ}\times\Gm)
\end{array}
$$

 By abuse of notations, for a homomoprhism $\kappa: \check{\QQ}\times\Gm\to \check{\GG}$ we will also denote by $\kappa: \check{\QQ}\times\Gm\to \check{\GG}\times\Gm$ the map $(\kappa, \pr)$, where $\pr:  \check{\QQ}\times\Gm\to\Gm$ is the projection. This convention is used for restriction functors $\gRes^{\kappa}: \D\Sph_{\GG}\to\D\Sph_{\QQ}$, they commute with cohomological shifts. 

\ssec{Twisted setting} 
\label{Sect_Twisted_setting_Hecke_functors}

\sssec{} 
\label{Sect_2.3.1}
The version of the twisted setting of the geometric Langlands proposed below is essentially the same as in (\cite{HNY}, Appendix B) and (\cite{Zh}, Appendix~A). 

Let $\Sigma$ be a finite group, $\pi: \tilde X\to X$ be a $\Sigma$-torsor in \'etale topology. Given a homomorphism $\rho: \Sigma\to\Aut(\GG)$, define the group scheme $G$ on $X$ as $(\GG\times \tilde X)/\Sigma$, where $\Sigma$ acts diagonally. We refer to $G$ as the twisting of $\GG$ by the $\Sigma$-torsor $\pi:\tilde X\to X$. Let $\tilde K$ be the ring of rational functions on $\tilde X$. For a closed point $\tilde x\in\tilde X$ write $K_{\tilde x}$ for the completion of $\tilde K$ at $\tilde x$, $\cO_{\tilde x}\subset K_{\tilde x}$ for its ring of integers, set $D_{\tilde x}=\Spec\cO_{\tilde x}$, $D_{\tilde x}^*=\Spec K_{\tilde x}$. The map $\pi:\tilde X\to X$ yields isomorphisms $D_{\tilde x}\,\iso\, D_{\pi(\tilde x)}$. Write $\Gr_{\GG,\tilde x}$ for the affine grassmanian $\GG(K_{\tilde x})/\GG(\cO_{\tilde x})$. 

 The group $\Sigma$ acts on $\Gr_{\GG}$ on the left by functoriality. For $\sigma\in\Sigma$ we still denote by $\sigma: \Gr_{\GG}\,\iso\,\Gr_{\GG}$ the corresponding isomorphism. It yields a right action of $\Sigma$ on $\Sph_{\GG}$, where $\sigma$ act as 
\begin{equation}
\label{def_of_sigma^*}
\sigma^*: \Sph_{\GG}\iso\Sph_{\GG}
\end{equation} 
 
 Write $\Bun_G$ for the stack of $G$-torsors on $X$. One defines the Hecke functors 
\begin{equation}
\label{Hecke_twisted_functors_2.2}
_{\tilde x}\H^{\la}_G, {_{\tilde x}\H^{\ra}_G}: \Sph_{\GG}\times\D(\Bun_G)\to \D(\Bun_G)
\end{equation}
as in (\cite{L}, Appendix~B). Namely, let $_x\cH_G$ be the Hecke stack classifying $G$-torsors $\cF_G, \cF'_G$ on $X$ and an isomorphism $\cF_G\,\iso\,\cF'_G\mid_{X-x}$. Let $\Gr_{G,x}$ be the ind-scheme classifying a $G$-torsor on $D_x$ with a trivialization over $D_x^*$. Note that over $D_x$ the group scheme $G$ is constant.

We have a diagram
$$
\Bun_G\getsup{h^{\la}} {_x\cH_G}\toup{h^{\ra}}\Bun_G,
$$
where $h^{\la}$ (resp., $h^{\ra}$) sends $(\cF_G,\cF'_G, x)$ to $\cF_G$ (resp., to $\cF'_G$). Let $\Bun_{G, x}$ be the stack classifying $\cF_G\in\Bun_G$ with a trivialization $\cF_G\mid_{D_x}\,\iso\, \cF^0_{G}$.  
One gets the isomorphisms
$$
\id^l, \id^r: {_x\cH_G}\,\iso\, \Bun_{G,x}\times^{G(\cO_x)}\Gr_{G,x}
$$
such that the projection to the first factor corresponds to $h^{\la}, h^{\ra}$ respectively. 

 A choice of $\tilde x\in\tilde X$ over $x=\pi(\tilde x)$ yields an isomorphism $\eta_{\tilde x}: \Gr_{\GG, \tilde x}\,\iso\, \Gr_{G, x}$. For $\cS\in\Sph_{\GG}$ we still denote by $\cS$ the corresponding $G(\cO_x)$-equivariant perverse sheaf on $\Gr_{G,x}$. Thus, to $\cS\in\Sph_{\GG}$, $K\in \D(\Bun_G)$ and $\tilde x$ over $x$ one attaches the twisted external product $(K\tboxtimes\cS)^l$ and $(K\tboxtimes\cS)^r$ on $_x\cH_G$, they are normalized to be perverse for $K,S$ perverse. The functors (\ref{Hecke_twisted_functors_2.2}) are defined by
$$
_{\tilde x}\H^{\la}_G(\cS,K)=h^{\la}_!(K \tboxtimes \!\ast\cS)^r\;\;\;\;\mbox{and}\;\;\;\;
_{\tilde x}\H^{\ra}_G(\cS,K)=h^{\ra}_!(K \tboxtimes \cS)^l
$$
We have canonically $_{\tilde x}\H^{\la}_G(\ast\cS, K)\,\iso\, {_{\tilde x}\H^{\ra}_G}(\cS, K)$. Letting $\tilde x$ move along $\tilde X$, one similarly defines Hecke functors
$$
\H^{\la}_G, \H^{\ra}_G: \Sph_{\GG}\times\D(\Bun_G)\to \D(\tilde X\times\Bun_G)
$$
They are compatible with the tensor structure on $\Sph_{\GG}$ and commute with the Verdier duality (cf. \cite{BG, L2}). As in Section~\ref{Sect_2.2.1}, we extend the above Hecke functors
to the action of $\D\Sph_{\GG}$, which is sometimes replaced by $\Rep(\check{\GG}\times\Gm)$ in view of $\Loc^{\gr}$.

\begin{Rem} 
\label{Rem_2.2.7}
The category of $G$-torsors on $X$ is equivalent to the category of $\GG$-torsors $\cF_{\GG}$ on $\tilde X$ equipped with a compatible system of  isomorphisms $\alpha_{\sigma}: \sigma^*\cF_{\GG}\,\iso\, \cF_{\GG}^{\sigma}$. Here $\cF_{\GG}^{\sigma}$ is the extension of scalars of $\cF_{\GG}$ via $\sigma: \GG\to\GG$ (cf. \cite{L}, Lemma~21). Compatibility means that for $\sigma, \tau\in \Sigma$, the diagram commutes
$$
\begin{array}{ccc}
\tau^*\sigma^*\cF_{\GG} & \toup{\alpha_{\sigma}} & \tau^*(\cF^{\sigma}_{\GG})=(\tau^*\cF_{\GG})^{\sigma}\\
\downarrow\lefteqn{\scriptstyle \alpha_{\sigma\tau}} && \downarrow\lefteqn{\scriptstyle \alpha_{\tau}}\\
\cF_{\GG}^{\sigma\tau}  & = & (\cF^{\tau}_{\GG})^{\sigma}
\end{array}
$$
\end{Rem}
 
\sssec{} 
\label{Sect_for_Out}
For $\sigma\in\Sigma$ consider the diagram
$$
\Gr_{\GG,\tilde x}\toup{\eta_{\tilde x}}\Gr_{G, x}\getsup{\eta_{\sigma\tilde x}} \Gr_{\GG, \sigma\tilde x}
$$ 
By Remark~\ref{Rem_2.2.7}, the functor $(\eta_{\tilde x})^*(\eta_{\sigma\tilde x})_*: \Sph_{\GG}\to \Sph_{\GG}$ identifies with (\ref{def_of_sigma^*}) naturally.  
So, for $\sigma\in \Sigma$ we have 
\begin{equation}
\label{iso_key_for_twisted_setting}
(\sigma\times\id)^*\comp \H^{\la}_G(\cS,\cdot)\,\iso\, \H^{\la}_G(\sigma^*\cS, \cdot)
\end{equation}
functorially in $\cS\in\Sph_{\GG}$.
  
% Let us fix a pinning of $(\check{\GG},\check{\BB}, \check{\TT})$. It yields a canonical section of the surjection $\Aut(\check{\GG})\to \Out(\check{\GG})$ as in (\cite{B}, Section~1.1). 

 From now on we make a simplifying assumption that either $\GG$ is a torus or $\Sigma$ acts trivially on $\GG/[\GG,\GG]$. This is sufficient for our purposes. The above action of $\Sigma$ on $\Sph_{\GG}$ preserves the fibre functor, so gives a homomorphism $\Sigma\to \Aut(\check{\GG})$. 

 Write $\Out(\GG)$ for the group of outer automorphisms of $\GG$. Recall that $\check{\GG}$ is equipped with the maximal torus and Borel subgroups $\check{\TT}\subset \check{\BB}\subset\check{\GG}$. Recall that a pinning is a collection of isomorphisms $\Ga\,\iso\, \UU_{\alpha}$ for each simple root $\alpha$ of $\check{\GG}$, where $\UU_{\alpha}\subset \check{\BB}$ is the corresponding 1-parameter subgroup (\cite{GP}, Expos\'e XXIII, 1.1). By (\cite{HNY}, Lemma~B.3), there is a natural pinning $\dagger$ of $(\check{\TT}, \check{\BB},\check{\GG})$ with the following property. Let $\Aut^{\dagger}(\check{\GG})$ be automorphism group of $\check{\GG}$ fixing this pinning. 
 
 In the case when $\Sigma$ acts trivially on $\GG/[\GG,\GG]$, the corresponding homomorphism $\Aut(\GG)\to \Aut(\check{\GG})$ factors as $\Aut(\GG)\to \Out(\GG)\toup{\tilde\iota} \Aut^{\dagger}(\check{\GG})\subset\Aut(\check{\GG})$, where $\tilde\iota$ is an isomorphism (compare also with (\cite{B}, Section~1.1). Let $\cL_{\det}$ be the determinant line bundle on $\Gr_G$. Recall that $\dagger$ is obtained from the principal nilpotent in $\Lie\check{\GG}$ given by the action of $\cL_{\det}$ on the fibre functor for $\Sph_{\GG}$. In this case the $\Sigma$-action on $\GG$ gives the $\Sigma$-action on $\check{\GG}$ via $\tilde\iota$. 
 
 In any case under our assumptions we get an action of $\Sigma$ on $\Rep(\check{\GG})$ by functoriality. For $V\in\Rep(\check{\GG})$ write $\sigma^*V$ for the representation of $\check{\GG}$, where $g\in\check{\GG}$ acts as $\sigma(g)$. Then the Satake equivalence $\Loc: \Rep(\check{\GG})\,\iso\, \Sph_{\GG}$ is $\Sigma$-equivariant, that is, $\Loc(\sigma^*V)\,\iso\, \sigma^*\Loc(V)$ naturally. Let $\Sigma$ also act on $\check{\GG}\times\Gm$ as the product of the above action on $\check{\GG}$ with the trivial action on $\Gm$. The above extends to a $\Sigma$-equivariant equivalence 
\begin{equation}
\label{Loc_gr_Sigma_equivariant}
\Loc^{\gr}: \Rep(\check{\GG}\times\Gm)\,\iso\,\D\Sph_{\GG}
\end{equation} 

 For a $\check{\GG}$-local system $E$ on $\tilde X$ and $V\in\Rep(\check{\GG})$ write $V_E$ for the twist of $V$ by $E$, this is a local system on $\tilde X$. 

\sssec{} For a $\check{\GG}$-local system $E$ on $\tilde X$ and $\sigma\in\Sigma$ write $E^{\sigma}$ for the $\check{\GG}$-local system on $X$ obtained from $E$ via extension of scalars via $\sigma: \check{\GG}\to \check{\GG}$. Write $LS_{\tilde X}(\check{\GG})$ for the category of $\check{\GG}$-local systems on $\tilde X$. We may twist the notion of a $\check{\GG}$-local system on $X$ using the $\Sigma$-torsor $\pi: \tilde X\to X$ as follows. Recall the homomorphism $\rho: \Sigma\to\Aut(\GG)$.

\begin{Def} A $\pi$-twisted $\check{\GG}$-local system on $\tilde X$ is a datum of a $\check{\GG}$-local system $E$ on $\tilde X$ and a compatible system of isomorphisms $\beta_{\sigma}: \sigma^*E\,\iso\, E^{\sigma}$. The compatibility means that for $\sigma,\tau\in\Sigma$ the diagram commutes
$$
\begin{array}{ccc}
\tau^*\sigma^*E & \toup{\beta_{\sigma}} & \tau^*(E^{\sigma})=(\tau^*E)^{\sigma}\\
\downarrow\lefteqn{\scriptstyle \beta_{\sigma\tau}} && \downarrow\lefteqn{\scriptstyle \beta_{\tau}}\\
E^{\sigma\tau}  & = & (E^{\tau})^{\sigma}
\end{array}
$$
Write $LS^{\pi,\rho}_X(\check{\GG})$ for the category of $\pi$-twisted $\check{\GG}$-local systems on $\tilde X$.
\end{Def}

\sssec{} Note that for $V\in\Rep(\check{\GG})$ and a $\check{\GG}$-local system $E$ on $\tilde X$, one has $\sigma^*(V_E)\,\iso\, V_{\sigma^*E}$ and $V_{E^{\sigma}}\,\iso\, (\sigma^*V)_E$ naturally.  

\begin{Def} Let $K\in\D(\Bun_G)$. Let $E$ be a $\pi$-twisted $\check{\GG}$-local system on $X$. A structure of a $E$-Hecke eigensheaf on $K$ is a collection of isomorphisms 
$$
\gamma_V: \H^{\la}_G(V, K)\,\iso\, V_E[1]\boxtimes K 
$$
compatible with the tensor structure on $\Rep(\check{\GG})$ as in (\cite{L3}, Section 2). We require in addition that for $\sigma\in\Sigma, V\in\Rep(\check{\GG})$ the diagram commutes
$$
\begin{array}{ccccc}
(\sigma\times\id)^*\H^{\la}_G(V, K)& \toup{(\sigma\times\id)^*\gamma_V} & (\sigma^*(V_E))[1]\boxtimes K & \iso & V_{\sigma^*E}[1]\boxtimes K\\
\downarrow &&&& \downarrow\lefteqn{\scriptstyle \beta_{\sigma}}\\
\H^{\la}_G(\sigma^*V, K) & \toup{\gamma_{\sigma^*V}} & (\sigma^*V)_E[1]\boxtimes K & \iso & V_{E^{\sigma}}[1]\boxtimes K,
\end{array}
$$ 
where the left vertical arrow is the isomorphism (\ref{iso_key_for_twisted_setting}). 
\end{Def}

\sssec{} For a category $\cA$ by a left action of $\Sigma$ on $\cA$ we mean a monoidal functor $\Sigma\to \Fun(\cA,\cA)$, (cf. \cite{HTT}, Definition~A.1.3.5). For such an action and $\sigma\in\Sigma$ we get the  functor $act_{\sigma}: \cA\to\cA$, and for $\sigma, \tau\in\Sigma$ the 2-isomorphism $act_{\gamma\sigma}\,\iso\, act_{\gamma}\act_{\sigma}$. The category of invariants $\cA^{\Gamma}$ is the category of collections: $(a\in\cA, \{c_{\sigma}\}_{\sigma\in\Sigma})$, where $c_{\sigma}: act_{\sigma}(a)\,\iso\, a$ is an isomorphism satisfying the corresponding cocyle condition, that is, the above 2-morphism yields for $\sigma,\tau\in\Sigma$ an isomorphism $c_{\tau}\comp act_{\tau}(c_{\sigma})\,\iso\,c_{\tau\sigma}$. 

 A right $\Sigma$-action on $\cA$ is a monoidal functor $\Sigma^{rm}\to \Fun(\cA,\cA)$, where $\Sigma^{rm}$ denotes the group $\Sigma$ with the reversed multiplication. 

\sssec{Example i)} For the above action of $\Sigma$ on $\Rep(\check{\GG})$ one has $\Rep(\check{\GG})^{\Sigma}\,\iso\, \Rep(\check{\GG}\rtimes\Sigma)$ by (\cite{Zh}, Lemma~A.3). 
 
\sssec{Example ii)} Consider the $\Sigma$ action on $LS_{\tilde X}(\check{\GG})$ such that $\sigma$ sends $E$ to $(\sigma^{-1})^*E$. It is well-known that $LS_{\tilde X}(\check{\GG})^{\Sigma}\,\iso\, LS_X(\check{\GG})$. 
 
\sssec{} Consider the new action of $\Sigma$ on $LS_{\tilde X}(\check{\GG})$ such that $\sigma\in\Sigma$ sends $E$ to $(\sigma^{-1})^*(E^{\sigma})$. By definition, $LS^{\pi,\rho}_X(\check{\GG})$ is the category of $\Sigma$-invariants of $LS_{\tilde X}(\check{\GG})$ with respect to this action. 

 Write $LS_S$ for the category of smooth sheaves on a base $S$, we view it as a symmetric monoidal abelian category. Note that $LS_{B(\Sigma)}\,\iso\,\Rep(\Sigma)$ naturally. The $\Sigma$-torsor $\pi: \tilde X\to X$ gives a morphism $X\to B(\Sigma)$. The pullback along this map defines a functor $\Rep(\Sigma)\to LS_X$. A choice of $x\in X$ gives the natural identification $\Rep(\pi_1(X, x))\,\iso\, LS_X$. A choice of $\tilde x\in \tilde X$ with $x=\pi(\tilde x)$ defines a homomorphism $\mu_{\tilde x}: \pi_1(X, x)\to \Sigma$ such that 
$$
\Rep(\Sigma)\to LS_X\,\iso\, \Rep(\pi_1(X, x))
$$ 
is the restriction along $\mu_{\tilde x}$. 
  
\begin{Rem} i) The category $LS^{\pi,\rho}_X(\check{\GG})$ is naturally equivalent to the category of exact symmetric monoidal $\Qlb$-linear functors $f$
that fit into the commutative diagram
$$
\begin{array}{ccc}
\Rep(\check{\GG})^{\Sigma} & \toup{f} & LS_X\\
\uparrow & \nearrow\\
\Rep(\Sigma)
\end{array}
$$
ii) For a choice of $\tilde x\in \tilde X$ with $x:=\pi(\tilde x)$, this identifies
$LS^{\pi,\rho}_X(\check{\GG})$ with the category of continuous homomorphisms $\bar f: \pi_1(X, x)\to \check{\GG}\rtimes\Sigma$ making the following diagram commute
$$
\begin{array}{ccc}
\pi_1(X, x) & \toup{\bar f} & \check{\GG}\rtimes\Sigma\\
\downarrow\lefteqn{\scriptstyle \mu_{\tilde x}} & \swarrow\\ 
\Sigma
\end{array}
$$
\end{Rem}
\begin{proof} Let $(E, \{\beta_{\sigma}\})$ be a $\pi$-twisted $\check{\GG}$-local system on $\tilde X$. Define the corresponding functor $f$ as follows. Let us be given an object of $\Rep(\check{\GG})^{\Sigma}$ given by $V\in\Rep(\check{\GG})$ and a collection of isomorphisms $c_{\sigma}: \sigma^*V\,\iso\, V$ for $\sigma\in\Sigma$. To these data we attach the following smooth sheaf $f(V, c)$ on $X$.

 For $\sigma\in\Sigma$ our data yield isomorphisms 
$$
\sigma^*(V_E)\,\iso\, V_{\sigma^*E}\toup{\beta_{\sigma}} V_{E^{\sigma}}\,\iso\, (\sigma^*V)_E\toup{c_{\sigma}} V_E
$$
of local systems on $\tilde X$, where the unnamed maps are canonical isomorphisms. These are descent data for $V_E$ giving rise to a smooth sheaf $f(V, c)$ on $X$. One checks that this is an equivalence of categories.
\end{proof}  

\sssec{More Hecke functors} Recall that an object of $(\Sph_{\GG})^{\Sigma}$ is a collection $(\cS, \{c_{\sigma}\}_{\sigma\in\Sigma})$, where $c_{\sigma}: \sigma^*\cS\,\iso\,\cS$ is such that $c_{\tau}\comp\tau^*c_{\sigma}=c_{\sigma\tau}$ for $\sigma, \tau\in\Sigma$. For $\sigma\in\Sigma$, $\cS\in\Sph_{\GG}$ one has naturally 
$$
\sigma^*(\ast\cS)\,\iso\, \ast(\sigma^*\cS)
$$ 
So, $(\Sph_{\GG})^{\Sigma}$ inherits an auto-equivalence still denoted 
$\ast: (\Sph_{\GG})^{\Sigma}\,\iso\,(\Sph_{\GG})^{\Sigma}$ by abuse of notations. The latter sends a collection $(\cS, \{c_{\sigma}\})$ to the collection $(\ast\cS, \{c'_{\sigma}\})$, where $c'_{\sigma}$ is the composition $\sigma^*(\ast\cS)\,\iso\, \ast(\sigma^*\cS)\toup{\ast(c_{\sigma})} \ast\cS$. 

 Let $\cH_G$ be the stack classifying $x\in X$, $\cF_G,\cF'_G\in\Bun_G$ and an isomorphism $\cF_G\,\iso\,\cF'_G\mid_{X-x}$. Let $\supp: \cH_G\to X$ be the pojection sending the above point to $x$. As in Section~\ref{Sect_2.3.1}, one defines a diagram
$$
\Bun_G \getsup{h^{\la}} \cH_G\toup{h^{\ra}} \Bun_G
$$ 
 
 Let $_{\tilde X}\Bun_G$ be the stack classifying $\cF_G\in\Bun_G, \tilde x\in \tilde X$ and an isomorphism of $G$-torsors $\beta: \cF_G\,\iso\, \cF^0_G\mid_{D_{\pi(\tilde x)}}$. Then $_{\tilde X}\Bun_G$ can equivalently be seen as the stack classifying $\tilde x\in \tilde X, \cF_G\in\Bun_G$ and an isomorphism $\tilde\beta: \pi^*\cF_G\,\iso\, \cF^0_{\GG}\mid_{D_{\tilde x}}$ of $\GG$-torsors on $D_{\tilde x}$. Indeed, $\pi$ gives an isomorphism $D_{\tilde x}\,\iso\, D_x$. 
 
 Let $G_{\tilde X}$ be the group scheme over $\tilde X$ with fibre $\GG(\cO_{\tilde x})$ at $\tilde x$. Let $\Gr_{\GG, \tilde X}$ be the ind-scheme over $\tilde X$ whose fibre at $\tilde x$ is $\Gr_{\GG,\tilde x}$. 
We have the isomorphisms
$$
\tilde\id^l, \tilde\id^r: \tilde X\times_X \cH_G\,\iso\, {_{\tilde X}\Bun_G\times^{G_{\tilde X}} \Gr_{\GG, \tilde X}}
$$
such that the projection to the first factor corresponds to $h^{\la}, h^{\ra}$ respectively (the product in the RHS is taken over $\tilde X$). 
 
So, for $\cS\in \Sph_{\GG}$, we may form the corresponding twisted exteriour products 
\begin{equation}
\label{ext_product_for_twisted_Hecke_functor}
(\IC(\tilde X\times\Bun_G)\tboxtimes\cS)^l, 
(\IC(\tilde X\times\Bun_G)\tboxtimes\cS)^r\in \D(\tilde X\times_X \cH_G)
\end{equation} 
normalized to be perverse. Assume in addition that $\cS\in (\Sph_{\GG})^{\Sigma}$. Then (\ref{ext_product_for_twisted_Hecke_functor}) are naturally equipped with the descent data for $\pi\times\id: \tilde X\times_X \cH_G\to \cH_G$, we denote by $\IC(\cH_G, V)^l, \IC(\cH_G, V)^r$ the perverse sheaves on $\cH_G$ so obtained. 

 Define the functors 
\begin{equation}
\label{functors_Hecke_bold}
\bH_G^{\la}, \bH_G^{\ra}: (\Sph_{\GG})^{\Sigma}\times \D(\Bun_G)\to \D(X\times\Bun_G)
\end{equation}
by 
$$
\bH^{\la}_G(\cS, K)=(\supp\times h^{\la})_!(\IC(\cH_G, \ast V)^r\otimes (h^{\ra})^*K)[-\dim\Bun_G]
$$
and 
$$
\bH^{\ra}_G(\cS, K)=(\supp\times h^{\ra})_!(\IC(\cH_G, V)^l\otimes (h^{\la})^*K)[-\dim\Bun_G]
$$

 For $\cS\in \Sph_{\GG}^{\Sigma}$ and $\pi\times\id: \tilde X\times\Bun_G\to X\times\Bun_G$ we get 
$$
(\pi\times\id)^*\bH^{\ra}_G(\cS, K)\,\iso\, \H^{\la}_G(\Res(\cS), K),
$$ 
where $\Res: \Sph_{\GG}^{\Sigma}\to \Sph_{\GG}$ is the restriction functor along $\check{\GG}\to \check{\GG}\rtimes\Sigma$.

 Recall the $\Sigma$-equivariant equivalence (\ref{Loc_gr_Sigma_equivariant}). It allows to extend $\bH_G^{\la}, \bH_G^{\ra}$ to functors still denoted
\begin{equation}
\label{functors_Hecke_bold_ext}
\bH_G^{\la}, \bH_G^{\ra}: (\D\Sph_{\GG})^{\Sigma}\times \D(\Bun_G)\to \D(X\times\Bun_G)
\end{equation}
and given by the same formula. 

\sssec{Functoriality} 
\label{Sect_2.3.14_functoriality}
The setting for the geometric Langlands functoriality in the non-twisted setting is proposed in (\cite{L3}, Section~2). In this section we explain the modification we need for our twisted setting.

Assume in addition we have a split connected reductive groups $\HH$ with a maximal torus and Borel $\TT_{\HH}\subset \BB_{\HH}\subset \HH$. Assume given an action of $\Sigma$ on $\HH$. Let $H=(\HH\times\tilde X)/\Sigma$. 

 As above, we get a $\Sigma$-action on $\check{\HH}$. Assume given a homomorphism $\bar\kappa: \check{\HH}\rtimes\Sigma\to \check{\GG}\rtimes\Sigma$ preserving the projections to $\Sigma$. Let $\kappa: \check{\HH}\to\check{\GG}$ be its restriction. Consider a complex $\cM\in \D(\Bun_G\times\Bun_H)$ giving rise to the functor 
$$
F_G: \D(\Bun_H)\to \D(\Bun_G)
$$ 
as in (\cite{L3}, Section~2). Namely, for the diagram of projections
$$
\Bun_H\getsup{p_H} \Bun_G\times\Bun_H \toup{p_G} \Bun_G
$$
we let 
$$
F_G(K)=(p_G)_!((p_H)^*K\otimes \cM)[-\dim\Bun_H]
$$ 
For a scheme $S$, we similarly get the functor $\id\boxtimes F_G: \D(S\times\Bun_H)\to \D(S\times \Bun_G)$ with the kernel $\pr^*\cM$ for the projection $\pr: S\times \Bun_H\times\Bun_G\to \Bun_H\times\Bun_G$. 
 
 A functoriality datum for $F_G$ is a collection of isomorphisms
\begin{equation}
\label{iso_epsilon_V_Sect_2.3.14}
\varepsilon_V: \H^{\la}_G(V, F_G(K))\,\iso\, (\id\boxtimes F_G)\H^{\la}_H(\Res^{\kappa}(V), K)
\end{equation}
in $\D(\tilde X\times \Bun_G)$ functorial in $V\in \Rep(\check{\GG})$, $K\in \D(\Bun_H)$ and compatible with the tensor structure on $\Rep(\check{\GG})$ as in (\cite{L3}, Section~2). 

  In the twisted setting we impose additional properties on the above isomorphisms. First, consider the simplifying assumption that $\kappa$ is $\Sigma$-equivariant. In this case we require in addition that for $\sigma\in\Sigma$ the diagram commutes
\begin{equation}
\label{diag_functoriality_Sigma-compatibility_2.3.8}
\begin{array}{ccc}
(\sigma\times\id)^*\H^{\la}_G(V, F_G(K))& \toup{(\sigma\times\id)^*\varepsilon_V} & (\id\boxtimes F_G)(\sigma\times\id)^*\H^{\la}_H(\Res^{\kappa}(V), K)\\
&& \downarrow\\
\downarrow && (\id\boxtimes F_G)\H^{\la}_H(\sigma^*\Res^{\kappa}(V), K)\\
&& \mid\mid\\
\H^{\la}_G(\sigma^*V, F_G(K)) & \toup{\varepsilon_{\sigma^*V}} & (\id\boxtimes F_G)\H^{\la}_H(\Res^{\kappa}(\sigma^*V), K),
\end{array}
\end{equation}
where the vertical arrows are the isomorphisms (\ref{iso_key_for_twisted_setting}).     

\sssec{} 
\label{Section_2.3.15}
Consider now the case of $\bar\kappa$ such that $\kappa: \check{\HH}\to\check{\GG}$ is not necessarily $\Sigma$-equivariant. Then restriction along $\bar\kappa$ gives a functor 
$$
\Res^{\bar\kappa}: \Rep(\check{\GG})^{\Sigma}\, \iso\,\Rep(\check{\GG}\rtimes\Sigma) \to \Rep(\check{\HH}\rtimes\Sigma)\,\iso\,\Rep(\check{\HH})^{\Sigma}
$$ 

 In this case a functoriality datum for $F_G$ is a collection of isomorphisms (\ref{iso_epsilon_V_Sect_2.3.14}) compatible with the tensor structure on $\Rep(\check{\GG})$ and a collection of isomorphisms 
\begin{equation}
\label{iso_for_Sect2.3.15_functoriality_varepsilon}
{\boldsymbol\varepsilon}_V: \bH^{\la}_G(V, F_G(K))\,\iso\, (\id\boxtimes F_G)\bH^{\la}_H(\Res^{\bar\kappa}(V), K)
\end{equation}
in $\D(\tilde X\times \Bun_G)$ functorial in $V\in \Rep(\check{\GG})^{\Sigma}, K\in \D(\Bun_H)$ and compatible with the tensor structure on 
$$
\Rep(\check{\GG})^{\Sigma}\,\iso\,\Rep(\check{\GG}\rtimes\Sigma)
$$ 
It is required in addition that for $V\in \Rep(\check{\GG})^{\Sigma}$ one has $(\pi\times\id)^*{\boldsymbol\varepsilon}_V\,\iso\, \varepsilon_V$.  

\sssec{} 
\label{Sect_2.3.16}
We need even more general twisted setting. Namely, assume given a homomorphism $\bar\kappa: (\check{\HH}\rtimes\Sigma)\times\Gm\to \check{\GG}\rtimes\Sigma$ for which we denote by $\kappa: \check{\HH}\times\Gm\to \check{\GG}$ its restriction. (It is expected that $\kappa$ always comes from a homomorphism $\check{\HH}\times\SL_2\to \check{\GG}$). Note that the restriction of $\bar\kappa$ to $\Gm$ factor is a map $\Gm\to \check{\GG}\hook{} \check{\GG}\rtimes\Sigma$. The restriction along $\bar\kappa$ gives a
functor 
$$
\Res^{\bar\kappa}: \Rep(\check{\GG})^{\Sigma}\,\iso\,\Rep(\check{\GG}\rtimes\Sigma)\to \Rep((\check{\HH}\rtimes \Sigma)\times\Gm)\,\iso\, \Rep(\check{\HH}\times\Gm)^{\Sigma},
$$
recall that $\Sigma$ acts trivially on $\Gm$ here. 

 In this case the functoriality datum for $F_G$ is a collection of isomorphisms (\ref{iso_epsilon_V_Sect_2.3.14}) compatible with the tensor structure on $\Rep(\check{\GG})$ and a collection of isomorphisms (\ref{iso_for_Sect2.3.15_functoriality_varepsilon}) in $\D(\tilde X\times \Bun_G)$ functorial in $V\in \Rep(\check{\GG})^{\Sigma}, K\in \D(\Bun_H)$ and compatible with the tensor structure on 
$$
\Rep(\check{\GG})^{\Sigma}\,\iso\,\Rep(\check{\GG}\rtimes\Sigma)
$$ 
It is required in addition that for $V\in \Rep(\check{\GG})^{\Sigma}$ one has $(\pi\times\id)^*{\boldsymbol\varepsilon}_V\,\iso\, \varepsilon_V$.  

 We used here the extension (\ref{functors_Hecke_bold_ext}) of the functors (\ref{functors_Hecke_bold}). 
 
\sssec{Example} Let $\HH$ be trivial and $\kappa:  \Gm\to \check{\GG}$ be equal to $2\check{\rho}_{\GG}$. Then $\kappa$ is $\Sigma$-invariant, because $\Sigma$ permutes the simple roots of $\check{\GG}$. So, it gives rise to a map $\bar\kappa: \Sigma\times\Gm\to \check{\GG}\rtimes\Sigma$. We get $\Bun_H=\Spec k$. The functor $F_G: \D(\Spec k)\to \D(\Bun_G)$ sending $K$ to $K\otimes\IC(\Bun_G)$ has a functoriality datum for $\bar\kappa$ in the sense of Section~\ref{Sect_2.3.16}. 

\sssec{Question} Assume both $\GG,\HH$ are tori, and we are given $\bar\kappa: \check{\HH}\rtimes\Sigma\to \check{\GG}\rtimes\Sigma$ as in Section~\ref{Sect_2.3.14_functoriality}. What is the corresponding functor $F_G: \D(\Bun_H)\to\D(\Bun_G)$?

\ssec{Theta-lifting functors}  
\label{Sect_Theta-lifting functors}

\sssec{} 
\label{Sect_2.4.1_now}
The following notations are borrowed from \cite{L}. Write $\Omega$ for the canonical line bundle on $X$ (everywhere except Sections~\ref{Sect_Local_theory},\ref{Sect_Dual_pair_GSp_GO}, where we work in the local setting).

For $r\ge 1$ let $G_r$ denote the sheaf of automorphisms of $\cO_X^r\oplus\Omega^r$ preserving the natural symplectic form $\wedge^2(
\cO_X^r\oplus\Omega^r)\to\Omega$. The stack $\Bun_{G_r}$ of $G_r$-torsors on $X$ classifies $M\in\Bun_{2r}$ equipped with a symplectic form $\wedge^2 M\to\Omega$. Write $\cA_{G_r}$ for the line bundle on $\Bun_{G_r}$ with fibre $\det\RG(X, M)$ at $M$, we view it as $\ZZ/2\ZZ$-graded of parity zero. Let $\Bunt_{G_r}\to\Bun_{G_r}$ denote the $\mu_2$-gerbe of square roots of $\cA_{G_r}$. Write $\Aut$ for the perverse theta-sheaf on $\Bunt_{G_r}$ (cf. \cite{L1}).  

\sssec{}  
\label{Sect_2.3.2_def_tilde_pi}
Pick an  \'etale degree 2 covering $\pi:\tilde X\to X$ with Galois group $\Sigma=\{\id,\sigma\}$. Let $\cE$ be the $\sigma$-anti-invariants in $\pi_*\cO$, it is equipped with a trivialization $\cE^2\,\iso\, \cO_X$. 
  
\sssec{}  
\label{Sect_2.3.3}
Let $n, m\in\NN$ and $\GG=\GSp_{2n}$. We realize $\GG$ as the subgroup of $\GL_{2n}$ preserving up to a scalar the symplectic form given by the matrix 
$$
\left(
\begin{array}{cc}
0  & E_n\\
-E_n & 0
\end{array}
\right),
$$ 
where $E_n\in\GL_n$ is the unity. Pick the maximal torus $\TT_{\GG}$ of diagonal matrices, and the Borel subgroup $\BB_{\GG}$ preserving for $i=1,\ldots, n$ the isotropic subspace generated by the first $i$ vectors $\{e_1,\ldots, e_i\}$.

 The stack $\Bun_{\GG}$ classifies $M\in \Bun_{2n}, \cA\in\Bun_1$ with symplectic form $\wedge^2 M\to\cA$. Write $\cA_{\GG}$ for the $\ZZ/2\ZZ$-graded line bundle on $\Bun_{\GG}$ with fibre $\det\RG(X, M)$ at $(M,\cA)$. 

 We let $\Sigma$ act trivially on $\GG$, so the corresponding twisted group is $G=\GG$. For $\cS\in\Sph_{\GG}, K\in \D(\Bun_G)$ by definition
$$
\H^{\la}_G(\cS, K)=(\pi\times\id)^*\H^{\la}_{\GG}(\cS, K)
$$
with $\pi\times\id: \tilde X\times\Bun_G\to X\times\Bun_G$ the projection.
  
  Write $\check{\omega}_0$ for the character of $\GG$ such that $\GG$ acts on the determinant of the standard representation as $n\check{\omega}_0$, this is the usual similitude character. It yields an isomorphism $\pi_1(\GG)\to\ZZ$. Write $_a\Sph_{\GG}\subset\Sph_{\GG}$ for the full subcategory of objects that vanish off the connected component $\Gr_{\GG}^{\theta}$ satisfying $\<\theta,\check{\omega}_0\>=-a$. Let $_a\Rep(\check{\GG})\subset \Rep(\check{\GG})$ be the preimage of $_a\Sph_{\GG}$ under $\Loc$. An irreducible representation of $\check{\GG}$ of highest weight $\lambda$ appear in $_a\Rep(\check{\GG})$ iff $\<\lambda, \check{\omega}_0\>=-a$.  
  
\sssec{} 
\label{Sect_2.3.4}
Let $\HH=\GSO_{2m}=(\Gm\times\SO_{2m})/(-1,-1)$. Realize $\HH$ as the subgroup of $\GL_{2m}$ preserving up to a multiple the symmetric form given by the matrix
$$
\left(
\begin{array}{cc}
0  & E_m\\
E_m & 0
\end{array}
\right),
$$ 
where $E_m\in\GL_m$ is the unity. Take $\TT_H$ to be the maximal torus of diagonal matrices, $\BB_{H}$ the Borel subgroup preserving for $i=1,\ldots,m$ the isotropic subspace generated by the first $i$ base vectors $\{e_1,\ldots,e_i\}$. 

 Pick an involution $\tilde\sigma$ of $\HH$ as follows. For $m=1$ we identify $\HH\,\iso\,\Gm\times\Gm$ in such a way that $\tilde\sigma$ permutes the two copies of $\Gm$. For $m\ge 2$ we take $\tilde\sigma_0\in \OO_{2m}(k)$ interchanging $e_m$ and $e_{2m}$ and acting trivially on the orthogonal complement to $\{e_m, e_{2m}\}$. Let $\tilde\sigma(h)=\tilde\sigma_0 h\tilde\sigma_0^{-1}$ for $h\in \HH$. Let $\Sigma$ act on $\HH$ so that $\sigma$ sends $h$ to $\tilde\sigma(h)$. Then $\tilde\sigma$ preserves $\TT_{\HH}$ and $\BB_{\HH}$. For $m\ge 2$ it induces the unique\footnote{except for $m=4$. The group $\GSO_8$ also has other outer forms, we do not consider them. We consider only the automorphisms coming from the semi-direct product $\GO_{2m}\,\iso\, \GSO_{2m}\ltimes\ZZ/2\ZZ$.} nontrivial automorphism of the Dynkin diagram of $\HH$.  

 For $m\ge 1$ we get an action of $\Sigma=\{1,\sigma\}$ on $\check{\HH}$ as in Section~\ref{Sect_for_Out}. We denote by $\sigma_{\HH}\in\Aut(\check{\HH})$ the corresponding involution of $\check{\HH}$. 

Let $H$ be the group scheme on $X$, the twisting of $\HH$ by the $\Sigma$-torsor $\pi:\tilde X\to X$. Let $T_{H}\subset B_{H}$ be the corresponding twists of $\TT_{\HH}$, $\BB_{\HH}$. Note that if $m=1$ then $H\,\iso\,\pi_*\Gm$. 
 
 Write $\check{\alpha}_0$ for the character of $\HH$ such that $\HH$ acts on the determinant of the standard representation by $m\check{\alpha}_0$, this is the usual similitude character. Write $_a\Sph_{\HH}\subset\Sph_{\HH}$ for the full subcategory of objects that vanish off the connected components $\Gr_{\HH}^{\theta}$ of $\Gr_{\HH}$ satisfying $\<\theta, \check{\alpha}_0\>=-a$. Denote by $_a\Rep(\check{\HH})$ the preimage of $_a\Sph_{\HH}$ under $\Loc$. An irreducible $\check{\HH}$-module of highest weight $\lambda$ appear in $_a\Rep(\check{\HH})$ iff $\<\lambda, \check{\alpha}_0\>=-a$. 

\sssec{} Let $\Bun_{\mu_2}$ denote the stack classifying $\mu_2$-torsors on $X$. Its connected components are indexed by 
$\H^1(X, \mu_2)$, each connected component is isomorphic to the classifying stack $B(\mu_2)$. 

 Set $\bar\HH=\GO_{2m}=\Gm\times\OO_{2m}/(-1, -1)$. We have an exact sequence $1\to \HH\to \bar\HH\to \mu_2\to 1$, this is a semi-direct product $\bar\HH\,\iso\, \HH\rtimes \mu_2$. The stack $\Bun_{\bar\HH}$ classifies $V\in\Bun_{2m}, \cC\in\Bun_1$ and a non-degenerate symmetric form $\Sym^2 V\to \cC$. 
 
  A point of $\Bun_{\bar\HH}$ yields a $\mu_2$-torsor on $X$ given by the line bundle $\cC^{-m}\otimes\det V$ on $X$ with the induced trivialization $(\cC^{-m}\otimes\det V)^2\,\iso\,\cO$. Let $\Bun_{\bar \HH}^{\pi}$ be the substack of $\Bun_{\bar\HH}$ given by requiring that the above $\mu_2$-torsor equals $\pi$ in $\H^1(X, \mu_2)$. 
\begin{Rem}  
\label{Rem_2.3.6}
As for any semi-direct product of groups, the stack $\Bun_{\bar \HH}$ classifies a $\mu_2$-torsor $\cF$ and a $\HH_{\cF}$-torsor. Here $\HH_{\cF}=(\HH\times \cF)/\mu_2$, where $\mu_2$ acts diagonally.   
\end{Rem}  
  
\sssec{} By Remark~\ref{Rem_2.3.6}
$$
\Bun_H\,\iso\, \Spec k\times_{\Bun_{\mu_2}} \Bun_{\bar\HH},
$$
where the map $\Spec k\to \Bun_{\mu_2}$ is given by the $\mu_2$-torsor $\pi: \tilde X\to X$. The projection $\Bun_{H}\to\Bun_{\bar\HH}^{\pi}$ is a $\mu_2$-torsor. 

Note that the stack $\Bun_H$ classifies: $V\in\Bun_{2m},  \, \cC\in\Bun_1$, a nondegenerate symmetric
form $\Sym^2 V\to\cC$, and a compatible trivialization $\gamma: \cC^{-m}\otimes\det V\,\iso\,\cE$. This means that the composition 
$$
\cC^{-2m}\otimes(\det V)^2\,\toup{\gamma^2}\,\cE^2\,\iso\,\cO
$$ 
is the isomorphism induced by $V\,\iso\, V^*\otimes\cC$.

Let $\cA_H$ be the $\ZZ/2\ZZ$-graded line bundle on $\Bun_H$ with fibre $\det\RG(X,V)$ at $(V,\cC)$. Set
$$
\Bun_{G,H}=\Bun_H\times_{\Pic X} \Bun_G, 
$$
where  the map $\Bun_H\to \Pic X$ sends $(V,\cC, \Sym^2 V\to\cC)$ to $\Omega\otimes\cC^{-1}$, and $\Bun_G\to\Pic X$ sends $(M, \wedge^2 M\to\cA)$ to $\cA$. So, we have an isomorphism
$\cC\otimes \cA\,\iso\, \Omega$ for a point of $\Bun_{G,H}$. Let
$$
\tau: \Bun_{G,H}\to \Bun_{G_{2nm}}
$$ 
be the map sending a point as above to $V\otimes M$ with the induced symplectic form $\wedge^2(V\otimes M)\to\Omega$. 

 By (\cite{L}, Proposition~2), for a point $(M,\cA, V, \cC)$ of $\Bun_{G,H}$
there is a canonical $\ZZ/2\ZZ$-graded isomorphism 
\begin{equation}
\label{iso_line_bundle_Bun_GH}
\det\RG(X, V\otimes M)\,\iso\, \frac{\det\RG(X, V)^{2n}\otimes \det\RG(X, M)^{2m}}{\det\RG(X, \cO)^{2nm}\otimes\det\RG(X,\cA)^{2nm}}
\end{equation}  
It yields a map $\tilde\tau: \Bun_{G,H}\to \Bunt_{G_{2nm}}$ sending $(\wedge^2 M\to\cA, \Sym^2 V\to\cC, \cA\otimes\cC\,\iso\,\Omega)$ to $(\wedge^2(M\otimes V)\to\Omega, \cB)$. Here 
$$
\cB=\frac{\det\RG(X,V)^n\otimes\det\RG(X, M)^m}{\det\RG(X,\cO)^{nm}\otimes \det\RG(X,\cA)^{nm}},
$$
and $\cB^2$ is identified with $\det\RG(X, M\otimes V)$ via (\ref{iso_line_bundle_Bun_GH}).  
  
\begin{Def} Set $\Aut_{G,H}=\tilde\tau^*\Aut[\dimrel(\tau)]$. 
For the diagram of projections
$$
\Bun_H\getsup{\gq}  \Bun_{G,H}\toup{\gp} \Bun_G
$$
define $F_G: \D(\Bun_H)\to\D(\Bun_G)$ and $F_H: \D(\Bun_G)\to\D(\Bun_H)$ by
$$
\begin{array}{c}
F_G(K)=\gp_!(\Aut_{G,H}\otimes \gq^*K)[-\dim\Bun_H]  \\ \\
F_H(K)=\gq_!(\Aut_{G,H}\otimes \gp^*K)[-\dim\Bun_G]
\end{array}
$$
\end{Def}  
  
\ssec{Morphism of dual groups} 

\sssec{} The Langlands dual groups $\check{\GG}, \check{\HH}$ are described in Section~\ref{Sect_Root data}. Consider the split group $\Spin_{2m}$ over $\Spec k$. For $m\ge 2$ let $i_{\HH}\in \Spin_{2m}$ be the central element of order 2 such that $\Spin_{2m}/\{i_{\HH}\}\,\iso\,\SO_{2m}$. For $m\ge 2$ set
$$
\GSpin_{2m}=\Gm\times\Spin_{2m}/\{(-1, i_{\HH})\}
$$
We declare $\GSpin_2\,\iso\, \Gm\times\Gm$. Then $\check{\HH}\,\iso\,\GSpin_{2m}$. We also have $\check{\GG}\,\iso\,\GSpin_{2n+1}$, where 
$$
\GSpin_{2n+1}=\Gm\times\Spin_{2n+1}/\{(-1, i_{\GG})\},
$$ 
and $i_{\GG}\in\Spin_{2n+1}$ is the nontrivial central element. Let $V_{\HH}$ (resp., $V_{\GG}$) denote the standard representation of $\SO_{2m}$ (resp., of $\SO_{2n+1}$). The image of $\check{\omega}_0: \Gm\to \check{\GG}$ is the center of $\check{\GG}$. For $m\ge 2$ the image of $\check{\alpha}_0: \Gm\to \check{\HH}$ is the connected component of the center of $\check{\HH}$.  

\sssec{{\scshape CASE} $m\le n$} 
\label{Sect_2.4.2_case_m_le_n}
Define the homomorphism $\bar\kappa: (\check{\HH}\rtimes\Sigma)\times\Gm\to \check{\GG}\times\Sigma$ as follows.
Pick an inclusion $V_{\HH}\hook{} V_{\GG}$ compatible with symmetric forms. It gives rise to the inclusion $\Spin(V_{\HH})\to \Spin(V_{\GG})$. Then the map 
$$
\Gm\times \Spin(V_{\HH})\to \Gm\times\Spin(V_{\GG}), \; (z,y)\mapsto (z^{-1}, y)
$$
descends to an inclusion $i_{\kappa}: \check{\HH}\to \check{\GG}$.
Our normalization is such that $i_{\kappa}\check{\alpha}_0(z)=\check{\omega}_0(z^{-1})$ for $z\in\Gm$. 

%Pick $\sigma_{\GG}\in \SO(V_{\GG})\,\iso\, \check{\GG}_{ad}$ 
%normalizing $\check{\TT}_{\GG}$ and preserving 
%$V_{\HH}$ and $\check{\TT}_{\HH}\subset \check{\BB}_{\HH}$. Let $\sigma_{\HH}\in\OO(V_{\HH})$ be its restriction to $V_{\HH}$. Then $\sigma_{\HH}$ acts on $\check{\HH}$ by conjugation preserving $\check{\TT}_{\HH}, \check{\BB}_{\HH}$. Let $\sigma\in \Sigma$ act on $\check{\HH}$ by $\sigma_{\HH}$. 

 Concretely, we take $V_{\GG}=\Qlb^{2n+1}$ with the symmetric form given by the matrix
$$
\left(
\begin{array}{ccc}
0 & E_n & 0\\
E_n & 0 & 0\\
0 & 0 & 1
\end{array}
\right),
$$
where $E_n\in\GL_n$ is the unity. We assume $\check{\TT}_{\GG}$ is the preimage of the torus of diagonal matrices under $\check{\GG}\to\SO_{2n+1}$. Let $V_{\HH}\subset V_{\GG}$ be generated by $\{e_1,\ldots, e_m, e_{n+1},\ldots,e_{n+m}\}$. We assume $\check{\TT}_{\HH}$ is the preimage under $\check{\HH}\to\SO(V_{\HH})$ of the torus of diagonal matrices, and $\check{\BB}_{\HH}$ the Borel subgroup preserving for $i=1,\ldots,m$ the isotropic subspace generated by $\{e_1,\ldots,e_i\}$. 

 Consider $\sigma_0\in \SO(V_{\GG})$ permuting $e_m$ and $e_{n+m}$, sending $e_{2n+1}$ to $-e_{2n+1}$ and acting trivially on the other base vectors. Let $\sigma_{0,\GG}\in \check{\GG}$ be an element of order two over $\sigma_0$. 
 
 We define $\Sigma\to \check{\GG}\times\Sigma$ by sending $\sigma$ to $(\sigma_{0, \GG}, \sigma)$. The adjoint action of $\sigma_{0, \GG}$ on $\check{\GG}$ preserves the above subgroup $\check{\HH}$. We may and do assume that the so obtained involution of $\check{\HH}$ coincides with $\sigma_{\HH}$ defined in Section~\ref{Sect_2.3.4}. This defines the desired map $\check{\HH}\rtimes\Sigma\to \check{\GG}\times\Sigma$ that fits into the framework of Section~\ref{Section_2.3.15} (the corresponding map $\check{\HH}\to \check{\GG}$ is not $\Sigma$-equivariant). 
 
  Let $\cC(\check{\HH})\subset \check{\GG}$ be the connected centralizer of $\check{\HH}$ in $\check{\GG}$. We define the component $\Gm\to \check{\GG}$ of $\bar\kappa$ as the composition 
$$
\Gm\hook{}\SL_2\toup{prin}\cC(\check{\HH}),
$$ 
where $prin$ is the homomorphism corresponding to the principal nilpotent in $\Lie \cC(\check{\HH})$, and $\Gm$ is the maximal torus of diagonal matrices in $\SL_2$. This completes the definition of $\bar\kappa$. 

 According to our convention of Section~\ref{Sect_2.3.14_functoriality}, we denote by $\kappa: \check{\HH}\times\Gm\to \check{\GG}$ the restriction of $\bar\kappa$. 
   
\sssec{{\scshape CASE} $m>n$} 
\label{Sect_2.4.3_m_bigger_n}
Define the homomorphism $\bar\kappa: \check{\GG}\times\Sigma\times\Gm\to \check{\HH}\rtimes\Sigma$ as follows.
Pick an inclusion $V_{\GG}\hook{} V_{\HH}$ compatible with symmetric forms. It yields an inclusion $\Spin(V_{\GG})\to \Spin(V_{\HH})$. Then the map
$$
\Gm\times \Spin(V_{\GG})\to \Gm\times\Spin(V_{\HH}),\; (z,y)\mapsto (z^{-1}, y)
$$
descends to an inclusion $i_{\kappa}: \check{\GG}\hook{}\check{\HH}$ compatible with the corresponding maximal tori. Our normalization is such that $i_{\kappa}\check{\omega}_0(z)=\check{\alpha}_0(z^{-1})$ for $z\in\Gm$.

Take the symmetric form on $V_{\HH}=\Qlb^{2m}$ given by the matrix
$$
\left(
\begin{array}{cc}
0  & E_m\\
E_m & 0
\end{array}
\right)
$$
Let $V_{\GG}$ be the subspace of $V_{\HH}$ generated by 
$
\{e_1,\ldots, e_n; e_{m+1},\ldots, e_{m+n}; e_{n+1}+e_{m+n+1}\}
$. 
We assume that $\check{\TT}_{\HH}$ is the preimage under $\check{\HH}\to\SO(V_{\HH})$ of the torus of diagonal matrices, and $\check{\BB}_{\HH}$ is the Borel subgroup preserving for $i=1,\ldots,m$ the isotropic subspace of $V_{\HH}$ generated by $\{e_1,\ldots, e_i\}$. We may assume $\check{\TT}_{\GG}$ is the preimage under $\check{\GG}\to\check{\HH}$ of $\check{\TT}_{\HH}$. Recall the involution $\sigma_{\HH}\in\Aut(\check{\HH})$ defined in Section~\ref{Sect_2.3.4}. It preserves and acts trivially on 
 $\check{\GG}$. So, the map $i_{\kappa}$ is $\Sigma$-equivariant. We require the restriction of $\bar\kappa$ to $\Sigma$ to be the map $\Sigma\to \check{\HH}\rtimes\Sigma$, $\sigma\mapsto (1, \sigma)$. The so obtained map $\check{\GG}\times\Sigma\to \check{\HH}\rtimes\Sigma$ fits in the framework of Section~\ref{Sect_2.3.14_functoriality}. 

 Let $\cC(\check{\GG})\subset \check{\HH}$ be the connected centralizer of $\check{\GG})$ in $\check{\HH}$. We define the component $\Gm\to \check{\HH}$ of $\bar\kappa$ as the composition
$$
\Gm\hook{}\SL_2\toup{prin} \cC(\check{\GG}),
$$
analogous to Section~\ref{Sect_2.4.2_case_m_le_n}. In particular, $prin$ is the homomorphism corresponding to the principal nilopotent in $\Lie \cC(\check{\GG})$. This yields the desired map $\bar\kappa$. 

 According to the above conventions, we denote by $\kappa: \check{\GG}\times\Gm\to \check{\HH}$ the restriction of $\bar\kappa$. 
   
\sssec{} Recall the definition of the functoriality datum from Section~\ref{Sect_2.3.16}. The following is our main result. 
  
\begin{Th}
\label{Th_theta-lifting_functors-general}
1) For $m\le n$ 
% the map $i_{\kappa}$ extends to a homomorphism $\kappa=(i_{\kappa},\delta_{\kappa}): \check{\HH}\times\Gm\to\check{\GG}$ with the following property. 
the functor $F_G$ is naturally equipped with a functoriality datum for $\bar\kappa:(\check{\HH}\rtimes\Sigma)\times\Gm\to \check{\GG}\times\Sigma$. So, there exists an isomorphism
\begin{equation}
\label{iso_theta_lifting_Hecke_Th1}
\H^{\la}_{G}(V, F_G(K))\,\iso\, (\id\boxtimes F_G)(\H^{\la}_{H}(\Res^{\bar\kappa}(V), K))
\end{equation}
in $\D(\tilde X\times \Bun_G)$ functorial in $V\in\Rep(\check{\GG})$ and $K\in \D(\Bun_{H})$. Here   
$$
\id\boxtimes F_G: \D(\tilde X\times \Bun_{\tilde H})\to \D(\tilde X\times\Bun_G)
$$ 
is the corresponding theta-lifting functor, and $\pi\times\id: \tilde X\times\Bun_G\to X\times\Bun_G$.

\smallskip\noindent
2) For $m>n$ the functor $F_H$ is naturally equipped with the functoriality datum for $\bar\kappa: \check{\GG}\times\Sigma\times\Gm\to \check{\HH}\rtimes\Sigma$.
% the map $i_{\kappa}$ extends to a $\Sigma$-equivariant homomorphism $\kappa=(i_{\kappa},\delta_{\kappa}): \check{\GG}\times\Gm\to\check{\HH}$ with the following property. 
So, there exists an isomorphism
\begin{equation}
\label{iso_theta_lifting_Hecke_Th1_m>n}
\H^{\la}_H(V, F_H(K))\,\iso\, (\id\boxtimes F_H)(\H^{\la}_G(\Res^{\kappa}(V), K))
\end{equation}
in $\D(\tilde X\times \Bun_H)$ functorial in $V\in\Rep(\check{\HH})$ and $K\in \D(\Bun_G)$. Here  
$$
\id\boxtimes F_H: \D(X\times \Bun_G)\to \D(X\times \Bun_H)
$$ 
is the corresponding theta-lifting functor, and $\pi\times\id: \tilde X\times \Bun_H\to X\times\Bun_H$.  
\end{Th} 
\begin{Rem} 
i) For $\pi: \tilde X\to X$ trivial, the functor $F_G$ (resp., $F_H$) sends a complex with a given central character to the complex with an opposite central character by (\cite{L}, Remark~2). This is why we used the maps $\Gm\to\Gm, z\mapsto z^{-1}$ in our definitions of $i_{\kappa}$ for both cases. 

\smallskip\noindent
ii) The explicit formulas for the restriction of $\bar\kappa$ to the $\Gm$-component are given in Sections~\ref{Sect_Proof_of_Th_local_case_m>n}, \ref{Sect_5.12.6_m_le_n}. If $m=n$ or $m=n+1$ then the restriction of $\bar\kappa$ to $\Gm$ is trivial. If $m\le n$ then $\kappa$ fits into the diagram
$$
\begin{array}{ccc}
\check{\HH}\times\Gm &\toup{\kappa} & \check{\GG}\\ 
\downarrow&& \downarrow\\
\SO_{2m}\times\Gm & \toup{\bar\kappa} & \SO_{2n+1},
\end{array}
$$
where $\bar\kappa$ is the map from (\cite{L2}, Theorem~3).
If $m>n$ then $\kappa$ fits into the diagram
$$
\begin{array}{ccc}
\check{\GG}\times\Gm &\toup{\kappa} & \check{\HH}\\ 
\downarrow&& \downarrow\\
\SO_{2n+1}\times\Gm & \toup{\bar\kappa} & \SO_{2m},
\end{array}
$$
where $\bar\kappa$ is the map from (\cite{L2}, Theorem~3). 
\end{Rem}

\sssec{} For $a\in\ZZ$ let $^a\Bun_{G,H}$ be the stack classifying $\tilde x\in \tilde X$, $(M,\cA)\in\Bun_G$, $(V,\cC,\gamma)\in\Bun_H$, and an isomorphism $\cA\otimes\cC\,\iso\,\Omega(a\pi(\tilde x))$. For 
$$
\cS\in {_{-a}\Sph_{\GG}}, \;\cS'\in {_a\Sph_{\GG}}, \;\cT\in {_{-a}\Sph_{\HH}}, \;\cT'\in {_a\Sph_{\HH}}
$$
we have the Hecke functors defined as in Section~\ref{Sect_Twisted_setting_Hecke_functors}
$$
\H^{\la}_G(\cS, \cdot), \H^{\ra}_G(\cS', \cdot):
\D(\Bun_{G,H})\to \D({^a\Bun_{G,H}})
$$ 
and
$$
\H^{\la}_H(\cT,\cdot), \H^{\ra}_H(\cT', \cdot)
:\D(\Bun_{G,H})\to \D({^a\Bun_{G,H}})
$$ 
Let $\Sigma$ act on $^a\Bun_{G,H}$ via its action on $\tilde X$. We will derive Theorem~\ref{Th_theta-lifting_functors-general} from the following Hecke property of $\Aut_{G,H}$ analogous to (\cite{L2}, Theorem~4). 
 
\begin{Th} 
\label{Th_Hecke_property_Aut_GH} Let $\kappa$ be as in Theorem~\ref{Th_theta-lifting_functors-general}, $a\in\ZZ$.\\
1) For $m\le n$ there exists an isomorphism
\begin{equation}
\label{iso_Th_Hecke_Aut_GH_m_less_n}
\H^{\la}_G(V, \Aut_{G,H})\,\iso\, \H^{\ra}_H(\Res^{\kappa}(V), \Aut_{G,H})
\end{equation}
in $\D(^a\Bun_{G,H})$ functorial in $V\in{_{-a}\Rep(\check{\GG})}$. These isomorphisms are compatible with the tensor structures on $\Rep(\check{\GG}), \Rep(\check{\HH})$ and with the $\Sigma$-actions. 

\smallskip\noindent
2) For $m>n$ there exists an isomorphism
\begin{equation}
\label{iso_Th_Hecke_Aut_GH_m_greater_n}
\H^{\la}_H(V, \Aut_{G,H})\,\iso\, \H^{\ra}_G(\Res^{\kappa}(V), \Aut_{G,H})
\end{equation}
in $\D(^a\Bun_{G,H})$ functorial in $V \in{_{-a}\Rep(\check{\HH})}$. \\
These isomorphisms are compatible with the tensor structures on $\Rep(\check{\GG}), \Rep(\check{\HH})$ and with the $\Sigma$-actions. 
\end{Th} 
 
\ssec{Application: automorphic sheaves on $\Bun_{\GSp_4}$.}  
\label{Section_autom_sheaves_Bun_GSp4}

\sssec{} Keep the notation of Section~\ref{Sect_Theta-lifting functors} assuming $m=n=2$, so $G=\GSp_4$. We get $\check{\HH}\,\iso\, \{(g_1, g_2)\in \GL_2\times\GL_2\mid \det g_1=\det g_2\}$, so that $\Sigma$ permutes the two copies of $\GL_2$. Let $\tilde E$ be an irreducible rank two smooth $\Qlb$-sheaf on $\tilde X$, $\chi$ a rank one local system on $X$ equipped with an isomorphism $\pi^*\chi\,\iso\,\det \tilde E$. To these data one associates the perverse sheaf $K_{\tilde E,\chi, H}$ on $\Bun_H$ introduced in (\cite{L}, Section~5.1).\footnote{In \select{loc.cit.} it was denoted $K_{\tilde E,\chi, \tilde H}$.}

 Recall the construction of  $K_{\tilde E,\chi, H}$. Let $\Bun_{2,\tilde X}$ be the stack of rank 2 vector bundles on $\tilde X$. We have a smooth and surjective map $\delta_{H}: \Bun_{2,\tilde X}\to\Bun_H$ sending $W\in\Bun_{2, \tilde X}$ to $(V,\cC,\gamma)$, where $V\in\Bun_4$ is the descent of $W\otimes\sigma^* W$ with the natural descent data under $\pi$, $\cC=N(\det W)$, and the symmetric form $\Sym^2 V\to\cC$ is the descent under $\pi$ of the natural symmetric form $\Sym^2(W\otimes \sigma^* W)\to \det W\otimes\sigma^*\det W$. Here $\gamma: \cC^{-2}\otimes\det V\,\iso\, \cE$ is a compatible trivialization.

 Let $\Aut_{\tilde E}$ be the perverse Hecke eigensheaf on $\Bun_{2,\tilde X}$ associated to $\tilde E$ normalized as in \cite{FGV}. Under our assumption $\Aut_{\tilde E}$ descends naturally under $\delta_H$ to a perverse sheaf $K_{\tilde E,\chi, H}$ on $\Bun_H$. %Since $\Aut_{\tilde E}$ on $\Bun_{2, \tilde X}$ is cuspidal, $K_{\tilde E,\chi, H}$ is the extension by zero from some open substack of finite type in $\Bun_H$. 

The local system $\pi_*\tilde E^*$ is equipped with a natural symplectic form $\wedge^2 (\pi_*\tilde E^*)\to \chi^{-1}$, so gives rise to a $\check{\GG}$-local system $E_{\check{\GG}}$ on $X$. Since $K_{\tilde E,\chi,\tilde H}$ is a Hecke eigensheaf, Theorem~\ref{Th_theta-lifting_functors-general} implies the following and partially establishes (\cite{L}, Conjecture~2). It also establishes a particular case of Conjecture~\ref{Con_2.2.3}.

\begin{Cor} 
\label{Cor_1}
The complex $F_G(K_{\tilde E,\chi, H})\in\D(\Bun_G)$ is a Hecke eigensheaf corresponding to the $\check{\GG}$-local system $E_{\check{\GG}}$. \QED
\end{Cor}
 
\begin{Rem} i) The nonvanishing of $F_G(K_{\tilde E,\chi, H})$ follows from the compatibility of the functor $F_G$ with the first Whittaker coefficient functors, namely a version of (\cite{LL1}, Theorem~2) for the dual pair $(H, G)$. The proof of (\cite{LL1}, Theorem~2) extends to our case of $(H, G)$ for $n=m=2$. The first Whittaker coefficient of $K_{\tilde E,\chi, H}$ is nonzero, as the first Whittaker coefficient of $\Aut_{\tilde E}$ is known to be $\Qlb$. We conjecture $F_G(K_{\tilde E,\chi, H})$ to be perverse. 

\noindent
ii) If $\pi:\tilde X\to X$ is trivial then fix a numbering of connected components of $\tilde X$. The local system $\tilde E$ becomes a pair of irreducible rank 2 local systems $E_1,E_2$ on $X$ with the isomorphisms $\det E_1\,\iso\,\det E_2\,\iso\,\chi$.
\end{Rem} 
 
\section{Root data} 
\label{Sect_Root data}
 
For the convenience of the reader, in this section we fix notations for the root data of $\HH$, $\GG$ and identify their dual groups with $\GSpin_{2m}$ and $\GSpin_{2n+1}$ respectively. For future references, we define Levi subgroups $Q(\HH)\subset\HH$, $Q(\GG)\subset\GG$ and automorphisms $\tau_{\HH}$ (resp., $\tau_{\GG}$) of $\check{\HH}$ (resp., of $\check{\GG}$).

\ssec{The group $\HH=\GSO_{2m}$}
\label{Sect_3.2_The_group_HH}

\sssec{} We realize $\SO_{2m}$ as the subgroup of $\SL_{2m}$ preserving in the standard representation the symmetric form is given by the matrix
$$
J=\left(
\begin{array}{cc}
0  & E_m\\
E_m & 0
\end{array}
\right),
$$
where $E_m$ is the unity. Set $\HH=\Gm\times\SO_{2m}/(-1, -1)$. We equip $\SO_{2m}$ with the maximal torus of diagonal matrices, and the Borel subgroup preserving for $i=1,\ldots, m$ the isotropic subspace generated by $\{e_1,\ldots, e_i\}$. We assume $a\in\Gm$ acts on the standard representation of $\HH$ as $a$. 

 Write $\Lambda_{\HH}$ (resp., $\check{\Lambda}_{\HH}$) for the coweights (resp., weights) lattice of $\HH$. So, $\check{\Lambda}_{\HH}\subset \ZZ\times \ZZ^m$ is the subgroup of index 2 given by $a+\sum_{i=1}^m b_i=0\mod 2$ for $\check{\lambda}=(a, b)$ with $a\in\ZZ, b\in \ZZ^m$. The subgroup $\Lambda_{\HH}\subset \QQ\times\QQ^m$ is generated by $\ZZ\times\ZZ^m$ and the element 
$
(\frac{1}{2}, \frac{1}{2}(1,\ldots,1))
$. 
The positive roots are 
$$
\begin{array}{ll}
\check{\alpha}_{ij}=(0, \check{e}_i-\check{e}_j),\;\;\mbox{for}\; 1\le i<j\le m; \\
\check{\beta}_{ij}=(0, \check{e}_i+\check{e}_j),\;\;\mbox{for}\; 1\le i<j\le m
\end{array}
$$
The corresponding coroots are $\alpha_{ij}=(0, e_i-e_j), \beta_{ij}=(0, e_i+e_j)$. 

 Let $i_{\HH}\in\Spin_{2m}$ be the nontrivial central element such that $\Spin_{2m}/(i_{\HH})\,\iso\, \SO_{2m}$. Let $\GSpin_{2m}=(\Gm\times\Spin_{2m})/(-1, i_{\HH})$. 
\begin{Lm} The Langlands dual group $\check{\HH}$ identifies canonically with $\GSpin_{2m}$.
\end{Lm}
\begin{proof} The root system for $\Spin_{2m}$ is as follows. The coroots lattice is 
$$
\Lambda_{Spin}=\{(a_1,\ldots, a_m)\in\ZZ^m\mid \sum_i a_i=0\mod 2\}
$$ 
The dual lattice $\check{\Lambda}_{\Spin}$ is generated by $\ZZ^m$, and the element $\frac{1}{2}(1,\ldots, 1)$. The positive roots are $\check{e}_i-\check{e}_j$ and $\check{e}_i+\check{e}_j$ for $1\le i<j\le n$, the corresponding coroots are $e_i-e_j$, $e_i+e_j$.

The root datum of $\GSpin_{2m}=(\Gm\times\Spin_{2m})/(-1, i_{\HH})$ is as follows. The weight lattice $\check{\Lambda}_{\GSpin}$ is a subgroup of $\ZZ\times \check{\Lambda}_{\Spin}$ containing $(2\ZZ)\times\ZZ^n$ and the element $(1, \frac{1}{2}(1,\ldots, 1))$. 
The coweight lattice $\Lambda_{\GSpin}$ is the subgroup of $(\frac{1}{2}\ZZ)\times\ZZ^m$ consisting of $(a,b)$ such that 
$$
a+\frac{1}{2}\sum_i b_i\in\ZZ
$$ 
The roots and coroots are as above. 

 Consider the isomorphism $\check{\Lambda}_{\GSpin}\to \Lambda$ sending $(a,b)$ to $(\frac{a}{2}, b)$. The corresponding isomorphism $\Lambda_{\GSpin}\,\iso\, \check{\Lambda}$ sends $(a, b)$ to $(2a, b)$. It identifies the root datum of $\GSpin_{2m}$ with the dual of that of $\HH$. 
\end{proof} 

\sssec{Levi subgroup $Q(\HH)$} 
\label{Sect_3.2.3_Levi_Q(HH)}
The character $\check{\alpha}_0$ is given by $(2,0)\in\check{\Lambda}_{\HH}$. For $m\ge 2$ the image of $\check{\alpha}_0: \Gm\to \check{\HH}$ is the connected center of $\check{\HH}$.  
 
 Let $Q(\HH)\subset\HH$ be the standard Levi subgroup, where we keep only the positive roots $\check{\alpha}_{ij}$ for $1\le i<j\le m$. Set 
$$
\check{\Lambda}_1=\{(a,b)\in\check{\Lambda}_{\HH}\mid \sum_i b_i=a\}
$$ 
and $\check{\Lambda}_2=\{(a,0)\in\check{\Lambda}_{\HH}\mid a\in 2\ZZ\}$. 
The decomposition $\check{\Lambda}_1\oplus\check{\Lambda}_2=\check{\Lambda}_{\HH}$ yields the dual decomposition $\Lambda_1\oplus \Lambda_2=\Lambda_{\HH}$. Here
$$
\Lambda_1=\{(0,b)\in\Lambda\mid b\in\ZZ^m\}
$$
and $\Lambda_2=\ZZ\bar\omega_{\HH}$, where $\bar\omega_{\HH}=(\frac{1}{2}, -\frac{1}{2}(1,\ldots,1))\in\Lambda$.

 The subgroup of $Q(\HH)$ corresponding to the root datum $(\Lambda_2, \check{\Lambda}_2)$ is identified with $\Gm$ via $\ZZ\,\iso\, \check{\Lambda}_2$, $1\mapsto \check{\alpha}_0$. Let $U_0$ be the irreducible $Q(\HH)$-module with highest weight $(1,b)\in\check{\Lambda}_1$ with $b=(1,0,\ldots, 0)$. Let $C_0$ be the one-dimensional $\HH$-module with highest weight $\check{\alpha}_0$. The above decomposition yields an isomorphism $Q(\HH)\,\iso\, \GL(U_0)\times\GL(C_0)$. Let $V_0$ be the irreducible $\HH$-module with highest weight $(1, b)\in \check{\Lambda}_{\HH}$ with $b=(1,0,\ldots,0)$. This is the standard representation of $\HH$, and $V_0\,\iso\, U_0\oplus U_0^*\otimes C_0$. 
 
  Write $\check{Q}(\HH)$ for the Langlands dual of $Q(\HH)$ over $\Qlb$. Define $\bar U_0$ as the irreducible $\check{Q}(\HH)$-module over $\Qlb$ with highest weight $(0,b)\in \Lambda_{\HH}$ with $b=(1,0,\ldots,0)$. This provides an isomorphism $\GL(\bar U_0)\,\iso\, \check{\GL}(U_0)$ over $\Qlb$. Let $V_{\HH}$ be the irreducible $\check{\HH}$-module over $\Qlb$ with highest weight $(0, b)\in\Lambda_{\HH}$, where $b=(1,0,\ldots, 0)$. Then $V_{\HH}\,\iso\, \bar U_0\oplus \bar U_0^*$ as $\check{Q}(\HH)$-modules. 
  
  Let $\omega_{\HH}\in\Lambda_{\HH}$ be the coweight $(1,0)$. For $m\ge 2$ the image of $\omega_{\HH}$ is the connected center of $\HH$, and $\<\omega_{\HH}, \check{\alpha}_0\>=2$. The element $\omega_{\HH}$ decomposes nontrivially as an element of $\Lambda_1\oplus\Lambda_2$. Write $\alpha_m$ for the character of $\check{Q}(\HH)$ on $\det\bar U_0$. Then $\alpha_m=\omega_{\HH}-2\bar\omega_{\HH}$. We will use $\alpha_m=(0, (1,\ldots, 1))\in\Lambda_{\HH}$ in Proposition~\ref{Pp_geom_theta_lifting_GL_n_GL_m}. 
  
\sssec{}  Since $J=J^{-1}$, there is an automorphism of $\SO(V_{\HH})$ given by $g\mapsto (^tg)^{-1}$, where $^tg$ is the transpose of the matrix $g$. It lifts to a unique automorphism of $\Spin(V_{\HH})$ preserving $i_{\HH}$ that we denote $\tau'$. Now the automorphism $(a,g)\mapsto (a^{-1}, \tau')$ of $\Gm\times\Spin(V_{\HH})$ descends to an automorphism of $\check{\HH}$ denoted $\tau_{\HH}$. It induces the equivalence $*: \Sph_{\HH}\,\iso\, \Sph_{\HH}$. The involution $\tau_{\HH}$ preserves $\check{\TT}_{\HH}\subset \check{Q}(\HH)$ and acts on $\Lambda_{\HH}$ as $-1$. 
  
\ssec{The group $\GG=\GSp_{2n}$}

\sssec{} Let $V_{\GG}$ be the standard representation of $\SO_{2n+1}$. We assume the symmetric form is given by the matrix
$$
J=\left(
\begin{array}{ccc}
0 & E_n & 0\\
E_n & 0 & 0\\
0 & 0 & 1
\end{array}
\right),
$$
where $E_n\in\GL_n$ is the unity. We pick the maximal torus of diagonal matrices in $\SO_{2n+1}$, and the Borel subgroup preserving for $i=1,\ldots, n$ the isotropic subspace generated by $\{e_1,\ldots, e_i\}$. The root system of $\Sp_{2n}$ is the dual to that of $\SO_{2n+1}$. 

 Recall our notation $\GG=\GSp_{2n}=(\Gm\times\Sp_{2n})/(-1, -1)$. We assume $a\in\Gm$ acts on the standard representation of $\Sp_{2n}$ as $a$. Write $\Lambda_{\GG}$ (resp., $\check{\Lambda}_{\GG}$) for the coweights (resp., weights) lattice of $\TT_{\GG}$. The subgroup $\check{\Lambda}_{\GG}\subset \ZZ\times \ZZ^n$ consists of $(a,b)\in\ZZ\times\ZZ^n$ such that $a+\sum_i b_i=0\mod 2$. The subgroup $\Lambda_{\GG}\subset \QQ\oplus\QQ^n$ is generated by $\ZZ\times\ZZ^n$ and the element $(\frac{1}{2}, \frac{1}{2}(1,\ldots, 1))$. The positive roots of $\GG$ are
$$
\begin{array}{l}
\check{\alpha}_{ij}=(0, \check{e}_i-\check{e}_j), \;\;\mbox{for}\; 1\le i<j\le n; \\
\check{\beta}_{ij}=(0, \check{e}_i+\check{e}_j),\;\;\mbox{for}\; 1\le i\le j\le n
\end{array}
$$
The corresponding coroots are $\alpha_{ij}=(0, e_i-e_j)$, $\beta_{ij}=(0, e_i+e_j)$ for $1\le i<j\le n$, and $\beta_{ii}=(0, e_i)$ for $1\le i\le n$. 
 
 Let $i_{\GG}\subset \Spin_{2n+1}$ be the central element such that $\Spin_{2n+1}/(i_{\GG})\,\iso\,\SO_{2n+1}$. Let $\GSpin_{2n+1}=(\Gm\times\Spin_{2n+1})/(-1, i_{\GG})$. 
  
\begin{Lm} The Langlands dual group $\check{\GG}$ identifies canonically with $\GSpin_{2n+1}$. 
\end{Lm}
\begin{proof} The root datum for $\Spin_{2n+1}$ is as follows. The coweights lattice is $\{b\in \ZZ^n\mid \sum_i b_i=0\mod 2\}$, the weights lattice is the subgroup of $\QQ^n$ generated by $\ZZ^n$ and the element $\frac{1}{2}(1,\dots, 1)$. The roots and coroots are obtained from those given above for $\Sp_{2n}$ by passing to the dual root datum.

 Let $\Lambda_{\GSpin}, \check{\Lambda}_{\GSpin}$ denote the coweights and weights lattice of $\GSpin_{2n+1}$. So, $\check{\Lambda}_{\GSpin}$ is the subgroup of $\QQ\times \QQ^n$ generated by $(2\ZZ)\times \ZZ^n$ and the element $(1, \frac{1}{2}(1,\dots, 1))$. The group $\Lambda_{\GSpin}\subset (\frac{1}{2}\ZZ)\times \ZZ^n$ consists of $(a,b)\in (\frac{1}{2}\ZZ)\times \ZZ^n$ such that 
$$
a+\frac{1}{2}\sum_i b_i\in\ZZ
$$ 
The roots and coroots are as above, they are of the form $(0, \cdot)$. 
 
  We have an isomorphism $\Lambda_{\GSpin}\,\iso\, \check{\Lambda}$ sending $(a,b)$ to $(2a, b)$. The isomorphism $\check{\Lambda}_{\GSpin}\,\iso\, \Lambda$ sends 
$(a,b)$ to $(\frac{a}{2}, b)$. This is an isomorphism of the root datum of $\GSpin_{2n+1}$ to the dual of the root datum of $\GG$.
\end{proof}  

\sssec{Levi subgroup $Q(\GG)$} 
\label{Sect_3.3.3_Levi_Q(GG)}
The character $\check{\omega}_0\in\check{\Lambda}_{\GG}$ is given by $(2,0)$. We write $A_0$ for the one-dimensional representation of $\GG$ with weight $\check{\omega}_0$, $M_0$ for the irreducible representation of $\GG$ with highest weight $(1,b)\in\check{\Lambda}_{\GG}$ with $b=(1,0\ldots, 0)\in\ZZ^n$. So, $M_0$ is the standard representation of $\GG$. The image of $\check{\omega}_0: \Gm\to \check{\GG}$ is the center of $\check{\GG}$. 

 Let $Q(\GG)\subset\GG$ be the standard Levi subgroup, where we keep only the positive roots $\check{\alpha}_{ij}$ for $1\le i<j\le n$. Let 
$$
\check{\Lambda}_1=\{(a,b)\in \check{\Lambda}_{\GG}\mid \sum_i b_i=a\}
$$  
and $\check{\Lambda}_2=\ZZ\check{\omega}_0$. Then $\check{\Lambda}_{\GG}=\check{\Lambda}_1\oplus\check{\Lambda}_2$. Let $\Lambda_{\GG}=\Lambda_1\oplus\Lambda_2$ be the dual decomposition. We get
$$
\Lambda_1=\{(0, b)\in \Lambda\mid b\in\ZZ^n\}
$$
and $\Lambda_2=\ZZ\bar\omega_{\GG}$ with $\bar\omega_{\GG}=(\frac{1}{2}, -\frac{1}{2}(1,\ldots,1))\in\Lambda$ normalized by $\<\bar\omega_{\GG}, \check{\omega}_0\>=1$. The group corresponding to the root datum $(\Lambda_2, \check{\Lambda}_2)$ is identified with $\Gm$ via $\ZZ\,\iso\,\check{\Lambda}_2$, $1\mapsto \check{\omega}_0$. 

  Let $L_0$ be the irreducible $Q(\GG)$-module with highest weight $(1, b)\in \check{\Lambda}_1$ with $b=(1,0,\ldots, 0)\in\ZZ^n$. This provides an isomorphism $Q(\GG)\,\iso\, \GL(L_0)\times\GL(A_0)$. We have $M_0\,\iso\, L_0\oplus L_0^*\otimes A_0$ as $Q(\GG)$-modules. 
  
   Write $\check{Q}(\GG)$ for the Langlands dual of $Q(\GG)$ over $\Qlb$. Write $\bar L_0$ for the irreducible $\check{Q}(\GG)$-module with highest weight $(0, b)\in \Lambda_1$ with $b=(1,0,\ldots,0)$. This yields an isomorphism $\GL(\bar L_0)\,\iso\, \check{\GL}(L_0)$ over $\Qlb$. Note that $V_{\GG}$ is the irreducible $\check{\GG}$-module with highest weight  $(0, b)\in\Lambda_{\GG}$, where $b=(1,0,\ldots, 0)$. We have $V_{\GG}\,\iso\, \bar L_0\oplus \bar L_0^*\oplus\Qlb$ as $\check{Q}(\GG)$-modules.

 Set $\omega_{\GG}=(1,0)\in\Lambda_{\GG}$. The image of $\omega_{\GG}: \Gm\to \GG$ is the center of $\GG$, and $\<\omega_{\GG}, \check{\omega}_0\>=2$. As for $\HH$, the decomposition of $\omega_{\GG}$ in $\Lambda_1\oplus\Lambda_2$ is nontrivial. Write $\omega_n$ for the character of $\check{Q}(\GG)$ on $\det \bar L_0$. Then $\omega_{\GG}-2\bar\omega_{\GG}=\omega_n$. We will use $\omega_n$ in Proposition~\ref{Pp_geom_theta_lifting_GL_n_GL_m}. 
 
\sssec{}  Since $J=J^{-1}$, there is an automorphism $g\mapsto (^tg)^{-1}$ of $\SO(V_{\GG})$, where $^tg$ is the transpose of  $g$. It lifts to a unique automorphism of $\Spin(V_{\GG})$ preserving $i_{\GG}$ that we denote by $\tau'$. Now the automorphism $(a, g)\mapsto (a^{-1}, \tau')$ of $\Gm\times \Spin(V_{\GG})$ descends to an automorphism of $\check{\GG}$ that we denote $\tau_{\GG}$. It induces the equivalence $*: \Sph_{\GG}\,\iso\, \Sph_{\GG}$. 
 
 The involution $\tau_{\GG}$ preserves $\check{\TT}_{\GG}\subset \check{Q}(\GG)$, and acts on $\Lambda_{\GG}$ as $-1$.  
 
\ssec{Case $m\le n$}
\label{Sect_3.4_m_le_n}

\sssec{} 
\label{Sect_3.4.1_case_m_le_n}
Consider the map $i_{\kappa}: \check{\HH}\to\check{\GG}$ defined in Section~\ref{Sect_2.4.2_case_m_le_n}. It fits into the commutative diagram
$$
\begin{array}{ccccc}
\check{\TT}_{\GG} & \subset & \check{Q}(\GG) &\subset & \check{\GG}\\
\uparrow\lefteqn{\scriptstyle i_T} && \uparrow\lefteqn{\scriptstyle i_Q} &&\uparrow\lefteqn{\scriptstyle i_{\kappa}}\\
\check{\TT}_{\HH} & \subset & \check{Q}(\HH) &\subset & \check{\HH}  
\end{array}
$$   
Moreover, $\tau_{\GG}$ preserves all the groups in the above diagram, and induces the automorphism $\tau_{\HH}$ on $\check{\HH}$. The corresponding map $\check{\Lambda}_{\HH}\to \check{\Lambda}_{\GG}$ sends $(a,b)$ to $(-a, b')$ with $b'=(b_1,\ldots, b_m, 0,\ldots, 0)$. 
 
\ssec{Case $m>n$}
\label{Sect_Case m>n_3.5}

\sssec{} 
\label{Sect_3.5.1_case_m>n}
Consider the map $i_{\kappa}: \check{\GG}\to\check{\HH}$ defined in Section~\ref{Sect_2.4.3_m_bigger_n}. It fits into the commutative diagram
$$
\begin{array}{ccccc}
\check{\TT}_{\HH} & \subset & \check{Q}(\HH) &\subset & \check{\HH}\\
\uparrow\lefteqn{\scriptstyle i_T} && \uparrow\lefteqn{\scriptstyle i_Q} &&\uparrow\lefteqn{\scriptstyle i_{\kappa}}\\
\check{\TT}_{\GG} & \subset & \check{Q}(\GG) &\subset & \check{\GG}  
\end{array}
$$   
\begin{Lm} 
\label{Lm_3.5.2_about_tau_GG}
The involution $\tau_{\HH}$ preserves all the groups in the above diagram, and induces an automorphism $\tau_{\GG}$ on $\check{\GG}$, which is a Chevalley involution. The corresponding map $\check{\Lambda}_{\GG}\to \check{\Lambda}_{\HH}$ sends $(a,b)\in \check{\Lambda}_{\GG}$ to $(-a, b')$ with $b'=(b_1,\ldots, b_n, 0,\ldots, 0)$. 
\end{Lm}
\begin{proof}
Clearly, $\tau_{\HH}$ preserves the maximal torus of diagonal matrices of $\check{\GG}$, and acts on it as $z\mapsto z^{-1}$.
To check that $\tau_{\HH}$ preserves the subgroup $\check{\GG}$, it suffices to see how $\tau_{\HH}$ acts on the roots of $\check{\GG}$. Since $\tau_{\HH}$ clearly preserves $\check{Q}(\GG)$, we only have to do the verification for roots $\beta_{ii}$. This calculation is left to a reader.
\end{proof}
 
\section{Local theory}
\label{Sect_Local_theory}

In this section we recall the construction of the Weil category geometrizing the local Weil representation from \cite{LL}. We also explain our approach to the geometric version of the representation obtained from the previous one by induction along $\Gm$. 

\ssec{Background on the Weil category}
\label{Sect_Background_Weil}

\sssec{} Remind the following constructions from \cite{LL}. 
Let $W$ be a symplectic Tate space over $k$. By definition (\cite{BD}, 4.2.13), $W$ is a complete topological $k$-vector space having a base of neighbourhoods of 0 consisting of commesurable vector subspaces (i.e., $\dim U_1/(U_1\cap U_2) <\infty$ for any $U_1, U_2$ from this base). It is equipped with a continuous symplectic form $\wedge^2 W\to k$, it induces a topological isomorphism $W\,\iso\, W^*$. 

 For a $k$-subspace $L\subset W$ write $L^{\perp}=\{w\in W\mid \<w,l\>=0\;\mbox{for all}\; l\in L\}$. Write $\cL_d(W)$ for the scheme of discrete lagrangian lattices in $W$. For a c-lattice $R\subset W$ let $\cL_d(W)_R\subset \cL_d(W)$ be the open subscheme of $L\in \cL_d(W)$ satisfying $L\cap R=0$. 

 For a $k$-point $L\in\cL_d(W)$ one defines the category $\cH_L$ as in (\cite{LL}, Section~6.1). Let us remind the definition. For a c-lattice $R\subset R^{\perp}\subset W$ with $R\cap L=0$ we have a lagrangian subspace $L_R:=L\cap R^{\perp}\in\cL(R^{\perp}/R)$ and the Heisenberg group $H_R=(R^{\perp}/R)\oplus k$. Let $\cH_{L_R}$ be the category of perverse sheaves on $H_R$, which are $(\bar L_R, \chi_{L,R})$-equivariant under the left multiplication on $H_R$. Here $\bar L_R=L_R\times\A^1\subset H_R$ and $\chi_{L,R}$ is the local system $\pr^*\cL_{\psi}$ for the projection $\pr: \bar L_R\to\A^1$ sending $(l,a)$ to $a$. Let $D\cH_{L_R}$ be the $\DG$-category of compexes on $H_R$ that are $(\bar L_R, \chi_{L,R})$-equivariant. The forgetful functor $D\cH_{L_R}\to \D(H_R)$ is fully faithful.
 
 For another c-lattice $S\subset R$ we have (an exact for the perverse t-structures) transition functor $T^L_{S,R}: \D\cH_{L_R}\to \D\cH_{L_S}$ (cf. \select{loc.cit.}, Section~6.1). Now $\D\cH_L=\colim \D\cH_{L_R}$ in $\DGCat_{cont}$ over the partially ordered set of c-lattices $R\subset R^{\perp}$ such that $R\cap L=0$.
 Then $\cH_L\subset \D\cH_L$ is the heart of the perverse t-structure on $\D\cH_L$. 
   
 Given a c-lattice $M$ in $W$, we have a $\ZZ/2\ZZ$-graded line bundle on $\cL_d(W)$, whose fibre at $L$ is $\det(M:L)$. Remind that 
$$
\det(M:L)=\det(M\oplus L\to W),
$$
where the complex $M\oplus L\to W$ is placed in cohomological degrees $0$ and $1$. If $S\subset M\subset S^{\perp}$ is a c-lattice with $S\cap L=0$ then $\det(M:L)\,\iso\, \det(M/S)\otimes\det L_S$, where $L_S:=L\cap S^{\perp}$. Note that $\det(M:L)\,\iso\, \det(M^{\perp}:L)$ canonically. If $M'\subset W$ is another c-lattice then we have $\det(M:L)\,\iso\, \det(M:M')\otimes\det(M':L)$ canonically. If $R'\subset W$ is a lagrangian c-lattice then, as $\ZZ/2\ZZ$-graded, $\det(M:L)$ is of parity $\dim(R':M)\mod 2$. 

 Fix a one-dimensional $\ZZ/2\ZZ$-graded space $\cJ_W$ placed in degree $\dim(R':M)\mod 2$. Let $\cA_d$ be the $\ZZ/2\ZZ$-graded purely of degree zero line bundle on $\cL_d(W)$ with fibre $\cJ_W\otimes \det(M:L)$ at $L$. Let $\wt\cL_d(W)$ be the $\mu_2$-gerbe of square roots of $\cA_d$.
 
 For $k$-points $N^0, L^0\in \wt\cL_d(W)$ one associates to them in a canonical way a functor $\cF_{N^0, L^0}: \D\cH_L\to \D\cH_N$ sending $\cH_L$ to $\cH_N$ (defined as in \cite{LL}, Section~6.2). Let us  precise some details. For a c-lattice $R\subset R^{\perp}$ in $W$ we have the projection 
$$
\cL_d(W)_R\to\cL(R^{\perp}/R)
$$
sending $L$ to $L_R$. Let $\cA_R$ be the $\ZZ/2\ZZ$-graded purely of degree zero line bundle on $\cL(R^{\perp}/R)$ whose fibre at $L_1$ is $\det L_1\otimes\det(M:R)\otimes\cJ_W$. Its pullback to $\cL_d(W)_R$ identifies canonically with $\cA_d$, hence a morphism of stacks
\begin{equation}
\label{map_L_d(W)_R_quotient_by_R}
\wt\cL_d(W)_R\to \wt\cL(R^{\perp}/R)
\end{equation}
where $\wt\cL(R^{\perp}/R)$ is the gerbe of square roots of $\cA_R$. Write $N^0_R, L^0_R$ for the images of $N^0, L^0$ under (\ref{map_L_d(W)_R_quotient_by_R}). By definition, the enhanced structure on $L_R$ and $N_R$ is given by one-dimensional spaces $\cB_L, \cB_N$ equipped with
$$
\cB_L^2\,\iso\, \det L_R\otimes\det(M:R)\otimes \cJ_W,\;\; \cB_N^2\,\iso\, \det N_R\otimes\det(M:R)\otimes \cJ_W,
$$
hence an isomorphism $\cB^2\,\iso\, \det L_R\otimes\det N_R$ for $\cB:=\cB_L\otimes\cB_N\otimes \det(M:R)^{-1}\otimes \cJ_W^{-1}$. Write 
\begin{equation}
\label{functor_CIO_N_L_finitedim}
\cF_{N_R^0, L^0_R}: \D\cH_{L_R}\to \D\cH_{N_R}
\end{equation}
for the canonical intertwining functor corresponding to $(N_R, L_R,\cB)$ (as in \select{loc.cit}, Section~6.2). Then $\cF_{N^0, L^0}$ is defined as the colimit of the functors (\ref{functor_CIO_N_L_finitedim}) over the partially ordered set of c-lattices $R\subset R^{\perp}$ such that $N,L\in \cL_d(W)_R$. 

 The proof of (Theorem~2, \cite{LL}) holds through, so for a $k$-point $L^0\in\wt\cL_d(W)$ we have the functor $\cF_{L^0}: \D\cH_L\to \D(\wt\cL_d(W))$ exact for the perverse t-structures. The $\DG$-category $\D(\wt\cL_d(W))$ was defined in (\cite{L2}, Appendix A,B). For two $k$-points
$L^0, N^0\in \wt\cL_d(W)$ the diagram is canonically commutative
$$
\begin{array}{ccc}
\D\cH_L\;\; & \;\;\toup{\cF_{L^0}} & \D(\wt\cL_d(W))\\
\downarrow\lefteqn{\scriptstyle \cF_{N^0, L^0}}\;\; & \;\;\;\nearrow\lefteqn{\scriptstyle \cF_{N^0}}\\
\D\cH_N\;\;
\end{array}
$$ 
Let $W(\wt\cL_d(W))$ denote the essential image of $\cF_{L^0}: \cH_L\to \P(\wt\cL_d(W))$ for any $k$-point $L^0\in \wt\cL_d(W)$. 
% maybe define here the smallest $\DG$-subcategory generated by the essential image?

\begin{Rem} 
\label{Rem_4.1.2_DG_Weil_version}
Probably, one could define the $\DG$-version of the abelian category $W(\wt\cL_d(W))$ as the full $\DG$-subcategory of $\D(\wt\cL_d(W))$ generated under colimits by the essential image of $\cF_{L^0}$ for any $k$-point $L^0\in \wt\cL_d(W)$. This category is independent of a choice of $L^0$. 
\end{Rem}

\ssec{Inducing along $\Gm$} 
\label{Sect_3.2}

\sssec{} 
\label{Sect_4.2.1_for_inducing_along}
Recall that $\cO$ denotes a complete discrete valuation $k$-algebra, $F$ its fraction field. Write $\Omega$ for the completed module of relative differentials of $\cO$ over $k$. For a free $\cO$-module $V$ of finite rank write $V(r)\subset V(F)$ for the $\cO$-submodule $t^{-r}V$, where $t\in\cO$ is a uniformizer. By a $\cO$-lattice in $V(F)$ we mean a $\cO$-submodule $V'\subset V(F)$ such that the natural map $V'(F)\to V(F)$ is an isomorphism.

 Fix $n\ge 1$. For $r\in\ZZ$ let $W_r$ be a free $\cO$-module of rank $2n$ with symplectic form $\wedge^2 W_r\to \Omega(r)$. Then $W_r(F)$ is a symplectic Tate space with the form $\wedge^2 W_r(F)\to \Omega(F)\toup{\Res} k$. Set 
$$
\cL_d^{ex}=\sqcup_{r\in\ZZ} \;\cL_d(W_r(F))
$$  

 Let $\cG_{b,a}$ be the variety of $F$-linear isomorphisms $g: W_a(F)\to W_b(F)$ of symplectic $F$-spaces. Let $G_r=\Sp(W_r)$ as a group scheme over $\cO$. We hope no confusion with our notation $G_r$ from Section~\ref{Sect_2.4.1_now} will occur, as they are in different contexts (global and local). 
 
  Note that $\det W_r\,\iso\,\Omega^n(nr)$. Fix a $\ZZ/2\ZZ$-graded line $\cJ_r$ placed in degree $nr\!\mod 2$. Let $\cA_{d,r}$ be the $\ZZ/2\ZZ$-graded purely of degree zero line bundle on $\cL_d(W_r(F))$ whose fibre at $L$ is $\cJ_r\otimes \det(W_r:L)$. Let $\wt\cL_d(W_r(F))$ be the $\mu_2$-gerbe of square roots of $\cA_{d,r}$. 
  
 Let $\wt\cG_{b,a}$ be the $\mu_2$-gerbe over $\cG_{b,a}$ classifying $g\in \cG_{b,a}$, a one-dimensional space $\cB$ and an isomorphism $\cB^2\,\iso\, \cJ_b\otimes \cJ_a^{-1}\otimes \det(W_b: gW_a)$. The composition $\cG_{c,b}\times \cG_{b,a}\to \cG_{c,a}$ lifts to a morphism $\wt\cG_{c,b}\times\wt\cG_{b,a}\to \wt\cG_{c,a}$ sending $(g_2,\cB_2)\in \wt\cG_{c,b}$, $(g_1,\cB_1)\in \wt\cG_{b,a}$ to $(g_2g_1,\cB)$, where $\cB=\cB_1\otimes\cB_2$. 
  
 Consider the action map
$$
\wt\cG_{b,a}\times \wt\cL_d(W_a(F))\to \wt\cL_d(W_b(F))
$$
sending $(g,\cB)\in\wt\cG_{b,a}$ and $(L, \cB_L)\in \wt\cL_d(W_a(F))$ to $(gL, \cB_1)$, where $\cB_1=\cB\otimes\cB_L$ is equipped with the induced isomorphism 
$$
\cB_1^2\,\iso\, \cJ_b\otimes \det(W_b: gL)
$$ 
In this way $\wt\cG^{ex}:=\sqcup_{a,b\in\ZZ} \;\wt\cG_{b,a}$ becomes a groupoid acting on 
$$
\wt\cL^{ex}_d:=\mathop{\sqcup}\limits_{r\in\ZZ}\wt\cL_d(W_r(F))
$$
 
 The gerbe $\wt\cG_{a,a}\to \cG_{a,a}$ has a canonical section over $G_a(\cO)\subset \cG_{a,a}$ sending $g\in G_a(\cO)$ to $(g,\cB=k)$ equipped with $\id: \cB^2\,\iso\, \det(W_a: W_a)$. One defines the $\DG$-category $\D_{G_a(\cO)}(\wt\cL_d(W_a(F)))$ as in (\cite{L2}, Section~8.2.2). 

 For $g\in \cG_{b,a}$ and a c-lattice $R\subset R^{\perp}\subset W_a(F)$ we have $(gR)^{\perp}=g(R^{\perp})$, and $g$ induces an isomorphism of symplectic spaces 
\begin{equation}
\label{iso_sympl_spaces} 
g: R^{\perp}/R\,\iso\, (gR)^{\perp}/(gR)
\end{equation} 
If $L\in \cL_d(W_a(F))_R$ then $g$ yields an equivalence $\cH_{L_R}\,\iso\, \cH_{gL_{gR}}$ sending $K$ to $g_*K$ for the map $g: H_R\,\iso\, H_{gR}$. Passing to the colimit by $R$, we further get equivalences $g: \cH_L\,\iso\, \cH_{gL}$ and $g: \D\cH_L\,\iso\, \D\cH_{gL}$.
 
\begin{Pp} 
\label{Pp_cG_ab_action}
Let $a,b\in\ZZ$, $\tilde g\in\wt\cG_{b,a}$ over $g\in\cG_{b,a}$ and $L^0\in\wt\cL_d(W_a(F))$ be $k$-points. Then the diagram is canonically commutative
$$
\begin{array}{ccc}
\D\cH_L & \toup{\cF_{L^0}} & \D(\wt\cL_d(W_a(F)))\\
\downarrow\lefteqn{\scriptstyle g} && \downarrow\lefteqn{\scriptstyle \tilde g}\\
\D\cH_{gL} & \toup{\cF_{\tilde g L^0}} & \D(\wt\cL_d(W_b(F)))
\end{array}
$$
\end{Pp}
\begin{Prf}
Let $R\subset R^{\perp}\subset W_a(F)$ be a c-lattice with $R\cap L=0$. We get an equivalence $g:\cH_{L_R}\,\iso\, \cH_{gL_{gR}}$. Let $\cA_R$ be the line bundle on $\cL(R^{\perp}/R)$ whose fibre at $L_1$ is 
$$
\cJ_a \otimes\det(W_a:R)\otimes\det L_1
$$ 
Let $\wt\cL(R^{\perp}/R)$ be the $\mu_2$-gerbe of square roots of $\cA_R$. We have the projection 
$$
\wt\cL_d(W_a(F))_R\to \wt\cL(R^{\perp}/R)
$$ 
sending $L^0$ to $L^0_R$. As in (\cite{LL}, Section~6.4), we have the functors $\cF_{L_R^0}: \D\cH_{L_R}\to \D(\wt\cL(R^{\perp}/R))$. It suffices to show that the diagram is canonically commutative
\begin{equation}
\label{diag_for_section_3.2}
\begin{array}{ccc}
\D\cH_{L_R} & \toup{\cF_{L_R^0}} & \D(\wt\cL(R^{\perp}/R))\\
\downarrow\lefteqn{\scriptstyle g} && \downarrow\lefteqn{\scriptstyle \tilde g}\\
\D\cH_{gL_{gR}} & \toup{\cF_{\tilde g L_{gR}^0}} & \D(\wt\cL((gR)^{\perp}/gR))
\end{array}
\end{equation}
The above expression $\tilde g L_{gR}^0$ is the image of $\tilde g (L^0)$ under $\wt\cL_d(W_b(F))_{gR}\to \wt\cL((gR)^{\perp}/(gR))$. Note that 
$\tilde g L_{gR}^0=\tilde g(L^0_R)$, where 
$$
\tilde g: \wt\cL(R^{\perp}/R)\,\iso\, \wt\cL((gR)^{\perp}/gR)
$$ 
sends $(L_1,\cB)$ to $(gL_1, \cB\otimes\cB_0)$. Here $\tilde g=(g,\cB_0)$. 

 Remind that $H_R$ denotes the Heisenberg group $(R^{\perp}/R)\times\A^1$. For the isomorphism
$$
\tilde g:  \wt\cL(R^{\perp}/R)\times \wt\cL(R^{\perp}/R)\times H_R\,\iso\, \wt\cL((gR)^{\perp}/gR)\times \wt\cL((gR)^{\perp}/gR)\times H_{gR}
$$
we have $\tilde g^* F\,\iso\, F$ canonically, where $F$ is the \select{canonical intertwining operators} sheaf on each side (introduced in \cite{LL}, Theorem~1). 
The commutativity of (\ref{diag_for_section_3.2}) follows.  
\end{Prf} 

\medskip

 By Proposition~\ref{Pp_cG_ab_action}, each $\tilde g\in\cG_{b,a}$ yields an equivalence 
$$
\tilde g: W(\wt\cL_d(W_a(F)))\,\iso\, W(\wt\cL_d(W_b(F)))
$$ 

\ssec{Models for Levi decompositions} 

\sssec{} Assume we are given for each $a\in\ZZ$ a decomposition $W_a=U_a\oplus U^*_a\otimes\Omega(a)$, 
where $U_a$ is a free $\cO$-module of rank $n$, $U_a$ and $U_a^*\otimes\Omega(a)$ are lagrangians, and 
the form $\omega: \wedge^2 W_a\to\Omega(a)$ is given by $\omega\<u, u^*\>=\<u, u^*\>$ for $u\in U_a, u^*\in U_a^*\otimes\Omega(a)$, where $\<\cdot, \cdot\>$ is the canonical pairing between $U_a$ and $U_a^*\otimes\Omega(a)$. 

\begin{Rem}
\label{Rem_Levi_relative_det}
 If $U_1$ is a free $\cO$-module of finite rank and $U_2\subset U_1(F)$ is a $\cO$-lattice then there is a canonical $\ZZ/2\ZZ$-graded isomorphism 
$$
\det(U_2: U_1)^*\,\iso\, \det(U_1^*\otimes\Omega: U_2^*\otimes\Omega)
$$  
Indeed, if $U_1\subset U_2$ then $U_2/U_1$ and $U_1^*\otimes\Omega/U_2^*\otimes\Omega$ are in duality.
\end{Rem}

\sssec{} 
\label{Sect_3.3.2}
For $a,b\in\ZZ$ let $\cU_{b,a}$ be the variety of $F$-linear isomorphisms $U_a(F)\to U_b(F)$. We have an inclusion $\cU_{b,a}\hook{} \cG_{b,a}$ given by $u\mapsto g=(u, (^tu)^{-1})$. Here $^tu\in\GL(U^*\otimes\Omega)(F)$ is the adjoint operator. By Remark~\ref{Rem_Levi_relative_det}, for $u\in \cU_{b,a}$ we have canonically
$$
\det(W_b: gW_a)\,\iso\, \det(U_b: uU_a)^2\otimes \frac{\det(U_a: U_a(-a))}{\det(U_b: U_b(-b))}
$$
For a free $\cO$-module $\cL$ write $\cL_x=\cL\otimes_{\cO} k$. %Do I need this here???

Assume given a one-dimensional $\ZZ/2\ZZ$-graded purely of degree zero vector space $\cJ_{U,a}$ equipped with $ \cJ_{U,a}^2\,\iso\, \cJ_a\otimes \det(U_a(-a): U_a)
$. This yields a section $\rho_{b,a}:\cU_{b,a}\to \wt\cG_{b,a}$ defined as follows. We send $u\in \cU_{b,a}$ to $(g=(u, (^tu)^{-1})\in \cG_{b,a}, \cB)$, where 
$$
\cB=\cJ_{U,b}\otimes\cJ_{U,a}^{-1}\otimes \det(U_b: uU_a)
$$ 
is equipped with the induced isomorphism 
$$
\cB^2\,\iso\, \cJ_b\otimes \cJ_a^{-1}\otimes\det(W_b: gW_a)
$$ 
The section $\rho$ is compatible with the groupoid structures on $\wt\cG^{ex}$ and $\cU^{ex}=\sqcup_{a,b}\; \cU_{b,a}$. We let $\cU^{ex}$ act on $\wt\cL^{ex}_d$ via $\rho$. 

 The $\DG$-category $\D(U_a^*\otimes\Omega(F))$ is defined as in (\cite{L2}). We use here that $U_a^*\otimes\Omega(F)$ is a placid ind-scheme. 

\begin{Pp} 
\label{Pp_functors_cF_U(F)_r}  For $n\ge 1, a\in\ZZ$ there is a canonical map
$$
\cF_{U_a(F)}: \D(U_a^*\otimes\Omega(F))\to \D(\wt\cL_d(W_a(F)))
$$ 
in $\DGCat_{cont}$ exact for the perverse t-structures. For $u\in \cU_{b,a}$ and $\tilde g=\rho_{b,a}(u)\in \wt\cG_{b,a}$ the diagram is canonically commutative
\begin{equation}
\label{diag_for_Pp2_M(F)_r}
\begin{array}{ccc}
\D(U_a^*\otimes\Omega(F)) & \toup{\cF_{U_a(F)}} & \D(\wt\cL_d(W_a(F)))\\
\downarrow\lefteqn{\scriptstyle u} && \downarrow\lefteqn{\scriptstyle \tilde g}\\
\D(U_b^*\otimes\Omega(F)) & \toup{\cF_{U_b(F)}} & \D(\wt\cL_d(W_b(F))),
\end{array}
\end{equation}
\end{Pp}
\begin{Prf}
\Step 1
Let $R_1\subset R_2\subset U_a(F)$ be c-lattices. Write $\<\cdot,\cdot\>_a$ for the symplectic form on the Tate space $W_a(F)$. For a c-lattice $\cR\subset U_a(F)$ set 
$
\cR'=\{w\in U_a^*\otimes\Omega(F)\mid \<w, r\>_a=0\;\mbox{for all}\; r\in R\}$,
this is a c-lattice in $U_a^*\otimes\Omega(F)$. 

Set $R=R_1\oplus R'_2$ then $R^{\perp}=R_2\oplus R'_1$. Let $U_R=R_2/R_1$ then $U_R\in \cL(R^{\perp}/R)$.
Set $U_R^0=(U_R, \cB)$ equipped with the canonical $\ZZ/2\ZZ$-graded isomorphism
$$
\cB^2\,\iso\, \cJ_a\otimes \det(U_R)\otimes \det(W_a:R),
$$ 
where $\cB=\cJ_{U,a}\otimes \det(U_a: R_1)$.  

 Remind the line bundle $\cA_R$ on $\cL(R^{\perp}/R)$ with fibre $\cJ_a\otimes \det L_1\otimes \det(W_a:R)$ at $L_1$ (cf. the proof of Proposition~\ref{Pp_cG_ab_action}). Let $\wt\cL(R^{\perp}/R)$ be the gerbe of square roots of $\cA_R$. So, $U_R^0\in \wt\cL(R^{\perp}/R)$. 

 Write $H_R$ for the Heisenberg group $(R^{\perp}/R)\times\A^1$ and $\cH_{U_R}$ for the corresponding category of $(\bar U_R, \chi_{U,R})$-equivariant perverse sheaves on $H_R$. Here $\bar U_R=U_R\times\A^1$ and $\chi_{U,R}$ is the local system $\pr^*\cL_{\psi}$ on $\bar U_R$, where $\pr: \bar U_R\to \A^1$ is the projection.
 
 Let $\cF_{U_R^0}: \D\cH_{U_R}\to \D(\wt\cL(R^{\perp}/R))$ be the corresponding functor (defined as in \cite{LL}, Section~3.6). The lattice $gR\subset W_b(F)$ satisfies the same assumptions, so we have $U_{gR}=gR_2/gR_1\in \cL(gR^{\perp}/gR)$, and $g(R^{\perp})=(gR)^{\perp}$. 
Further, $U^0_{gR}=(U_{gR}, \cB_1)$ with 
$$
\cB_1=\cJ_{U,b}\otimes \det(U_b: uR_1)
$$ 
equipped with the canonical isomorphism
$\cB_1^2\,\iso\, \cJ_b\otimes\det (U_{gR})\otimes\det(W_b:gR)$.  
 
 We have $\tilde g=(g, \cB_0)$, where 
$$
\cB_0=\cJ_{U,b}\otimes \cJ_{U,a}^{-1}\otimes\det(U_b: uU_a)
$$ 
is equipped with $\cB_0^2\,\iso\, \cJ_b\otimes \cJ_a^{-1}\otimes\det(W_b: gW_a)$. It follows that $\tilde g(U^0_R)\,\iso\, U^0_{gR}$ canonically. 

 Further, $g$ yields an equivalence $g: \D\cH_{U_R}\,\iso\, \D\cH_{U_{gR}}$, and the diagram is canonically commutative
\begin{equation}
\label{diag_action_Levi_similitudes}
\begin{array}{ccc}
\D\cH_{U_R} & \toup{\cF_{U_R^0}} & \D(\wt\cL(R^{\perp}/R))\\
\downarrow\lefteqn{\scriptstyle g} && \downarrow\lefteqn{\scriptstyle \tilde g}\\
\D\cH_{U_{gR}} & \toup{\cF_{U_{gR}^0}} & \D(\wt\cL(gR^{\perp}/gR))
\end{array}
\end{equation}
Indeed, this is a consequence of the following isomorphism. We have 
$$
\tilde g: \wt\cL(R^{\perp}/R)\times  \wt\cL(R^{\perp}/R)\times H_R\,\iso\, \wt\cL(gR^{\perp}/gR)\times  \wt\cL(gR^{\perp}/gR)\times H_{gR},
$$
and for this isomorphism $\tilde g^* F\,\iso\, F$ canonically, where $F$ on both sides is the corresponding \select{canonical intertwining operators} sheaf (introduced in \cite{LL}, Theorem~1). 
 
\smallskip 
\Step 2 Given c-lattices $S_1\subset R_1\subset R_2\subset S_2$ in $U_a(F)$, similarly define $S=S_1\oplus S'_2$ and $U_S^0\in \wt\cL(S^{\perp}/S)$ for $S\subset S^{\perp}\subset W_a(F)$. We have a canonical transition functor $T^U_{S,R}: \D\cH_{U_R}\to \D\cH_{U_S}$ defined as in (\cite{LL}, Section~6.6). 
Let $j: \cL(S^{\perp}/S)_R\subset \cL(S^{\perp}/S)$ be the open subscheme of $L$ satisfying $L\cap (R/S)=0$.  
We have a projection 
$$
p_{R/S}: \wt\cL(S^{\perp}/S)_R\to \wt\cL(R^{\perp}/R)
$$ 
sending $(L, \cB_S)$ to $(L_R, \cB_S)$, where $L_R:=L\cap R^{\perp}$. 
It is understood that $\cB_S$ is equipped with 
$$
\cB_S^2\,\iso\, \cJ_a\otimes \det L\otimes \det(W_a: S),
$$
and we used the canonical isomorphism 
$\det L\otimes \det(W_a: S)\,\iso\, \det L_R\otimes \det(W_a:R)$. 

  Set $P_{R/S}=p_{R/S}^*[\dimrel(p_{R/S})]$. Then 
the following diagram is canonically commutative
$$
\begin{array}{ccccc}
\D\cH_{U_R} &  \toup{\cF_{U^0_R}} & \D(\wt\cL(R^{\perp}/R)) & \toup{P_{R/S}} & \D(\wt\cL(S^{\perp}/S)_R)\\
\downarrow\lefteqn{\scriptstyle T^U_{S,R}} &&& \nearrow\lefteqn{\scriptstyle j^*}\\
\D\cH_{U_S} & \toup{\cF_{U^0_S}} & \D(\wt\cL(S^{\perp}/S))
\end{array}
$$
Define $_R\cF_{U_a(F)}$ as the composition
$$
\D\cH_{U_R}  \toup{\cF_{U^0_R}}  \D(\wt\cL(R^{\perp}/R))
\to \D(\wt\cL_d(W_a(F))_R),
$$
where the second arrow is the pullback (exact for the perverse t-structures) with respect to the projection
$\wt\cL_d(W_a(F))_R\to \wt\cL(R^{\perp}/R)$.

 The above diagram shows that the following diagram is also commutative
$$
\begin{array}{ccc}
\D\cH_{U_R} &  \toup{_R\cF_{U_a(F)}} & \D(\wt\cL_d(W_a(F))_R)\\
\downarrow\lefteqn{\scriptstyle T^U_{S,R}} && \uparrow\lefteqn{\scriptstyle j_{S,R}^*}\\
\D\cH_{U_S} & \toup{_S\cF_{U_a(F)}} & \D(\wt\cL_d(W_a(F))_S),
\end{array}
$$ 
where $j_{S,R}: \wt\cL_d(W_a(F))_R\subset \wt\cL_d(W_a(F))_S$ is the natural open immersion. 

 So, define 
$$
\cF_{U_a(F), R}: \D\cH_{U_R}\to \D(\wt\cL_d(W_a(F)))
$$ 
as the functor sending $K_1$ to the following object $K_2$. For c-lattices $S_1\subset R_1\subset R_2\subset S_2$ as above and $S=S_1\oplus S'_2$ declare the restriction of $K_2$ to $\wt\cL_d(W_a(F))_S$ to be 
$$
(_S\cF_{U_a(F)} \comp T^U_{S,R})(K_1)
$$
The corresponding projective system (indexed by such $S$) defines an object $K_2\in \D(\wt\cL_d(W_a(F)))$. 

 Further, passing to the colimit by $R$ (of the above form) the functors $\cF_{U_a(F),R}$ yield the desired functor $\cF_{U_a(F)}: \D(U_a^*\otimes\Omega(F))\to \D(\wt\cL_d(W_a(F)))$. The commutativity of (\ref{diag_for_Pp2_M(F)_r}) follows from the commutativity of (\ref{diag_action_Levi_similitudes}). 
\end{Prf}

\medskip
\begin{Rem} 
\label{Rem_another_prf_Pp_2}
We could also argue differently in Proposition~\ref{Pp_functors_cF_U(F)_r}. For each $a\in\ZZ$ and $L^0\in \wt\cL_d(W_a(F))$ we could first construct an equivalence 
$$
\cF_{U_a(F), L^0}: \D(U_a^*\otimes\Omega(F))\,\iso\, \D\cH_L
$$ 
as in (\cite{LL}, Proposition~5) such that for any $g\in \cU_{b,a}$ the diagram is commutative
$$
\begin{array}{ccc}
\D(U_a^*\otimes\Omega(F)) & \toup{\cF_{U_a(F), L^0}} & \D\cH_L\\
\downarrow\lefteqn{\scriptstyle g} && \downarrow\lefteqn{\scriptstyle g}\\
\D(U_b^*\otimes\Omega(F)) & \toup{\cF_{U_b(F), \tilde g(L^0)}} & \D\cH_{gL}
\end{array}
$$
with $\tilde g=\rho_{b,a}(g)$. Here $\tilde g(L^0)\in \wt\cL_d(W_b(F))$. 
Then we could define $\cF_{U_a(F)}$ as the composition $$
\D(U_a^*\otimes\Omega(F)) \toup{\cF_{U_a(F), L^0}} \D\cH_L \toup{\cF_{L^0}} \D(\wt\cL_d(W_a(F)))
$$ 
The resulting functor would be (up to a canonical isomorphism) independent of $L^0\in \wt\cL_d(W_a(F))$. 
\end{Rem}

\section{Local theory for the dual pair $\GSp_{2n}, \GSO_{2m}$}
\label{Sect_Dual_pair_GSp_GO}

In this section we formulate and prove Theorem~\ref{Th_local_main}, which is our main local result. 

\ssec{} 
\label{Sect_5.1_local_G_H}
As in Section~\ref{Sect_3.2}, let $\cO$ be a complete discrete valuation $k$-algebra, $F$ its fraction field, $\Omega$ the completed module of relative differentials of $\cO$ over $k$. For a free $\cO$-module $M$ we write $M_x=M\otimes_{\cO} k$ for its geometric fibre. 

Fix $n,m\ge 1$, $a\in\ZZ$. Fix free $\cO$-modules $M_a$ of rank $2n$, $V_a$ of rank $2m$, and $A_a, C_a$ of rank one with symplectic form $\wedge^2 M_a\to A_a$, a nondegenerate symmetric form $\Sym^2 V_a\to C_a$, and a compatible trivialization $\det V_a\,\iso\, C_a^m$. Fix an isomorphism $A_a\otimes C_a\,\iso\, \Omega(a)$. For $a=0$ there is an ambiguity with the notations for $V_0, C_0, M_0, A_0$ of Section~\ref{Sect_Root data}. We hope the sense is clear from the context. 

 Set $W_a=M_a\otimes V_a$, it is equipped with the symplectic form $\wedge^2 W_a\to \Omega(a)$ over $\Spec\cO$. For $a\in\ZZ$ as in Section~\ref{Sect_3.2} set 
\begin{equation} 
\label{def_of_cJ_a_Sect_Dual_pair_GSp_GO}
\cJ_a=C_{a,x}^{-anm}
\end{equation}
Now it is of parity zero as $\ZZ/2\ZZ$-graded, because $\rk(W_a)=4nm$. Define $\wt\cL_d(W_a(F))$, $\cG_{b,a}$, $G_a$, $\cA_{d,a}$ and $\wt\cG_{b,a}$ as in Section~\ref{Sect_4.2.1_for_inducing_along}. Recall that $G_a=\Sp(W_a)$ is a group scheme over $\Spec\cO$. 
 
 Let $\GG$, $\HH$ be as in Section~\ref{Sect_Root data}. We view $(M_a, A_a)$ (resp., $(V_a, C_a)$) as a $\GG$-torsor (resp., $\HH$-torsor) on $\Spec\cO$. 
 
 Let $\GG_{b,a}$ be the variety of isomorphisms $M_a(F)\to M_b(F)$ of $\GG$-torsors over $\Spec F$. Let $\HH_{b,a}$ be the variety of isomorphisms $V_a(F)\to V_b(F)$ of $\HH$-torsors over $\Spec F$.
 Let $\cT_{b,a}$ be the variety of pairs $g=(g_1, g_2)$, where $g_1\in\GG_{b,a}$, $g_2\in\HH_{b,a}$ such that $g\in\cG_{b,a}$. That is, the composition 
$$
\Omega(F)\,\iso\, A_a\otimes C_a(F)\toup{g_1\otimes g_2} A_b\otimes C_b(F)\,\iso\, \Omega(F)
$$ 
must equal to the identity. The natural composition map 
 $\cT_{c,b}\times \cT_{b,a}\to \cT_{c,a}$ makes $\cT=\sqcup_{a,b} \;\cT_{b,a}$ into a groupoid. The natural maps $\cT_{b,a}\to \cG_{b,a}$ are compatible with the groupoid structures on $\cT,\cG$. 
 
\begin{Lm} 
\label{Lm_too_general_det_rel}
Let $M_i, V$ be free $\cO$-modules of finite rank, where $M_2\subset M_1(F_x)$ is a $\cO$-lattice. Set $\dim(M_1: M_2)=\dim(M_1/R)-\dim(M_2/R)$ for a $\cO$-lattice $R\subset M_1\cap M_2$. Then we have a canonical $\ZZ/2\ZZ$-graded isomorphism
$$
\det(M_1\otimes V: M_2\otimes V)\,\iso\, \det(M_1:M_2)^{\rk V}\otimes (\det V_x)^{\dim(M_1: M_2)}
\eqno{\square}
$$
%Here $\det(M_1:M_2)$ is of parity $\dim(M_1: M_2)\!\!\mod 2$, and $\det V_x$ is of parity zero as $\ZZ/2\ZZ$-graded.
\end{Lm}
%\begin{Prf} 
%Pick a $\cO$-lattice $R\subset M_1\cap M_2$. It suffices to establish a canonical $\ZZ/2\ZZ$-graded isomorphism
%$$
%\det(M_1\otimes V: R\otimes V)\,\iso\, \det(M_1/R)^{\rk V}\otimes (\det V_x)^{\dim(M_1/R)}
%$$
%To do so, it suffices to pick a flag $R=R_0\subset R_1\subset\ldots\subset R_s=M_1$ of $\cO_x$-lattices with $\dim(R_i/R_{i-1})=1$. 
%\end{Prf}

\sssec{} 
\label{Sect_4.1.2}
For $e\in\ZZ$ set $\GG_{b,a}^e=\{g\in\GG_{b,a}\mid gA_a=A_b(e)\}$ and $\HH_{b,a}^e=\{g\in\HH_{b,a}\mid g C_a=C_b(e) \}$. 
 
  Let us construct a canonical section $\nu_{b,a}: \cT_{b,a}\to\wt\cG_{b,a}$ compatible with the groupoids structures. Let $g=(g_1, g_2)\in\cT_{b,a}$ with $g_1\in \GG_{b,a}^e$, $g_2\in \HH_{b,a}^c$, so  $e+c=a-b$. We have a canonical $\ZZ/2\ZZ$-graded isomorphism
\begin{multline*}
\det(M_b\otimes V_b: (g_1 M_a)\otimes(g_2 V_a))\,\iso\\   \det(M_b\otimes V_b: (g_1 M_a)\otimes V_b)\otimes \det((g_1 M_a)\otimes V_b: (g_1 M_a)\otimes(g_2 V_a))
\end{multline*}
The element $g_1$ gives an isomorphism $\det((g_1M_a)_x)\,\iso\,\det((M_a)_x)$. Applying Lemma~\ref{Lm_too_general_det_rel} we get a canonical $\ZZ/2\ZZ$-graded isomorphism
\begin{multline*}
\det(M_b\otimes V_b: (g_1 M_a)\otimes(g_2 V_a))\,\iso\, \\ \det(M_b: g_1 M_a)^{2m}\otimes \det(V_b: g_2 V_a)^{2n}\otimes (\det V_b)_x^{\dim(M_b: g_1  M_a)}\otimes (\det M_a)_x^{\dim(V_b: g_2 V_a)}\,\iso\\
\det(M_b: g_1 M_a)^{2m}\otimes \det(V_b: g_2 V_a)^{2n}\otimes C_{b,x}^{-mne}\otimes A_{a,x}^{-mnc}\,\iso\\
\det(M_b: g_1 M_a)^{2m}\otimes \det(V_b: g_2 V_a)^{2n}\otimes C_{b,x}^{-mne}\otimes \cC_{a,x}^{mnc}\otimes \cO((1-a)c)_x^{mn}
\end{multline*}  
We used that $\dim(M_b: g_1 M_a)=-ne$, $\dim(V_b: g_2 V_a)=-mc$. Identifying further $C_a\toup{g_2} C_b(c)$, we get
$$
\cJ_b\otimes \cJ_a^{-1}\otimes  \det(W_b: gW_a)\,\iso\, 
\det(M_b: g_1 M_a)^{2m}\otimes \det(V_b: g_2 V_a)^{2n}\otimes C_{b,x}^{2cnm}\otimes \cO(c(1+c))_x^{nm}
$$
Let $\nu_{b,a}(g)=(g,\cB)$, where 
$$
\cB=\det(M_b: g_1 M_a)^m\otimes \det(V_b: g_2 V_a)^n\otimes C_{b,x}^{cnm}\otimes \cO(c(1+c)/2)_x^{nm}
$$
is equipped with the induced isomorphism $\cB^2\,\iso\, \cJ_b\otimes \cJ_a^{-1}\otimes  \det(W_b: gW_a)$. 

 We let $\cT$ act on $\wt\cL_d^{ex}$ via $\nu$. 
 
\ssec{Categories $\D_{\cT_a}(\wt\cL_d(W_a(F)))$}
\label{Sect_5.2_stack_a_cXL}

\sssec{} Let $\GG_a=\GSp(M_a)$ and $\HH_a=\GSO(V_a)$, the connected component of unity of the group scheme $\GO(V_a)$ over $\Spec\cO$. One should not confuse $\GG_a$ with the additive group, for which we use the notation $\AA^1$.
Set 
$$
\cT_a=\{(g_1,g_2)\in (\GG_a\times\HH_a)(\cO)\mid g_1\otimes g_2\;\mbox{acts trivially on}\; A_a\otimes C_a\}
$$ 
Recall the line bundle $\cA_{d,a}$ on $\cL_d(W_a(F))$ with fibre $\cJ_a\otimes\det(W_a: L)$ at $L$. Note that $\cA_{d,a}$ is naturally $\cT_a$-equivariant, %(we underline that $\cT_a$ acts nontrivially on $\cJ_a$). 
so it gives rise to a line bundle denoted $^a\cA_{\cXL}$
on the quotient stack 
$$
^a\cXL=\cL_d(W_a(F))/\cT_a
$$ 
We also have the $\mu_2$-gerbe of square roots of this line bundle denoted
$$
^a\wt\cXL=\wt\cL_d(W_a(F))/\cT_a
$$ 
The $\DG$-category $\D_{\cT_a}(\wt\cL_d(W_a(F)))$ is defined as in (\cite{L2}, Section~8.2.2) or \cite{LL}.

 The stack $^a\cXL$ classifies: a $\GG$-torsor $(M,A)$ over $\Spec\cO$, a $\HH$-torsor $(V,C)$ over $\Spec\cO$ (so, we have a compatible isomorphism $\det V\,\iso\, C^m$), an isomorphism $A\otimes C\,\iso\, \Omega(a)$, and a discrete lagrangian subspace $L\subset M\otimes V(F)$. 
The fibre of $^a\cA_{\cXL}$ at $(M,A,V,C,L)$ is $C_x^{-anm}\otimes \det(M\otimes V: L)$.
 
\ssec{Hecke functors} 
\label{Sect_5.3_Hecke_functors}

\sssec{} Denote by $^{a,a'}\cH_{\GG, \cXL}$ the stack classifying: a point $(L, M,A, V, C)\in {^a\cXL}$, a $\cO$-lattice
 $M'\subset M(F)$ such that for $A'=A(a'-a)$ the induced form $\wedge^2 M'\to A'$ is regular and nondegenerate. We get a diagram
\begin{equation}
\label{diag_Hecke_cXL_for_G}
^a\cXL\; \getsup{h^{\la}}\; {^{a, a'}\cH_{\GG,\cXL}} \;\toup{h^{\ra}} \;{^{a'}\cXL},
\end{equation} 
where $h^{\la}$ (resp., $h^{\ra}$) sends a point of $^{a,a'}\cH_{\cXL}$ to $(L, M,A, V, C)$ (resp., to $(L, M', A', V, C))$. 
 
\begin{Lm} 
\label{Lm_det_calculation_one}
For a point $(L, M,A, M', A', V, C)$ of $^{a,a'}\cH_{\cXL}$ there is a canonical $\ZZ/2\ZZ$-graded isomorphism
$$
C_x^{-a'nm}\otimes\det(M'\otimes V: L)\,\iso\, 
C_x^{-anm}\otimes\det(M\otimes V: L) \otimes \det(M': M)^{2m}
$$
\end{Lm}
\begin{proof}
This follows from Lemma~\ref{Lm_too_general_det_rel}.
\end{proof}

 Let $^{a,a'}\wt\cH_{\GG, \cXL}\toup{\tilde h^{\ra}} {^{a'}\wt\cXL}$ be obtained from $h^{\ra}$ by the base change $^{a'}\wt\cXL\to {^{a'}\cXL}$. By Lemma~\ref{Lm_det_calculation_one}, we get a diagram
\begin{equation}
\label{diag_tilde_Hecke_cXL_for_G}
^a\wt\cXL \;\getsup{\tilde h^{\la}} \;{^{a,a'}\wt\cH_{\GG,\cXL}}\; \toup{\tilde h^{\ra}} \;{^{a'}\wt\cXL}
\end{equation}
Here a point of ${^{a,a'}\wt\cH_{\GG,\cXL}}$ is given by a collection $(L,M,A, M', A', V,C)\in  {^{a,a'}\cH_{\GG,\cXL}}$ together with a one-dimensional space $\cB$ equipped with
$$
\cB^2\,\iso\, C_x^{-a'nm}\otimes\det(M'\otimes V: L)
$$ 
The map $\tilde h^{\la}$ sends this point to $(L,M,A,V,C)\in {^a\cXL}$ together with the one-dimensional space $\cB_1=\cB\otimes \det(M':M)^{-m}$ with the induced isomorphism 
$$
\cB_1^2\,\iso\, C_x^{-anm}\otimes\det(M\otimes V: L)
$$ 
\sssec{} 
\label{Sect_4.3.3}
The affine grassmanian $\Gr_{\GG_a}=\GG_a(F)/\GG_a(\cO)$ is the ind-scheme classifying $\cO$-lattices $R\subset M_a(F)$ such that for some $r\in\ZZ$ the induced form $\wedge^2 R\to A_a(r)$ is regular and nondegenerate. Write $\Gr_{\GG_a}^r$ for the connected component of $\Gr_{\GG_a}$ given by fixing such $r$.  Write $_b\Sph_{\GG_{a'}}\subset \Sph_{\GG_{a'}}$ for the full subcategory of objects that vanish off $\Gr_{\GG_{a'}}^b$. This notation is compatible with the notation $_a\Sph_{\GG}$ from Section~\ref{Sect_2.3.3}. 
 
 Trivializing a point of $^{a'}\cXL$ (resp., of $^a\cXL$) one gets isomorphisms
$$
\id^r: {^{a,a'}\cH_{\GG, \cXL}}\,\iso\, (\cL_d(W_{a'}(F))\times \Gr_{\GG_{a'}}^{a-a'})/\cT_{a'}
$$
and 
$$
\id^l:  {^{a,a'}\cH_{\GG, \cXL}}\,\iso\, (\cL_d(W_a(F))\times \Gr_{\GG_a}^{a'-a})/\cT_a,
$$
where the corresponding action of $\cT_{a'}$ (resp., of $\cT_a$) is diagonal. They lift naturally to a $\cT_{a'}$-torsor 
$$
\wt\cL_d(W_{a'}(F))\times \Gr_{\GG_{a'}}^{a-a'} \to \, {^{a,a'}\wt\cH_{\GG, \cXL}}
$$
and a $\cT_a$-torsor
$$
\wt\cL_d(W_a(F))\times \Gr_{\GG_a}^{a'-a}\to\, {^{a,a'}\wt\cH_{\GG, \cXL}}
$$
 
 So, for $K\in \D_{\cT_a}(\wt\cL_d(W_a(F)))$, $K'\in \D_{\cT_{a'}}(\wt\cL_d(W_{a'}(F)))$, $\cS\in \Sph_{\GG_a}$, $\cS'\in \Sph_{\GG_{a'}}$ we can form their external products
$$
(K\tboxtimes\cS)^l, (K'\tboxtimes\cS')^r
$$
on ${^{a,a'}\wt\cH_{\GG, \cXL}}$. 

 The Hecke functors 
$$
\H^{\ra}_{\GG}
: {_{a'-a}\Sph_{\GG_{a'}}}
\times \D_{\cT_{a}}(\wt\cL_d(W_{a}(F)))\to \D_{\cT_{a'}}(\wt\cL_d(W_{a'}(F)))
$$ 
and
$$
\H^{\la}_{\GG}: {_{a'-a}\Sph_{\GG_{a'}}}
\times \D_{\cT_{a'}}(\wt\cL_d(W_{a'}(F)))\to \D_{\cT_a}(\wt\cL_d(W_a(F)))
$$
are defined by 
$$
\H^{\la}_{\GG}(\cS, K')=(\tilde h^{\la})_! (K'\tboxtimes\ast\cS)^r,\;\;\;\; \H^{\ra}_{\GG}(\cS, K)=(\tilde h^{\ra})_!(K\tboxtimes \cS)^l
$$ 
We used the fact that $\cS\in\Sph_{\GG}$ is the extension by zero from $\Gr_{\GG}^r$ iff $\ast\cS$ is the extension by zero from $\Gr_{\GG}^{-r}$. Since $\D_{\cT_a}(\wt\cL_d(W_a(F)))$ is defined by some limit procedure, the above definition needs further precisions to make it rigorous, this is done as in (\cite{L2}, Section~4.3). We have $\H^{\ra}_{\GG}(\cS, K)\,\iso\, \H^{\la}_{\GG}(\ast \cS, K)$. 
 
\sssec{} Let $^{a,a'}\cH_{\HH,\cXL}$ be the stack classifying: a point $(L,M,A,V,C)\in {^a\cXL}$, a lattice $V'\subset V(F)$ such that for $C'=C(a'-a)$ the induced form $\Sym^2 V'\to C'$ is regular and nondegenerate (we also get the isomorphism $C'^{-m}\otimes\det V'\,\iso\, C^{-m}\otimes \det V\,\iso\,\cO$). As for $\GG$, we get a diagram
\begin{equation}
\label{diag_Hecke_wt_HH_cXL}
\begin{array}{ccccc}
^a\wt\cXL & \getsup{\tilde h^{\la}} &{^{a,a'}\wt\cH_{\HH,\cXL}} &\toup{\tilde h^{\ra}} &{^{a'}\wt\cXL}\\
\downarrow && \downarrow && \downarrow\\
^a\cXL &\getsup{h^{\la}} &{^{a,a'}\cH_{\HH,\cXL}} &\toup{h^{\ra}} & {^{a'}\cXL},
\end{array}
\end{equation}
where $h^{\la}$ (resp. $h^{\ra}$) sends 
\begin{equation}
\label{point_of_aa'cH_HHcXL}
(L,M,A,V,C,V',C')
\end{equation} 
to $(L,M,A,V,C)$ (resp., to $(L,M,A,V',C')$), the vertical arrows are $\mu_2$-gerbes, and the right square is cartesian (thus defining the stack $^{a,a'}\wt\cH_{\HH,\cXL}$). 

 A point of $^{a,a'}\wt\cH_{\HH,\cXL}$ is given by $(L,M,A,V,C,V',C')\in {^{a,a'}\cH_{\HH, \cXL}}$ and a one-dimensional space $\cB$ equipped with
$$
\cB^2\,\iso\, (C'_x)^{-a'nm}\otimes\det(M\otimes V':L)
$$ 
The map $\tilde h^{\la}$ in (\ref{diag_Hecke_wt_HH_cXL}) sends this point to $(L,M,A,V,C)\in {^a\cXL}$, the one-dimensional space $\cB_1$ with the isomorphism
$\cB_1^2\,\iso\, C_x^{-anm}\otimes\det(M\otimes V: L)$ given by Lemma~\ref{Lm_4.3.5} below, where
$$
\cB_1=\cB\otimes C_x^{nm(a'-a)}\otimes \det(V:V')^n\otimes \cO(\frac{1}{2}nm(a-a')(a-a'-1))_x
$$
\begin{Lm} 
\label{Lm_4.3.5}
For a point (\ref{point_of_aa'cH_HHcXL}) of $^{a,a'}\cH_{\HH,\cXL}$ as above one has canonically
$$
\frac{(C'_x)^{-a'nm}\otimes\det(M\otimes V': L)}{(C_x)^{-anm}\otimes\det(M\otimes V: L)}\,\iso\, \cC_x^{2(a-a')nm}\otimes\det(V': V)^{2n}\otimes \cO((a-a')(a'-a+1)mn)_x
$$
\end{Lm}
\begin{proof}
We have $\dim(V': V)=(a'-a)m$. Using $\det M\,\iso\, A^n\,\iso\,(\Omega(a)\otimes\cC^{-1})^n$, from Lemma~\ref{Lm_too_general_det_rel} we get 
\begin{multline*}
\frac{\det(M\otimes V': L)}{\det(M\otimes V: L)}\,\iso\,\det(V': V)^{2n}\otimes \det(M_x)^{\dim(V': V)}\,\iso\\ \det(V': V)^{2n}\otimes(\cO(a-1)\otimes\cC^{-1})_x^{(a'-a)nm}
\end{multline*}
Using $C'=C(a'-a)$, one gets the desired isomorphism.
\end{proof}

\sssec{} The affine grassmanian $\Gr_{\HH_a}$ classifies lattices $V'\subset V_a(F)$ such that the induced symmetric form $\Sym^2 V'\to C_a(b)$ is regular and nondegenerate for some $b\in\ZZ$. Write $\Gr_{\HH_a}^b$ for the locus of $\Gr_{\HH_a}$ given by fixing this $b$. For $m\ge 2$ there is an exact sequence $0\to \ZZ/2\ZZ\to \pi_1(\HH_a)\to \ZZ\to 0$, so if $m\ge 2$ then $\Gr_{\HH_a}^b$ is a union of two connected components of $\Gr_{\HH_a}$. Write $_b\Sph_{\HH_a}\subset\Sph_{\HH_a}$ for the full subcategory of objects that vanish off $\Gr_{\HH_a}^b$. This notation is compatible with the notation $_b\Sph_{\HH}$ from Section~\ref{Sect_2.3.4}. 
 
 The Hecke functors
$$
\H^{\la}_{\HH}: {_{a'-a}\Sph_{\HH_{a'}}}\times  \D_{\cT_{a'}}(\wt\cL_d(W_{a'}(F)))\to \D_{\cT_a}(\wt\cL_d(W_a(F)))
$$
and
$$
\H^{\ra}_{\HH}
: {_{a'-a}\Sph_{\HH_{a'}}}\times  \D_{\cT_{a}}(\wt\cL_d(W_{a}(F)))\to \D_{\cT_{a'}}(\wt\cL_d(W_{a'}(F)))
$$
are defined as in Section~\ref{Sect_4.3.3} using the diagram (\ref{diag_Hecke_wt_HH_cXL}). We have $\H^{\ra}_{\HH}(\cS, K)\,\iso\, \H^{\la}_{\HH}(\ast\cS, K)$. 
 
\sssec{} 
\label{Sect_4.3.7}
For each $a\in\ZZ$ a trivialization $\alpha$ of the $\GG$-torsor $(M_a, A_a)$ on $\Spec\cO$ yields an isomorphism $\bar\alpha:\Gr_{\GG_a}\,\iso\, \Gr_{\GG}$. The induced equivalences $\bar\alpha^*: \Sph_{\GG}\,\iso\, \Sph_{\GG_a}$ are isomorphic for different $\alpha$'s. In what follows we sometimes identify these two categories in this way. Similarly, we identify $\Sph_{\HH_a}\,\iso\, \Sph_{\HH}$.

 The above Hecke actions are extended to the actions of $\D\Sph_{\GG}\,\iso\,\Rep(\check{\GG}\times\Gm)$
and $\DD\Sph_{\HH}\,\iso\,\Rep(\check{\HH}\times\Gm)$
respectively as in Section~\ref{Sect_2.2.1}.

\ssec{} Let $S_{W_0(F)}\in \P_{\cT_0}(\wt\cL_d(W_0(F)))$ be the theta-sheaf introduced in (\cite{LL}, Section~6.5). This is a $\cT_0$-equivariant object of the Weil category $W(\wt\cL_d(W_0(F)))$. Here is the main result of Section~\ref{Sect_Dual_pair_GSp_GO}. 
 
\begin{Th} 
\label{Th_local_main} Let $a\in\ZZ$, $\kappa$ be as in Theorem~\ref{Th_theta-lifting_functors-general}.\\ 
% In particular, $\kappa$ is $\Sigma$-equivariant.\\
1) For $m\le n$ there is an isomorphism in $\D_{\cT_a}(\wt\cL_d(W_a(F)))$
\begin{equation}
\label{iso_Th_local_main_m_le_n}
\H^{\la}_{\GG}(V, S_{W_0(F)})\,\iso\, \H^{\ra}_{\HH}(\Res^{\kappa}(V), S_{W_0(F)})
\end{equation} functorial in $V\in{_{-a}\Rep(\check{\GG})}$ and compatible with the tensor structures on $\Rep(\check{\GG}), \Rep(\check{\HH})$.\\
2) For $m>n$ there is an isomorphism in $\D_{\cT_a}(\wt\cL_a(W_a(F)))$
\begin{equation}
\label{iso_Th_local_main_m_bigger_n}
\H^{\la}_{\HH}(V, S_{W_0(F)})\,\iso\, \H^{\ra}_{\GG}(\Res^{\kappa}(V), S_{W_0(F)})
\end{equation}
functorial in $V\in {_{-a}\Rep(\check{\HH})}$ and compatible with the tensor structures on $\Rep(\check{\GG}), \Rep(\check{\HH})$. \\
\end{Th} 

The proof occupies the rest of Section~\ref{Sect_Dual_pair_GSp_GO}. The explicit formulas for $\kappa$ are found in Sections~\ref{Sect_Proof_of_Th_local_case_m>n}, \ref{Sect_5.12.6_m_le_n}.

\ssec{Levi subgroups} 
\label{Sect_5.5_Levi_subgroups}

\sssec{} Recall the Levi subgroups $Q(\GG)\subset\GG$ and $Q(\HH)\subset\HH$ defined in Sections~\ref{Sect_3.3.3_Levi_Q(GG)}, \ref{Sect_3.2.3_Levi_Q(HH)}. 

Similarly, assume given a decomposition $M_a=L_a\oplus (L_a^*\otimes A_a)$, where $L_a$ is a free $\cO$-module of rank $n$, $L_a$ and $L_a^*\otimes A_a$ are lagrangians, and the form is given by $(l, l^*)\mapsto \<l, l^*\>$, where $\<\cdot,\cdot\>$ is the canonical pairing between $L_a$ and $L_a^*$. Assume given a similar decomposition $V_a=U_a\oplus (U_a^*\otimes C_a)$ for $V_a$, here $U_a$ is a free $\cO$-module of rank $m$. Here $U_a, U_a^*\otimes C_a$ are isotropic, the symmetric form is given by the canonical pairing between $U_a$ and $U_a^*$. 

 Write $Q(\GG_a)\subset \GG_a$ and $Q(\HH_a)\subset\HH_a$ for the Levi subgroups preserving the above decompositions. 
Set
$$
\begin{array}{l}
Q \GG\HH_a=\{g=(g_1,g_2)\in Q(\GG_a)\times Q(\HH_a)\mid g\in \cT_a\} \\ \\
\GG Q\HH_a=\{g=(g_1,g_2)\in \GG_a\times Q(\HH_a)\mid g\in \cT_a\} \\  \\
\HH Q\GG_a=\{g=(g_1,g_2)\in \HH_a\times Q(\GG_a)\mid g\in \cT_a\} 
\end{array}
$$
We view all of them as group schemes over $\Spec\cO$. 

 The affine grassmanian $\Gr_{Q(\GG_a)}$ classifies pairs of lattices $L'\subset L_a(F)$, $A'\subset A_a(F)$. For $b\in\ZZ$ write $\Gr_{Q(\GG_a)}^b$ for the locus of $\Gr_{Q(\GG_a)}$ given by $A'=A_a(b)$. Write $_b\Sph_{Q(\GG_a)}\subset\Sph_{Q(\GG_a)}$ for the full subcategory of objects that vanish off $\Gr_{Q(\GG_a)}^b$. As in Section~\ref{Sect_4.3.7}, we identify $\Sph_{Q(\GG)}\,\iso\,\Sph_{Q(\GG_a)}$. We write $\Gr_{Q(\HH)}^a$ for the union of the connected components $\Gr_{Q(\HH)}^{\theta}$ satisfying $\<\theta, \check{\alpha}_0\>=-a$.

 In view of $\Loc$, the inclusion of the Langlands dual groups $\check{Q}(\GG)\hook{} \check{\GG}$ yields a faithful functor $_b\Sph_\GG\to {_b\Sph_{Q(\GG)}}$ for each $b$, and similarly for $\HH$. 

\sssec{} For $b,a\in\ZZ$ write $Q(\GG_{b,a})$ for the variety of isomorphisms $(L_a(F)\to L_b(F),\; A_a(F)\to A_b(F))$ of $\GL_n\times\Gm$-torsors over $\Spec F$. Let $Q(\HH_{b,a})$ be the variety of isomorphisms $(U_a(F)\to U_b(F),\; C_a(F)\to C_b(F))$ of $\GL_m\times\Gm$-torsors over $\Spec F$. Set 
$$
\begin{array}{l}
Q\GG\HH_{b,a}=\{g=(g_1,g_2)\in Q(\GG_{b,a})\times Q(\HH_{b,a})\mid g\in\cG_{b,a}\} \\ \\
\GG Q\HH_{b,a}=\{g=(g_1,g_2)\in \GG_{b,a}\times Q(\HH_{b,a})\mid g\in\cG_{b,a}\} \\ \\  
\HH Q\GG_{b,a}=\{g=(g_1,g_2)\in \HH_{b,a}\times Q(\GG_{b,a})\mid g\in\cG_{b,a}\}
\end{array}
$$ 

\sssec{} 
\label{Sect_4.5.3}
Set $\Upsilon_a=L_a^*\otimes A_a\otimes V_a$ and $\Pi_a=U_a^*\otimes C_a\otimes M_a$. For $a\in\ZZ$ and any $L^0\in \wt\cL_d(W_a(F))$ we have the equivalences 
$$
\cF_{L_a\otimes V_a(F), L^0}: \D(\Upsilon_a(F))\,\iso\, \D\cH_L
$$ 
and 
$$
\cF_{U_a\otimes M_a(F), L^0}: \D(\Pi_a(F))\,\iso\, \D\cH_L
$$ 
defined as in Remark~\ref{Rem_another_prf_Pp_2}. 

Remind that for free $\cO$-modules of finite type $\cV,\cU$ one has the partial Fourier transform 
$$
\Four_{\psi}: \D(\cV(F)\oplus \cU(F))\,\iso\, \D(\cV^*\otimes\Omega(F)\oplus \cU(F))
$$
normalized to preserve perversity (and purity in the case of a finite base field), cf. (\cite{L2}, Section~4.8) for the definition. 
Thus, the decompositions 
$$
\Pi_a\,\iso\, U^*_a\otimes C_a\otimes L_a\;\oplus\; U^*_a\otimes L^*_a\otimes\Omega(a)
$$
and 
$$
\Upsilon_a\,\iso\, L^*_a\otimes A_a\otimes U_a\;\oplus\;
U_a^*\otimes L_a^*\otimes \Omega(a)
$$
yield the partial Fourier transform, which we denote 
$$
\zeta_a: \D(\Upsilon_a(F))\,\iso\, \D(\Pi_a(F))
$$
One checks that $\zeta_a$ is canonically isomorphic to the functor $\cF_{U_a\otimes M_a(F), L^0}^{-1}
\comp\cF_{L_a\otimes V_a(F), L^0}$ for any $L^0\in \wt\cL_d(W_a(F))$. 

\sssec{} It is convenient to denote $\bar\Upsilon_a=L_a\otimes V_a$ and $\bar \Pi_a=U_a\otimes M_a$. Recall the line $\cJ_a$ given by (\ref{def_of_cJ_a_Sect_Dual_pair_GSp_GO}).  
For the decomposition $W_a=\bar\Pi_a\oplus \bar\Pi_a^*\otimes\Omega(a)$ define a $\ZZ/2\ZZ$-graded line (purely of parity zero)
$$
\cJ_{\bar\Pi, a}=\cO((1-a)a/2)_x^{nm}\otimes \det (U_a(-a): U_a)^n  
$$
Using Lemma~\ref{Lm_too_general_det_rel} we equip it with a natural $\ZZ/2\ZZ$-graded isomorphism 
$$
\cJ_{\bar\Pi,a}^2\,\iso\, \cJ_a\otimes \det(\bar\Pi_a(-a): \bar\Pi_a)
$$
It yields a section $_{\bar\Pi}\rho_{b,a}: \GG Q\HH_{b,a}\to \wt\cG_{b,a}$ defined as in Section~\ref{Sect_3.3.2}. Namely, it sends $g$ to $(g,\cB)$ with
$$
\cB=\cJ_{\bar\Pi,b}\otimes\cJ_{\bar\Pi, a}^{-1}\otimes\det(\bar\Pi_b: g\bar\Pi_a)
$$
\sssec{} For the decomposition 
$
W_a=\bar\Upsilon_a\oplus \bar\Upsilon_a^*\otimes\Omega(a)
$
define a $\ZZ/2\ZZ$-graded line (purely of parity zero)
$$
\cJ_{\bar\Upsilon, a}=C_{a,x}^{-mna}\otimes \det (L_a(-a): L_a)^m
$$
Using Lemma~\ref{Lm_too_general_det_rel} we equip it with a
natural $\ZZ/2\ZZ$-graded  isomorphism
$$
\cJ_{\bar\Upsilon, a}^2\,\iso\, \cJ_a\otimes \det(\bar\Upsilon_a(-a): \bar\Upsilon_a) 
$$
It yields a section $_{\bar\Upsilon}\rho_{b,a}: \HH Q\GG_{b,a}\to \wt\cG_{b,a}$ defined as in Section~\ref{Sect_3.3.2}. Namely, it sends $g$ to $(g,\cB)$ with
$$
\cB=\cJ_{\bar\Upsilon, b}\otimes \cJ_{\bar\Upsilon, a}^{-1}\otimes \det(\bar\Upsilon_b: g\bar\Upsilon_a)
$$

The following is straightforward 
from definitions.

\begin{Lm} For $a,b\in\ZZ$ the following diagrams are canonically commutative
$$
\begin{array}{ccc}
\cT_{b,a} & \toup{\nu_{b,a}} & \wt\cG_{b,a}\\
\uparrow & \nearrow\lefteqn{\scriptstyle {_{\bar\Pi}\rho_{b,a}}}\\
\GG Q\HH_{b,a}
\end{array}
\hskip 5em
\begin{array}{ccc}
\cT_{b,a} & \toup{\nu_{b,a}} & \wt\cG_{b,a}\\
\uparrow & \nearrow\lefteqn{\scriptstyle {_{\bar\Upsilon}\rho_{b,a}}}\\
\HH Q\GG_{b,a}
\end{array}
\eqno{\square}
$$
Here $\nu_{b,a}$ was defined in Section~\ref{Sect_4.1.2}. 
\end{Lm}

\sssec{} For $a\in\ZZ$ we have the functors $\cF_{\bar\Upsilon_a(F)}: \D(\Upsilon_a(F))\to \D(\wt\cL_d(W_a(F)))$ and 
$\cF_{\bar\Pi_a(F)}: \D(\Pi_a(F))\to \D(\wt\cL_d(W_a(F)))$ defined in Proposition~\ref{Pp_functors_cF_U(F)_r}. Note that the diagram is canonically commutative
$$
\begin{array}{ccc}
\D(\Upsilon_a(F)) & \toup{\cF_{\bar\Upsilon_a(F)}} & \D(\wt\cL_d(W_a(F))) \\
\downarrow\lefteqn{\scriptstyle \zeta_a} & \nearrow\lefteqn{\scriptstyle \cF_{\bar\Pi_a(F)}}\\
\D(\Pi_a(F))
\end{array}
$$

\begin{Rem} 
\label{Rem_groupoid_action_two_Shrodinger_models}
For each $g\in \cT_{b,a}$ there are canonical equivalences $g_{\Pi}, g_{\Upsilon}$ that fit into a commutative diagram
\begin{equation}
\label{diag_for_Rem_5.5.8}
\begin{array}{ccc}
\D(\Pi_a(F)) & \toup{g_{\Pi}} & \D(\Pi_b(F))\\
\uparrow\lefteqn{\scriptstyle \zeta_a} && \uparrow\lefteqn{\scriptstyle \zeta_b}\\
\D(\Upsilon_a(F)) & \toup{g_{\Upsilon}} & \D(\Upsilon_b(F))
\end{array}
\end{equation}
with the following properties: 
\begin{itemize}
\item[i)] they are compatible with the groupoid structure on $\cT$; 
\item[ii)] if $g\in \GG Q\HH_{b,a}$ then $g_{\Pi}=g_*$ for the isomorphism $g: \Pi_a(F)\,\iso\,\Pi_b(F)$; 
\item[iii)] if $g\in\HH Q\GG_{b,a}$ then $g_{\Upsilon}=g_*$ for the isomorphism $g: \Upsilon_a(F)\,\iso\, \Upsilon_b(F)$.  
\end{itemize}
\end{Rem}
\begin{Prf}
For $g\in \GG Q\HH_{b,a}$ and for $g\in \HH Q\GG_{b,a}$ we define $g_{\Pi}$ and $g_{\Upsilon}$ by the properties ii), iii) and the commutativity of (\ref{diag_for_Rem_5.5.8}).  

Any $g\in\cT_{b,a}$ writes as a composition $g=g''\comp g'$ with 
$g''\in \HH Q\GG_{b,b}$ and $g'\in \GG Q\HH_{b,a}$. The compatibility with the groupoid structure on $\cT$ implies the uniqueness of $g_{\Pi}, g_{\Upsilon}$. The existence follows from the fact that the Fourier transform $\zeta_a$ commutes with the actions of $Q\GG\HH_a(F)=\GG Q\HH_a(F)\cap \HH Q\GG_a(F)$.
\end{Prf}  %%% ubrat???? Ne  nuzhna?  

\ssec{} As in (\cite{L2}, Section~6.2), one checks that we have the full subcategories (stable under subquotients) 
$$
\P_{\HH Q\GG_a(\cO)}(\Upsilon_a(F)) \subset \P_{Q\GG\HH_a(\cO)}(\Upsilon_a(F))\subset \P(\Upsilon_a(F))
$$
$$
\P_{\GG Q\HH_a(\cO)}(\Pi_a(F)) \subset \P_{Q\GG\HH_a(\cO)}(\Pi_a(F))\subset \P(\Pi_a(F)),
$$
and $\zeta_a$ yields an equivalence 
$$
\zeta_a: \P_{Q\GG\HH_a(\cO)}(\Upsilon_a(F))\,\iso\,
\P_{Q\GG\HH_a(\cO)}(\Pi_a(F))
$$ 
The above categories of equivariant perverse sheaves (and the corresponding $\DG$-categories) are defined as in (\cite{L2}, Section~4.2). 

\begin{Def} 
\label{Def_Weil_a}
For $a\in\ZZ$ let $\Weil_a$ be the category of triples $(\cF_1, \cF_2, \beta)$, where 
$$
\cF_1\in \P_{\HH Q\GG_a(\cO)}(\Upsilon_a(F)), \;\;\;
\cF_2\in \P_{\GG Q\HH_a(\cO)}(\Pi_a(F)),
$$
and $\beta: \zeta_a(\cF_1)\,\iso\, \cF_2$ is an isomorphism in $\P_{Q\GG\HH_a(\cO)}(\Pi_a(F))$. Write $\D\!\Weil_a$ for the category obtained by replacing everywhere in the above definition $\P$ by $\D\!\P$. 
\end{Def}

 So, any object of $\D\!\Weil_a$ is a direct sum of objects of the form $K[i]$ with $K\in\Weil_a$, $i\in\ZZ$. For $K,K'\in\Weil_a$ we have
$$
\Hom_{\D\!\Weil_a}(K[i], K'[j])=\left\{
\begin{array}{ll}
\Hom_{\Weil_a}(K,K'), & \mbox{for}\; i=j\\
0, & \mbox{for}\; i\ne j
\end{array}
\right.
$$
The category $\D\!\Weil_a$ is abelian $\Qlb$-linear. 
\sssec{}
Set 
$$
\cT_a(F)=\{(g_1, g_2)\in \GG_a(F)\times \HH_a(F)\mid g_1\otimes g_2\;\mbox{acts trivially on}\; A_a\otimes C_a(F)\},
$$ 
this is a group ind-scheme.

Clearly, $\Weil_a$ is an abelian category, and the forgetful functors $f_{\HH}: \Weil_a\to \P_{\HH Q\GG_a(\cO)}(\Upsilon_a(F))$ and $f_{\GG}:\Weil_a\to \P_{\GG Q\HH_a(\cO)}(\Pi_a(F))$ are fully faithful. As in (\cite{L2}, Section~6.2.1) one shows that the full embeddings $f_{\HH}, f_{\GG}$ are stable under subquotients. 

  The coproduct of $\Weil_a$ over $a\in\ZZ$ (in the category of abelian categories) plays a role of the category of $(\GG\times\HH)(\cO)$-invariants in the compactly induced representation 
$$
c-ind_{\cT_0(F)}^{(\GG\times\HH)(F)}\cS,
$$
where $\cS$ is the Weil representation of $\cT_0(F)$. 
Its definition is analogous to (\cite{L2}, Definition~4). 

\begin{Rem} Recall that in $\DGCat_{cont}$ the products are canonically equivalent to coproducts. The $\DG$-categories $Shv(\cT_a(F))$, $Shv(\GG(F)\times\HH(F))$ with their convolution monoidal structures are defined as in (\cite{G5}, Section~C.3.1). 

 We hope the structure appeared in Remark~\ref{Rem_groupoid_action_two_Shrodinger_models} could be extended as follows. One should be able to define an action\footnote{in the sense of (\cite{G5}, Section~1.4.2).} of $Shv(\GG(F)\times\HH(F))$ on $\prod_{a\in\ZZ} \D(\Pi_a(F))$, where the product is taken in $\DGCat_{cont}$. In particular, this would give an action of $Shv(\cT_a(F))$ on $\D(\Pi_a(F))$. The $\DG$-version of the category $\Weil_a$ would be the category of invariants $\D(\Pi_a(F))^{\cT_a}$, compare with Remark~\ref{Rem_4.1.2_DG_Weil_version}. We do not need this in the present paper.
\end{Rem}

\sssec{} 
By Proposition~\ref{Pp_functors_cF_U(F)_r}, we get a functor 
$$
\cF_{\Weil_a}: \Weil_a\to \P_{\cT_a}(\wt\cL_d(W_a(F)))
$$ 
sending $(\cF_1,\cF_2,\beta)$ to $\cF_{\bar\Upsilon_a(F)}(\cF_1)$.  

\sssec{} 
\label{Sect_5.6.3_about_I_0}
Let $I_0\in \P_{\HH Q\GG_0(\cO)}(\Upsilon_0(F))$ denote the constant perverse sheaf on $\Upsilon_0$ extended by zero to $\Upsilon_0(F)$. Remind that $\zeta_0(I_0)$ is the constant perverse sheaf on $\Pi_0$ extended by zero to $\Pi_0(F)$. The object $\zeta_0(I_0)$ will also be denoted $I_0$ by abuse of notation. So, $I_0\in\Weil_0$ naturally. By definition of the theta-sheaf, we have canonically in $\P_{\cT_0}(\wt\cL_d(W_0(F)))$
$$
\cF_{\Weil_0}(I_0)\,\iso\, S_{W_0(F)}
$$ 

\ssec{Hecke functors for the Shr\"odinger models}
\label{Sect_Hecke_functors_for_the_Shrodinger_models}

\sssec{} 
\label{Sect_4.7.1}
For $a\in\ZZ$ let $^a\cX\Pi$ be the stack classifying: a $\GL_m\times\Gm$-torsor $(U,C)$ over $\Spec\cO$, $\GG$-torsor $(M,A,\wedge^2 M\to A)$ over $\Spec\cO$, an isomorphism $A\otimes C\,\iso\, \Omega(a)$, and a section $s\in U^*\otimes M^*\otimes\Omega(F)$. 

 Informally, we may view $\D_{\GG Q\HH_a(\cO)}(\Pi_a(F))$ as the derived category on $^a\cX\Pi$. 
 For $a,a'\in \ZZ$ we are going to define Hecke functors
\begin{equation}
\label{def_Hecke_functor_for_Pi}
\begin{array}{l}
\H^{\la}_{\GG}: {_{a'-a}\Sph_{\GG}}\times \D_{\GG Q\HH_{a'}(\cO)}(\Pi_{a'}(F))\to \D_{\GG Q\HH_a(\cO)}(\Pi_a(F))\\  \\
\H^{\ra}_{\GG}: {_{a'-a}\Sph_{\GG}}\times \D_{\GG Q\HH_a(\cO)}(\Pi_a(F))\to \D_{\GG Q\HH_{a'}(\cO)}(\Pi_{a'}(F))
\end{array}
\end{equation}
To do so, consider the stack $^{a,a'}\cH_{\cX\Pi}$ classifying: a point of $^a\cX\Pi$ as above, a lattice $M'\subset M(F)$ such that for $A'=A(a'-a)$ the induced form $\wedge^2 M'\to A'$ is regular and nondegenerate. 

 We get a diagram
$$
^a\cX\Pi \; \getsup{h^{\la}}\; {^{a,a'}\cH_{\cX\Pi}}\; \toup{h^{\ra}}\; {^{a'}\cX\Pi},
$$ 
where $h^{\la}$ sends the above collection to $(U,C,M,A,s)$, the map $h^{\ra}$ sends the above collection to $(U,C,M', A',s')$, where $s'$ is the image of $s$ under $U^*\otimes M^*\otimes\Omega(F)\,\iso\, U^*\otimes M'^*\otimes\Omega(F)$. 

 Trivializing a point of $^{a'}\cX\Pi$ (resp., of $^a\cX\Pi$), one gets isomorphisms
$$
\id^r: {^{a,a'}\cH_{\cX\Pi}}\,\iso\, (\Pi_{a'}(F)\times \Gr_{\GG_{a'}}^{a-a'})/\GG Q\HH_{a'}(\cO)
$$
and
$$
\id^l: {^{a,a'}\cH_{\cX\Pi}}\,\iso\, (\Pi_a(F)\times \Gr_{\GG_{a}}^{a'-a})/\GG Q\HH_a(\cO)
$$
So for 
$$
K\in \D_{\GG Q\HH_a(\cO)}(\Pi_a(F)), \; K'\in \D_{\GG Q\HH_{a'}(\cO)}(\Pi_{a'}(F)), \; \cS\in {_{a'-a}\Sph_{\GG}}, \; \cS'\in {_{a-a'}\Sph_{\GG}}
$$ 
one can form the twisted exteriour products
$(K\tboxtimes\cS)^l$ and $(K'\tboxtimes\cS')^r$ 
on $^{a,a'}\cH_{\cX\Pi}$. The functors (\ref{def_Hecke_functor_for_Pi}) are defined by 
$$
\H^{\la}_{\GG}(\cS, K')=h^{\la}_! (K'\tboxtimes \ast\cS)^r,\;\;\;\; \H^{\ra}_{\GG}(\cS, K)=h^{\ra}_!(K\tboxtimes \cS)^l
$$
It is understood that this informal definition is made rigorous as in (\cite{L2}, Section~4.3).  

\sssec{} For $a\in\ZZ$ let $^a\cX\Upsilon$ be the stack classifying: a $\GL_n\times\Gm$-torsor $(L,A)$ over $\Spec\cO$, a $\HH$-torsor $(V,C)$ over $\cO$ (so, we are also given a compatible trivialization $\det V\,\iso\, C^m$), an isomorphism $A\otimes C\,\iso\,\Omega(a)$, and a section $s\in L^*\otimes V^*\otimes\Omega(F)$. 

 We may view $\D_{\HH Q\GG_a(\cO)}(\Upsilon_a(F))$ as the derived category on $^a\cX\Upsilon$. For $a,a'\in\ZZ$ we define Hecke functors
\begin{equation}
\label{def_Hecke_functor_for_Upsilon}
\begin{array}{l}
\H^{\la}_{\HH}: {_{a'-a}\Sph_{\HH}}\times \D_{\HH Q\GG_{a'}(\cO)}(\Upsilon_{a'}(F))\to \D_{\HH Q\GG_a(\cO)}(\Upsilon_a(F))\\ \\
\H^{\ra}_{\HH}: {_{a'-a}\Sph_{\HH}}\times \D_{\HH Q\GG_{a}(\cO)}(\Upsilon_{a}(F))\to \D_{\HH Q\GG_{a'}(\cO)}(\Upsilon_{a'}(F))
\end{array}
\end{equation}
as follows. Let $^{a,a'}\cH_{\cX\Upsilon}$ be the stack classifying: a point of $^a\cX\Upsilon$ as above, a lattice $V'\subset V(F)$ such that for $C'=C(a'-a)$ the induced form $\Sym^2 V'\to C'$ is regular and nondegenerate (we also get a compatible trivialization 
$$
C'^{-m}\otimes\det V'\,\iso\, C^{-m}\otimes\det V\,\iso\,\cO,
$$ 
so $(V',C')$ is a $\HH$-torsor over $\Spec\cO$). 
 
 As in Section~\ref{Sect_4.7.1}, we get the diagram 
$$
^a\cX\Upsilon\;\getsup{h^{\la}}\;
^{a,a'}\cH_{\cX\Upsilon} \;\toup{h^{\ra}} \;{^{a'}\cX\Upsilon}
$$
and the desired functors (\ref{def_Hecke_functor_for_Upsilon}). 

\sssec{} We need the following lemma. Write $^a\cX\bar\Pi$ for the stack classifying: a $\GL_m\times\Gm$-torsor $(U,C)$ over $\Spec\cO$, a $\GG$-torsor $(M,A, \wedge^2 M\to A)$ over $\Spec\cO$, an isomorphism $A\otimes C\,\iso\, \Omega(a)$, and a section $s_1\in U\otimes M(F)$. View $\D_{\GG Q\HH_a(\cO)}(\bar\Pi_a(F))$ as the derived category on $^a\cX\bar\Pi$. 

 For $a,a'\in\ZZ$ define the Hecke functor 
\begin{equation}
\label{def_Hecke_functor_for_check_Pi}
\H^{\la}_{\GG}: {_{a'-a}\Sph_{\GG}}\times\D_{\GG Q\HH_{a'}}(\bar\Pi_{a'}(F))\to \D_{\GG Q\HH_a}(\bar\Pi_a(F))
\end{equation}
as follows. Let $^{a,a'}\cH_{\cX\bar\Pi}$ be the stack classifying: a point of $^a\cX\bar\Pi$ as above, a lattice $M'\subset M(F)$ such that for $A'=A(a'-a)$ the induced form $\wedge^2 M'\to A'$ is regular and nondegenerate. 

 As above, we get a diagram
$$
^a\cX\bar\Pi \getsup{h^{\la}} {^{a,a'}\cH_{\cX\bar\Pi}} \toup{h^{\ra}} {^{a'}\cX\bar\Pi},
$$ 
where $h^{\la}$ sends the above point to $(U,C,M,A, s_1)$, the map $h^{\ra}$ sends the above point to $(U,C,M',A',s'_1)$, where $s'_1$ is the image of $s_1$ under $U\otimes M(F)\,\iso\, U\otimes M'(F)$. Now (\ref{def_Hecke_functor_for_check_Pi}) is defined in a way similar to (\ref{def_Hecke_functor_for_Pi}).

Write 
$$
\Four_{\psi}:\D_{\GG Q\HH_{a'}}(\Pi_{a'}(F))\,\iso\, \D_{\GG Q\HH_{a'}}(\bar\Pi_{a'}(F))
$$ 
for the Fourier transform (normalized as in Section~\ref{Sect_4.5.3}). The following is standard, cf. also (\cite{L2}, Lemma~11). 

\begin{Lm} 
\label{Lm_Fourier_Hecke_commute}
We have a canonical isomorphism in $\D_{\GG Q\HH_a}(\bar\Pi_a(F))$
$$
\H_{\GG}^{\la}(\cS, \Four_{\psi}(K))\,\iso\, \Four_{\psi}\H^{\la}_{\GG}(\cS, K)
$$
functorial in $\cS\in {_{a'-a}\Sph_{\GG}}$, $K\in \D_{\GG Q\HH_{a'}}(\Pi_{a'}(F))$. \QED
\end{Lm}

\sssec{Parabolic subgroups}
\label{Sect_Parabolic subgroups}
Write $P_{\HH_a}\subset \HH_a$ (resp., $P^-_{\HH_a}\subset \HH_a$) for the parabolic subgroup preserving $U_a$ (resp., $U_a^*\otimes C_a$). Let $U_{\HH_a}\subset P_{\HH_a}$ and $U^-_{\HH_a}\subset P^-_{\HH_a}$ denote their unipotent radicals. We view all of them as group schemes over $\Spec\cO$. Then $U_{\HH_a}\,\iso\, C_a^*\otimes \wedge^2 U_a$ and $U^-_{\HH_a}\,\iso\, C_a\otimes \wedge^2 U_a^*$ canonically. 

 Similarly, let $P_{\GG_a}\subset \GG_a$ (resp., $P^-_{\GG_a}\subset\GG_a$) be the parabolic subgroup preserving $L_a$ (resp., $L_a^*\otimes A_a$). Write $U_{\GG_a}\subset P_{\GG_a}$ and $U^-_{\GG_a}\subset P^-_{\GG_a}$ for their unipotent radicals. All of them are group schemes over $\Spec\cO$. We have canonically 
$$
U_{\GG_a}\,\iso\, A_a^*\otimes \Sym^2 L_a,\;\;\;\;\; U^-_{\GG_a}\,\iso\, A_a\otimes \Sym^2 L_a^*
$$ 
 
 View $v\in\Pi_a(F)$ as a map $v: C_a^*\otimes U_a(F)\to M_a(F)$. For $v\in \Pi_a(F)$ let $s_{\Pi}(v)$ denote the composition
$$
\wedge^2 (U_a\otimes C_a^{-1})(F)\toup{\wedge^2 v} \wedge^2 M_a(F)\to A_a(F) 
$$ 
Let $\Char(\Pi_a)\subset \Pi_a(F)$ denote the ind-subscheme of $v\in \Pi_a(F)$ such that $s_{\Pi}(v): C_a^*\otimes\wedge^2 U_a\to \Omega$ is regular.
An object $K\in \P(\Upsilon_a(F))$ is $U_{\HH_a}(\cO)$-equivariant iff $\zeta_a(K)$ is the extension by zero from $\Char(\Pi_a)$. This follows from the explicit formulas for the Shrodinger model of the Weil representations as in (\cite{L2}, Proposition~6). 

 View $v\in \Upsilon_a(F)$ as a map $v: L_a\otimes A_a^*(F)\to V_a(F)$. For $v\in \Upsilon_a(F)$ let $s_{\Upsilon}(v)$ denote the composition
$$
\Sym^2 (A_a^*\otimes L_a)\toup{\Sym^2 v} \Sym^2 V_a(F)\to C_a(F)
$$
Write $\Char(\Upsilon_a)\subset \Upsilon_a(F)$ for the ind-subscheme of $v\in \Upsilon_a(F)$ such that $s_{\Upsilon}(v): A_a^*\otimes \Sym^2 L_a\to \Omega$ is regular. An object $K\in \P(\Pi_a(F))$ is $U_{\GG_a}(\cO)$-equivariant iff $\zeta_a^{-1}(K)$ is the extension by zero from $\Char(\Upsilon_a)$. 

 The next result follows from (\cite{L2}, Lemma~8).
 
\begin{Lm}  
\label{Lm_4.7.6}
The full subcategory $\P_{\HH Q\GG_a(\cO)}(\Upsilon_a(F))\subset \P(\Upsilon_a(F))$ is the intersection of the full subcategories 
$$
\P_{U_{\HH_a}(\cO)}(\Upsilon_a(F))\;\cap\;
\P_{U^-_{\HH_a}(\cO)}(\Upsilon_a(F))\;\cap\; \P_{Q\GG\HH_a(\cO)}(\Upsilon_a(F))
$$ 
inside $\P(\Upsilon_a(F))$. \QED
\end{Lm}
 
\begin{Pp} 
\label{Pp_Hecke_functors_factor}
For $a\in\ZZ$ the functor $_{-a}\Sph_{\GG}\to \D_{\GG Q\HH_a(\cO)}(\Pi_a(F))$ sending $\cS$ to $\H^{\la}_{\GG}(\cS, I_0)$ factors naturally into 
$$
_{-a}\Sph_{\GG}\to\D\!\Weil_a\to \D_{\GG Q\HH_a(\cO)}(\Pi_a(F))
$$  
For $a\in\ZZ$ the functor $_{-a}\Sph_{\HH}\to \D_{\HH Q\GG_a(\cO)}(\Upsilon_a(F))$ sending $\cS$ to $\H^{\la}_{\HH}(\cS, I_0)$ factors naturally into 
$$
_{-a}\Sph_{\HH}\to\D\!\Weil_a\to \D_{\HH Q\GG_a(\cO)}(\Upsilon_a(F))
$$
\end{Pp} 
\begin{Prf}
The argument is similar for both claims, we prove only the first one. For a finite subfield $k'\subset k$ pick a $k'$-structure on $\cO$. Then $I_0$ admits a $k'$-structure and, as such, is pure of weight zero. So, by the decomposition theorem (\cite{BBD}), one has $\H^{\la}_{\GG}(\cS, I_0)\in \D\P^{ss}_{\GG Q\HH_a}(\Pi_a(F))$. 

 It remains to show that each perverse cohomology sheaf $K$ of $\zeta_a^{-1}\H^{\la}_{\GG}(\cS, I_0)$ lies in the full subcategory $\P_{\HH Q\GG_a(\cO)}(\Upsilon_a(F))$ of $\P_{Q\GG\HH_a(\cO)}(\Upsilon_a(F))$. 
 
 By definition of the Hecke functors, $\H^{\la}_{\GG}(\cS, I_0)$ is the extension by zero from $\Char(\Pi_a)$, so $\zeta_a(K)$ also satisfies this property. This yields a $U_{\HH_a}(\cO)$-action on $K$. 
 
 To get a $U^-_{\HH_a}(\cO)$-action on $K$, consider the commutative diagram of equivalences
$$
\begin{array}{ccc}
\P_{Q\GG\HH_a(\cO)}(\bar\Pi_a(F)) & \getsup{\zeta_{1,a}} &
\P_{Q\GG\HH_a(\cO)}(\Upsilon_a(F)) \\
\uparrow\lefteqn{\scriptstyle \Four_{\psi}} & \swarrow\lefteqn{\scriptstyle \zeta_a}\\
\P_{Q\GG\HH_a(\cO)}(\Pi_a(F)),
\end{array}
$$
where $\Four_{\psi}$ is the complete Fourier transform, and $\zeta_{1,a}$ is the corresponding partial one. 
 
  For $v\in\bar\Pi_a(F)$ write $s_{\bar\Pi}(v)$ for the composition 
$$
\wedge^2 U_a^*(F)\toup{\wedge^2 v} \wedge^2 M_a(F)\to A_a(F)
$$
Write $\Char(\bar\Pi_a)\subset \bar\Pi_a(F)$ for the ind-subscheme of $v$ such that $s_{\bar\Pi}(v): C_a\otimes\wedge^2 U_a^*\to \Omega$ is regular. The $U^-_{\HH_a}(\cO)$-equivariance of $K$ is equivalent to the fact that $\zeta_{1,a}(K)$ is the extension by zero from $\Char(\bar\Pi_a)$. 

 By Lemma~\ref{Lm_Fourier_Hecke_commute}, we have $\Four_{\psi}\H^{\la}_{\GG}(\cS, I_0)\,\iso\, 
\H^{\la}_{\GG}(\cS, \check{I}_0)$, where $\check{I}_0:=\Four_{\psi}(I_0)$ is the constant perverse sheaf on $\bar\Pi_0$ extended by zero to $\bar\Pi_0(F)$. Clearly, $\H^{\la}_{\GG}(\cS, \check{I}_0)$ is the extension by zero from $\Char(\bar\Pi_a)$, and our assertion follows from Lemma~\ref{Lm_4.7.6}. 
\end{Prf} 

\begin{Rem} As in (\cite{L2}, Remark~9) Proposition~\ref{Pp_Hecke_functors_factor} could be extended as follows. For $a,b\in\ZZ$ there is a natural functor $_{-a}\Sph_{\GG}\times\D\!\Weil_b\to \D\!\Weil_{a+b}$ given informally by $(\cS, K)\mapsto \H^{\la}_{\GG}(\cS, K)$. Besides, there is a natural functor 
$$
_{-a}\Sph_{\HH}\times\D\!\Weil_b\to \D\!\Weil_{a+b}
$$ 
given informally by $(\cS, K)\mapsto \H^{\la}_{\HH}(\cS, K)$.
\end{Rem}

\sssec{} According to Proposition~\ref{Pp_Hecke_functors_factor}, in what follows we  write $\H^{\la}_{\GG}(\cdot, I_0): {_{-a}\Sph_{\GG}}\to \D\Weil_a$ and $\H^{\la}_{\HH}(\cdot, I_0): {_{-a}\Sph_{\HH}}\to \D\Weil_a$ for the corresponding functors. 

 From Proposition~\ref{Pp_functors_cF_U(F)_r} one derives the following. 

\begin{Cor} For $a\in\ZZ$, $\cS\in{_{-a}\Sph_{\GG}}$, $\cT\in{_{-a}\Sph_{\HH}}$ there are canonical isomorphisms
$$
\cF_{\Weil_a}\H^{\la}_{\GG}(\cS, I_0)\,\iso\, \H^{\la}_{\GG}(\cS, S_{W_0(F)})
$$ 
and 
$$
\cF_{\Weil_a}\H^{\la}_{\HH}(\cT, I_0)\,\iso\, \H^{\la}_{\HH}(\cT, S_{W_0(F)})
$$ 
in $\D\P^{ss}_{\cT_a}(\wt\cL_d(W_a(F)))$. They compatible with the tensor structures on $\Sph_{\GG}, \Sph_{\HH}$. 
\QED
\end{Cor}

 Here $\D\P^{ss}_{\cT_a}(\wt\cL_d(W_a(F)))$ is an abelian category defined as in Section~\ref{Sect_2.2.1_revised}, namely the abelian category of cohomologically graded objects of $\P^{ss}_{\cT_a}(\wt\cL_d(W_a(F)))$. 

Thus, Theorem~\ref{Th_local_main} is reduced to the following. 
\begin{Th} 
\label{Th2_local} Let the maps $\kappa$ be as in Theorem~\ref{Th_local_main}, $a\in\ZZ$.\\ 
1) Assume $m\le n$. The two functors $_{-a}\Rep(\check{\GG})\to \D\Weil_a$ given by
$$
V\mapsto \H^{\la}_{\GG}(V, I_0)\;\;\;\mbox{and}\;\;\;
V\mapsto \H^{\ra}_{\HH}(\Res^{\kappa}(V), I_0)
$$
are canonically isomorphic. This isomorphism is compatible with the tensor structures on $\Rep(\check{\HH})$, $\Rep(\check{\GG})$. \\
2) Assume $m>n$. The two functors $_{-a}\Rep(\check{\HH})\to \D\Weil_a$ given by
$$
V\mapsto \H^{\la}_{\HH}(V, I_0)\;\;\;\mbox{and}\;\;\;
V\mapsto \H^{\ra}_{\GG}(\Res^{\kappa}(V), I_0)
$$
are canonically isomorphic. This isomorphism is compatible with the tensor structures on $\Rep(\check{\HH})$, $\Rep(\check{\GG})$.
\end{Th}  

\begin{Rem} For $a=0$ Theorem~\ref{Th2_local} is nothing but (\cite{L2}, Theorem~7). The case of $a$ even reduces to the case $a=0$ in view of Remark~\ref{Rem_Weil_a_depends_on_parity}.
\end{Rem} 

\ssec{Hecke functors for Levi subgroups}

\sssec{}
\label{Sect_4.8.1}
For $a\in\ZZ$ set $Q\Pi_a=U_a^*\otimes C_a\otimes L_a\subset \Pi_a$ and $Q\Upsilon_a=L_a^*\otimes A_a\otimes U_a\subset \Upsilon_a$.

We are going to define for $a,a'\in\ZZ$ Hecke functors
\begin{equation}
\label{Hecke_functor_for_Q(GG_a)}
\H^{\la}_{Q(\GG)}: {_{a'-a}\Sph_{Q(\GG)}}\times \D_{Q\GG\HH_{a'}(\cO)}(Q\Pi_{a'}(F))\to \D_{Q\GG\HH_a(\cO)}(Q\Pi_a(F))
\end{equation}
in a way compatible with the Hecke functors defined in Section~\ref{Sect_Hecke_functors_for_the_Shrodinger_models}. 

 Let $^a\cX Q\Pi$ be the stack classifying a $Q(\HH)$-torsor $(U,C)$ over $\Spec\cO$, a $Q(\GG)$-torsor $(L,A)$ over $\Spec\cO$, an isomorphism $A\otimes C\,\iso\, \Omega(a)$, and a section $s\in U^*\otimes C\otimes L(F)$. 
We think of $\D_{Q\GG\HH_a(\cO)}(Q\Pi_a(F))$ as the derived category on $^a\cX Q\Pi$. 

 Let $^{a,a'}\cH_{\cX Q\Pi, Q(\GG)}$ be the stack classifying: a point of $^a\cX Q\Pi$ as above, a $\cO$-lattice $L'\subset L(F)$, for which we set $A'=A(a'-a)$. We get a diagram
$$
^a\cX Q\Pi\; \getsup{h^{\la}} \; {^{a,a'}\cH_{\cX Q\Pi, Q(\GG)}}\; \toup{h^{\ra}} \; {^{a'}\cX Q\Pi},
$$
where $h^{\la}$ sends the above collection to $(U,C, L,A,s)$, the map $h^{\ra}$ sends the above collection to
$(U,C, L', A', s')$, where $s'$ is the image of $s$ under
$U^*\otimes C\otimes L(F)\,\iso\, U^*\otimes C\otimes L'(F)$. 

 Trivializing a point of $^{a'}\cX Q\Pi$ (resp., of $^a\cX Q\Pi$), one gets isomorphisms
$$
\id^r: {^{a,a'}\cH_{\cX Q\Pi, Q(\GG)}}\,\iso\, (Q\Pi_{a'}(F) \times \Gr_{Q(\GG_{a'})}^{a-a'})/Q\GG\HH_{a'}(\cO)
$$
and
$$
\id^l: {^{a,a'}\cH_{\cX Q\Pi, Q(\GG)}}\,\iso\, (Q\Pi_a(F) \times \Gr_{Q(\GG_a)}^{a'-a})/Q\GG\HH_a(\cO)
$$
So, for 
$$
K\in \D_{Q\GG\HH_a(\cO)}(Q\Pi_a(F)), \; K'\in \D_{Q\GG\HH_{a'}(\cO)}(Q\Pi_{a'}(F)), \; \cS\!\in{_{a'-a}\Sph_{Q(\GG)}}, \; \cS'\!\in {_{a-a'}\Sph_{Q(\GG)}}
$$ 
one can form their twisted exteriour products $(K\tboxtimes \cS)^l$ and $(K'\tboxtimes\cS')^r$ on $^{a,a'}\cH_{\cX Q\Pi, Q(\GG)}$. The functor (\ref{Hecke_functor_for_Q(GG_a)}) is defined by
$$
\H^{\la}_{Q(\GG)}(\cS', K')= h^{\la}_!(K'\tboxtimes \ast \cS')^r
$$
\sssec{} Let $^a\cX Q\Upsilon$ be the stack classifying a $Q(\HH)$-torsor $(U,C)$ over $\Spec\cO$, a $Q(\GG)$-torsor $(L,A)$ over $\Spec\cO$, an isomorphism $A\otimes C\,\iso\,\Omega(a)$, and a section $s\in U\otimes A\otimes L^*(F)$. Informally, we think of $\D_{Q\GG\HH_a(\cO)}(Q\Upsilon_a(F))$ as the derived category on $^a\cX Q\Upsilon$. One defines the Hecke functor
\begin{equation}
\label{Hecke_functor_for_Q(HH_a)}
\H^{\la}_{Q(\HH)}: {_{a'-a}\Sph_{Q(\HH)}}\times  \D_{Q\GG\HH_{a'}(\cO)}(Q\Upsilon_{a'}(F))\to \D_{Q\GG\HH_a(\cO)}(Q\Upsilon_a(F))
\end{equation}
using a similar diagram 
$$
^a\cX Q\Upsilon \,\getsup{h^{\la}}\, {^{a,a'}\cH_{\cX Q\Upsilon, Q(\HH)}}\, \toup{h^{\ra}}\, {^{a'}\cX Q\Upsilon}
$$

By abuse of notation, we also write $I_0$ for the constant perverse sheaf on $Q\Upsilon_0$ and on 
$Q\Pi_0$, the exact meaning is easily understood from the context. The next result is a straightforward consequence of (\cite{L2}, Corollary~4). 
 
\begin{Pp} 
\label{Pp_description_modules_for_Levi}
1) Assume $m>n$. The functor 
$$
_{-a}\Sph_{Q(\GG)}\to  \D_{Q\GG\HH_a(\cO)}(Q\Pi_a(F))
$$ 
given by $\cS\mapsto \H^{\la}_{Q(\GG)}(\cS, I_0)$ takes values in $\P^{ss}_{Q\GG\HH_a(\cO)}(Q\Pi_a(F))$ and induces an equivalence
$$
_{-a}\Sph_{Q(\GG)}\,\iso\; \P^{ss}_{Q\GG\HH_a(\cO)}(Q\Pi_a(F))
$$
2) Assume $m\le n$. The functor
$$
_{-a}\Sph_{Q(\HH)}\to \D_{Q\GG\HH_a(\cO)}(Q\Upsilon_a(F))
$$
given by $\cS\mapsto \H^{\la}_{Q(\HH)}(\cS, I_0)$ takes values in $\P^{ss}_{Q\GG\HH_a(\cO)}(Q\Upsilon_a(F))$ and induces an equivalence
$$
_{-a}\Sph_{Q(\HH)}\,\iso\, 
\P^{ss}_{Q\GG\HH_a(\cO)}(Q\Upsilon_a(F))
$$
\QED
\end{Pp}

\sssec{} 
\label{Sect_4.8.4}
For $a,a'\in\ZZ$ we will use in Section~\ref{Sect_Proof of Th2_local} the following Hecke functor
\begin{equation}
\label{Hecke_last_not_least_Q(G)}
\H^{\la}_{Q(\GG)}: {_{a'-a}\Sph_{Q(\GG)}}\times \D_{Q\GG\HH_{a'}(\cO)}(Q\Upsilon_{a'}(F))\to \D_{Q\GG\HH_a(\cO)}(Q\Upsilon_a(F))
\end{equation}
%%% zdes' napisat' analog of (\cite{L2}, Proposition~4) dlja operatorov na Q\Upsilon_a(F). 

 Consider the stack $^{a,a'}\cH_{\cX Q\Upsilon, Q(\GG)}$ classifying: a point $(U,C,L,A,s)$ of $^a\cX Q\Upsilon$ as above, a $\cO$-lattice $L'\subset L(F)$ for which we set $A'=A(a'-a)$. We get a diagram
$$
^a\cX Q\Upsilon \getsup{h^{\la}} {^{a,a'}\cH_{\cX Q\Upsilon, Q(\GG)}} \toup{h^{\ra}} {^{a'}\cX Q\Upsilon},
$$
where $h^{\la}$ sends the above collection to $(U,C,L,A,s)$, and $h^{\ra}$ sends the same collection to $(U,C, L',A', s')$, where $s'$ is the image of $s$ under 
$U\otimes A\otimes L^*(F)\,\iso\, U\otimes A'\otimes L'^*(F)$. The functor (\ref{Hecke_last_not_least_Q(G)}) is defined as in Section~\ref{Sect_4.8.1} for the above diagram. 

\begin{Lm} For $\cS\in {_{a'-a}\Sph_{Q(\GG)}}$ the diagram of functors is canonically commutative
$$
\begin{array}{ccc}
\D_{Q\GG\HH_{a'}(\cO)}(Q\Upsilon_{a'}(F)) & \toup{\Four_{\psi}} & \D_{Q\GG\HH_{a'}(\cO)}(Q\Pi_{a'}(F))\\
\downarrow\lefteqn{\scriptstyle \H^{\la}_{Q(\GG)}(\cS,\cdot)} && \downarrow\lefteqn{\scriptstyle \H^{\la}_{Q(\GG)}(\cS,\cdot)}\\
\D_{Q\GG\HH_a(\cO)}(Q\Upsilon_a(F)) & \toup{\Four_{\psi}} & \D_{Q\GG\HH_a(\cO)}(Q\Pi_a(F))
\end{array}
$$
\end{Lm} 
\begin{proof} This is an immediate consequence of (\cite{L2}, Lemma~9).
\end{proof}

\ssec{Weak Jacquet functors}
\label{Sect_Weak_Jacquet}

\sssec{}
As in (\cite{L2}, Section~4.7) for $a\in\ZZ$ we define the weak Jacquet functors 
\begin{equation}
\label{functor_J_P_HH_a}
J^*_{P_{\HH_a}}, J^!_{P_{\HH_a}}: \D_{\HH Q\GG_a(\cO)}(\Upsilon_a(F))\to \D_{Q\GG\HH_a(\cO)}(Q\Upsilon_a(F))
\end{equation}
and 
\begin{equation}
\label{functor_J_P_GG_a}
J^*_{P_{\GG_a}}, J^!_{P_{\GG_a}}: \D_{\GG Q\HH_a(\cO)}(\Pi_a(F))\to \D_{Q\GG\HH_a(\cO)}(Q\Pi_a(F))
\end{equation}
Both definitions being similar, we recall the definition of (\ref{functor_J_P_HH_a}) only. 

 For a free $\cO$-module of finite type $M$ and $N,r\in\ZZ$ with $N+r\ge 0$ write $_{N,r}M=M(N)/M(-r)$.

 For $N+r\ge 0$ consider the natural embedding $i_{N,r}: {_{N,r}Q\Upsilon_a}\hook{} {_{N,r} \Upsilon_a}$. Set 
$$
P Q\GG_a=\{g=(g_1,g_2)\in P_{\HH_a}\times Q(\GG_a)\mid g\in\cT_a\},
$$
this is a group scheme over $\Spec\cO$.
We have a diagram of stack quotients
$$
\begin{array}{ccccc}
P Q\GG_a(\cO/t^{N+r})\backslash (_{N,r} Q\Upsilon_a) & \toup{i_{N,r}} & 
P Q\GG_a(\cO/t^{N+r})\backslash (_{N,r} \Upsilon_a) & \toup{p} & \HH Q\GG_a(\cO/t^{N+r})\backslash (_{N,r} \Upsilon_a)\\
\downarrow\lefteqn{\scriptstyle q}\\
Q\GG\HH_a(\cO/t^{N+r})\backslash (_{N,r} Q\Upsilon_a),
\end{array}
$$
where $t\in \cO$ is a uniformizer, $p$ comes from the inclusion $P_{\HH_a}\subset \HH_a$, and $q$ is the natural quotient map. First, define functors
\begin{equation}
\label{functors_J_P^*_and_J_P^!__for_Nr}
J^*_{P_{\HH_a}}, J^!_{P_{\HH_a}}: \D_{\HH Q\GG_a(\cO/t^{N+r})}(_{N,r}\Upsilon_a)\to \D_{Q\GG\HH_a(\cO/t^{N+r})}(_{N,r}Q\Upsilon_a)
\end{equation}
by
$$
\begin{array}{c}
q^*\comp J^*_{P_{\HH_a}}[\dimrel(q)]= i_{N,r}^* p^*[\dimrel(p) - rnm]\\
\\
q^*\comp J^!_{P_{\HH_a}}[\dimrel(q)]= i_{N,r}^! p^*[\dimrel(p) + rnm]
\end{array}
$$
Since 
$$
q^*[\dimrel(q)]: \D_{Q\GG\HH_a(\cO/t^{N+r})}(_{N,r}Q\Upsilon_a)\to \D_{PQ\GG_a(\cO/t^{N+r})}(_{N,r}Q\Upsilon_a)
$$ 
is an equivalence (exact for the perverse t-structures), the functors (\ref{functors_J_P^*_and_J_P^!__for_Nr}) are well-defined. Further, (\ref{functors_J_P^*_and_J_P^!__for_Nr}) are compatible with the transition functors in the definition of the corresponding $\DG$-categories, so give rise to the functors (\ref{functor_J_P_HH_a}) in the colimit as $N,r$ go to infinity. Note that for (\ref{functor_J_P_HH_a}) we get $\DD\comp J^*_{P_{\HH_a}}\,\iso\, J^!_{P_{\HH_a}}\comp \DD$ naturally on the subcategories of constructible objects. 

%  Denote by $_a\D\Sph_{\HH}\subset \D\Sph_{\HH}$ the full subcategory of objects of the form $\oplus_{i\in\ZZ} \, K_i[i]$ with $K_i\in {_a\Sph_{\HH}}$. Define $_a\D\Sph_{\GG}$ similarly. 
  
\sssec{}  
\label{Sect_4.9.2}
We identify $\HH\,\iso\, \HH_0$ and $Q(\HH)\,\iso\, Q(\HH_0)$. Let $\check{\mu}_{\HH}=\det V_0=m\check{\alpha}_0$ and $\check{\nu}_{\HH}=\det U_0$ viewed as characters of $Q(\HH)$ or, equivalently, as cocharacters of the center $Z(\check{Q}(\HH))$ of the Langlands dual group $\check{Q}(\HH)$ of $Q(\HH)$. Let $\kappa_{\HH}: \check{Q}(\HH)\times\Gm\to\check{\HH}$ be the map, whose first component is the natural inclusion of the Levi subgroup, and the second one is $2(\check{\rho}_{\HH}-\check{\rho}_{Q(\HH)})+n(\check{\mu}_{\HH}-\check{\nu}_{\HH})$. 

Recall the definition of the geometric restriction functor $\gRes^{\kappa_{\HH}}$ given in Section~\ref{Sect_Geometric restriction functors}. In view of the Satake equivalences (\ref{def_Loc_gr}), it identifies with the restriction functor $\Res^{\kappa_{\HH}}: \Rep(\check{\HH})\to \Rep(\check{Q}(\HH)\times\Gm)$ with respect to $\kappa_{\HH}$. 
   
\begin{Lm} 
\label{Lm_Jacquet_working_HH}
For $a,a'\in\ZZ$ and $\cS\in {_{a'-a}\Sph_{\HH}}$, $K\in \D_{\HH Q\GG_{a'}(\cO)}(\Upsilon_{a'}(F))$ there is a filtration in $\D_{Q\GG\HH_a(\cO)}(Q\Upsilon_a(F))$ on 
$$
J^*_{P_{\HH_a}}\H^{\la}_{\HH}(\cS, K)
$$ 
such that the corresponding graded object identifies canonically with
$$
\H^{\la}_{Q(\HH)}(\gRes^{\kappa_{\HH}}(\cS), J^*_{P_{\HH_{a'}}}(K))
$$
\end{Lm} 
\begin{Prf}
For the notion of a filtration in a stable $\infty$-category we refer to (\cite{HA}, Definition~1.2.2.2). Let $I_{a'}$ be the constant perverse sheaf on $\Upsilon_{a'}$ extended by zero to $\Upsilon_{a'}(F)$. The proof is similar to (\cite{L2}, Lemma~6), we only have to determine the corresponding map $\check{Q}(\HH)\times\Gm\to\check{\HH}$. To do so, it suffices to perform the calculation for $K=I_{a'}$.

 For $s_1,s_2\ge 0$ let $_{s_1,s_2}\Gr_{\HH_a}\subset \Gr_{\HH_a}$ be the closed subscheme of $h\HH_a(\cO)\in\Gr_{\HH_a}$ such that 
$$
V_a(-s_1)\subset hV_a\subset V_a(s_2)
$$ 
Assume that $s_1,s_2$ are large enough so that $\cS$ is the extension by zero from $_{s_1,s_2}\Gr_{\HH_a}$. Then $\H^{\la}_{\HH}(\cS, I_{a'})\in \D_{\HH Q\GG_a(\cO)}(_{s_2,s_1}\Upsilon_a)$ is as follows.
 Write $_{0,s_1}\Upsilon\ttimes {_{s_1,s_2}\Gr_{\HH_a}}$ for the scheme classifying pairs 
$$
h\HH_a(\cO)\in {_{s_1,s_2}\Gr_{\HH_a}}, \; v\in L^*_a\otimes A_a\otimes (hV_a)/V_a(-s_1)
$$
Let $\pi: {_{0,s_1}\Upsilon\ttimes {_{s_1,s_2}\Gr_{\HH_a}}}\to {_{s_2,s_1}\Upsilon_a}$ be the map sending $(h\HH_a(\cO), v)$ to $v$. By definition,  
\begin{equation}
\label{complex_Hecke_operator_on_Upsilon_a}
\H^{\la}_{\HH}(\cS, I_{a'})\,\iso\, \pi_!(\Qlb\tboxtimes \cS),
\end{equation}
where $\Qlb\tboxtimes \cS$ is normalized to be perverse. If $\theta\in\pi_1(\HH)$ then 
 $_{0,s_1}\Upsilon\ttimes {_{s_1,s_2}\Gr_{\HH_a}}$ is a vector bundle over $_{s_1,s_2}\Gr_{\HH_a}^{\theta}$ of rank $2s_1nm-\<\theta, n\check{\mu}_{\HH}\>$. 

 Let $_{s_1,s_2}P_{\HH_a}=\{p\in P_{\HH_a}(F)\mid V_a(-s_1)\subset p V_a\subset V_a(s_2)\}$. Then 
$$
_{s_1,s_2}\Gr_{P_{\HH_a}}=(_{s_1,s_2}P_{\HH_a}(F))/P_{\HH_a}(\cO)
$$ 
is closed in $\Gr_{\PP_{\HH_a}}$. The natural map $_{s_1,s_2}\Gr_{P_{\HH_a}}\to {_{s_1,s_2}\Gr_{\HH_a}}$ at the level of reduced schemes yields a stratification of $_{s_1,s_2}\Gr_{\HH_a}$ by the connected components of $_{s_1,s_2}\Gr_{P_{\HH_a}}$. Calculate (\ref{complex_Hecke_operator_on_Upsilon_a}) with respect to this stratification. Denote by $_{s_1,s_2}\Gr_{Q(\HH_a)}\subset \Gr_{Q(\HH_a)}$ the closed subscheme of $hQ(\HH_a)\in\Gr_{Q(\HH_a)}$ satisfying 
$$
U_a(-s_1)\subset hU_a\subset U_a(s_2),
$$
write $\gt_P: \Gr_{P_{\HH_a}}\to \Gr_{Q(\HH_a)}$ for the natural map. We have the diagram
$$
\begin{array}{ccccccc}
_{0,s_1} Q\Upsilon\ttimes {_{s_1,s_2}\Gr_{Q(\HH_a)}} & \getsup{\id\times\gt_P} & {_{0,s_1} Q\Upsilon\ttimes {_{s_1,s_2}\Gr_{P_{\HH_a}}}} & \hook{} &  {_{0,s_1} \Upsilon\ttimes {_{s_1,s_2}\Gr_{P_{\HH_a}}}} & \to &
{_{0,s_1} \Upsilon\ttimes {_{s_1,s_2}\Gr_{\HH_a}}}
\\
 & \searrow\lefteqn{\scriptstyle \pi_Q} & \downarrow && \downarrow & \swarrow\lefteqn{\scriptstyle \pi}\\
 && _{s_2,s_1}Q\Upsilon_a & \hook{} &_{s_2,s_1}\Upsilon_a,
\end{array}
$$
where the square is cartesian.
Here $_{0,s_1} \Upsilon\ttimes {_{s_1,s_2}\Gr_{P_{\HH_a}}}$ is the scheme classifying pairs
$$
hP_{\HH_a(\cO)}\in {_{s_1,s_2}\Gr_{P_{\HH_a}}}, v\in L^*_a\otimes A_a\otimes (hV_a)/V_a(-s_1),
$$
and 
$
{_{0,s_1} Q\Upsilon\ttimes {_{s_1,s_2}\Gr_{P_{\HH_a}}}}
$ 
is its closed subscheme given $v\in L^*_a\otimes A_a\otimes (hU_a)/U_a(-s_1)$. 

 By definition, given $\cT\in {_{a'-a}\Sph_{Q(\HH)}}$ we have 
$$
\H^{\la}_{Q(\HH)}(\cT, I_{a'})\,\iso\, \pi_{Q !}(\Qlb\tboxtimes \cT),
$$ 
where $\Qlb\tboxtimes \cT$ is normalized to be perverse. 
If $\theta\in\pi_1(Q(\HH))$ then $_{0,s_1} Q\Upsilon\ttimes {_{s_1,s_2}\Gr_{Q(\HH_a)}}$ is a vector bundle over $_{s_1,s_2}\Gr_{Q(\HH_a)}^\theta$ of rank $s_1nm-\<\theta, n\check{\nu}_{\HH}\>$. Our assertion follows.
\end{Prf}
 
\sssec{} 
\label{Sect_4.9.4}
We identify $\GG\,\iso\, \GG_0$, $Q(\GG)\,\iso\, Q(\GG_0)$. Write $\check{\mu}_{\GG}=\det M_0$ and $\check{\nu}_{\GG}=\det L_0$ as cocharacters of the center $Z(\check{Q}(\GG))$ of the Langlands dual group $\check{Q}(\GG)$ of $Q(\GG)$. Let $\kappa_{\GG}: \check{Q}(\GG)\times\Gm\to\check{\GG}$ be the map whose first component is the natural inclusion of the Levi subgroup, and the second one is $2(\check{\rho}_{\GG}-\check{\rho}_{Q(\GG)})+m(\check{\mu}_{\GG}-\check{\nu}_{\GG})$. The corresponding geometric restriction functor is denoted $\gRes^{\kappa_{\GG}}$. 

\begin{Lm} 
\label{Lm_Jacquet_working_GG}
For $a,a'\in\ZZ$, $\cS\in {_{a'-a}\Sph_{\GG}}$, and $K\in \D_{\GG Q\HH_{a'}(\cO)}(\Pi_{a'}(F))$ there is a filtration in $\D_{Q\GG\HH_a(\cO)}(Q\Pi_a(F))$ on 
$$
J^*_{P_{\GG_a}}\H^{\la}_{\GG}(\cS, K)
$$
such that the corresponding graded object identifies canonically with $$
\H^{\la}_{Q(\GG)}(\gRes^{\kappa_{\GG}}(\cS), J^*_{P_{\GG_{a'}}}(K))$$ 
\end{Lm}
\begin{proof}
Similar to Lemma~\ref{Lm_Jacquet_working_HH}.
\end{proof}

 We use Lemmas~\ref{Lm_Jacquet_working_HH} and \ref{Lm_Jacquet_working_GG} in the following form. 
\begin{Cor}  
\label{Cor_Braden_rules}
Let $a,a'\in\ZZ$, $\cS\in {_{a'-a}\Sph_{\HH}}$, and $k_0\subset k$ is a finite subfield. Assume $K\in \P_{\HH Q\GG_{a'}(\cO)}(\Upsilon_{a'}(F))$ admits a $k_0$-structure and, as such, is pure of weight zero. Then $J^*_{P_{\HH_{a'}}}(K)$ is also pure of weight zero over $k_0$, the filtration on it constructed in Lemma~\ref{Lm_Jacquet_working_HH} splits canonically. So,  
there is a canonical isomorphism in $\D_{Q\GG\HH_a(\cO)}(Q\Upsilon_a(F))$
$$
J^*_{P_{\HH_a}}\H^{\la}_{\HH}(\cS, K)\,\iso\, \H^{\la}_{Q(\HH)}(\gRes^{\kappa_{\HH}}(\cS), J^*_{P_{\HH_{a'}}}(K))
$$
(Similar strengthened version of Lemma~\ref{Lm_Jacquet_working_GG} also holds.)
\end{Cor}
\begin{proof} First, a non-canonical splitting of the filtration is constructed precisely as in (\cite{L2}, Corollary~3).

 Let now $\bar q_{Q\Upsilon_a(F), *}: \D_{Q\GG\HH_a(\cO)}(Q\Upsilon_a(F))\to\D(\Spec k)$ be the functor defined as in (\cite{L2}, Section~4.5.4). It factors as the composition
$$
\D_{Q\GG\HH_a(\cO)}(Q\Upsilon_a(F))\to \D_{Q\GG\HH_a(\cO)}(\Spec k)\to \D(\Spec k),
$$ 
where the first functor is the cohomlogy and the second is the pullback. From (\cite{L2}, Corollary~4) we see that for any irreducible object $\cK$ of $\P_{Q\GG\HH_a(\cO)}(Q\Upsilon_a(F))$ one has $\bar q_{Q\Upsilon_a(F), *}(\cK)\ne 0$. Now as in (\cite{L2}, Lemma~8) we see that there is a canonical splitting of this filtration, because there is one after applying the functor $\bar q_{Q\Upsilon_a(F), *}$.
\end{proof}

\ssec{Action of $\Sph_{\GG}$}
\label{Sect_Action_of_Sph_GG}

\sssec{}
Recall our choices of the maximal torus and a Borel subgroup $T_{\GG}\subset B_{\GG}\subset \GG$, and similarly for $\HH$
(cf. Sections~\ref{Sect_2.3.3}, \ref{Sect_2.3.4}, \ref{Sect_Root data}). A trivialization of the $\GG_a$-torsor $(M_a, A_a)$ over $\Spec\cO$ yields a maximal torus and a Borel subgroup in $\GG_a$, an equivalence $\Sph_{\GG_a}\,\iso\, \Sph_{\GG}$ and a bijection $\Lambda^+_{\GG_a}\,\iso\, \Lambda^+_{\GG}$ as in Section~\ref{Sect_4.3.7} (and similarly for $\HH_a$). 
Write $w_0^{\GG}$ for the longest element of the Weyl group of $\GG$. 

 Write $\check{\omega}_i$ for the highest weight of the fundamental representation of $\GG_a$ that appear in $\wedge^i M_a$ for $i=1,\ldots,n$. All the weights of $\wedge^i M_a$ are $\le \check{\omega}_i$. As above, $\check{\omega}_0$ denotes the weight of the $\GG_a$-module $A_a$.  
 
 For $\lambda\in\Lambda^+_{\GG}$ set $a=\<\lambda,\check{\omega}_0\>$ then $\cA^{\lambda}_{\GG}\in {_{-a}\Sph_{\GG}}$. By definition, 
$$
\H^{\lambda}_{\GG}(I_0)=\H^{\la}_{\GG}(\cA^{\lambda}, I_0)\in \D_{\GG Q\HH_a}(\Pi_a(F))
$$ 
is as follows. Set $r=\<\lambda,\check{\omega}_1\>$ and $N=\<-w_0^{\GG}(\lambda),\check{\omega}_1\>$. 
Let $_{0,r}\Pi\ttimes \Grb_{\GG_a}^{\lambda}$ be the scheme classifying $g\in\Grb^{\lambda}_{\GG_a}$, $x\in U_a^*\otimes C_a\otimes ((gM_a)/M_a(-r))$. Let \begin{equation}
\label{map_pi_for_Hecke_GG}
\pi: {_{0,r}\Pi}\ttimes \Grb_{\GG_a}^{\lambda}\to 
{_{N,r}\Pi_a}
\end{equation}
 be the map sending $(x, g\GG_a(\cO))$ to $x$. Then $\H^{\lambda}_{\GG}(I_0)\,\iso\, \pi_!(\Qlb\tboxtimes \cA^{\lambda}_{\GG})$ canonically. Here $\Qlb\tboxtimes \cA^{\lambda}_{\GG}$ is normalized to be perverse.  
 
 Recall the ind-scheme $\Char(\Pi_a)$ from Section~\ref{Sect_Parabolic subgroups}. Define the closed subscheme $_{\lambda}\Pi_a\subset \Pi_a(N)$ as follows. A point $v\in \Pi_a(N)$ lies in $_{\lambda}\Pi_a$ if the following  holds:
 
\begin{itemize} 
\item[C1)] $v\in\Char(\Pi_a)$;
\item[C2)] for $i=1,\ldots,n$ the map $\wedge^i v: \wedge^i(U_a\otimes C_a^{-1})\to (\wedge^i M_a)(\<-w_0^{\GG}(\lambda), \check{\omega}_i\>)$
is regular.
\end{itemize} 

 The subscheme $_{\lambda}\Pi_a$ is stable under translations by $\Pi_a(-r)$, so there is a closed subscheme $_{\lambda, N}\Pi_a\subset {_{N,r}\Pi_a}$ such that $_{\lambda}\Pi_a$ is the preimage of $_{\lambda,N}\Pi_a$ under the projection $\Pi_a(N)\to {_{N,r}\Pi_a}$. Since all the weights of $\wedge^i M_a$ are $\le \check{\omega}_i$ and $\ge w_0^{\GG}(\check{\omega}_i)$, the map (\ref{map_pi_for_Hecke_GG}) factors through the closed subscheme $_{\lambda, N}\Pi_a\subset {_{N,r}\Pi_a}$. 
 
 For each $v\in \Char(\Pi_a)$ let us define a $\cO$-lattice $M_v\subset M_a(F)$ as follows. View $v$ as a map $U_a\otimes C_a^{-1}\to M_a(F)$. For a $\cO$-lattice $R\subset M_a(F)$ set 
$$
R^{\perp}=\{m\in M_a(F)\mid \<m, x\>\in A_a(-a)\;\;\mbox{for all}\; x\in R\}
$$ 
Consider two cases. 
 
\medskip 
\noindent 
{\scshape CASE:} $a$ is even. For $v\in \Char(\Pi_a)$ set $R_v=v(U_a\otimes C_a^{-1})+M_a(-\frac{a}{2})$ and $M_v=v(U_a\otimes C_a^{-1})+R_v^{\perp}$. 
Then $R_v^{\perp}\subset M_v\subset R_v$, and the induced form $\wedge^2 M_v\to A_a(-a)$ is regular and nondegenerate. So, $M_v\in \Gr_{\GG_a}^{-a}$. 

\medskip 
\noindent 
{\scshape CASE:} $a$ is odd. Let $b=(-a-1)/2$. Note that $(M_a(b))^{\perp}=M_a(b+1)$. Set $R_v=v(U_a\otimes C_a^{-1})+M_a(b+1)$ and $M_v=v(U_a\otimes C_a^{-1})+R_v^{\perp}$. Clearly, the induced form $\wedge^2 M_v\to A_a(-a)$ is regular, but still can be degenerate. We call $v$ \select{generic} if the form $\wedge^2 M_v\to A_a(-a)$ is nondegenerate. In this case $M_v\in\Gr_{\GG_a}^{-a}$. 

\smallskip

 For $a$ even we get a stratification of $\Char(\Pi_a)$ indexed by $\{\lambda\in\Lambda^+_{\GG}\mid \<\lambda, \check{\omega}_0\>=a\}$, the stratum $_{\lambda}\Char(\Pi_a)$ is given by the condition that $M_v\in\Gr_{\GG_a}^{\lambda}$. This condition is also equivalent to requiring that there is an isomorphism of $\cO$-modules 
$$
R_v/(M_a(-a/2))\,\iso\, \cO/t^{a_1-\frac{a}{2}}\oplus\ldots\oplus \cO/t^{a_n-\frac{a}{2}},
$$ 
where $t\in\cO$ is a uniformizer. 

 Clearly, $_{\lambda}\Char(\Pi_a)\subset {_{\lambda}\Pi_a}$. There is a unique open subscheme $_{\lambda,N}\Pi^0_a\subset {_{\lambda,N}\Pi_a}$ whose preimage under the projection $_{\lambda}\Pi_a\to {_{\lambda,N}\Pi_a}$ equals $_{\lambda}\Char(\Pi_a)$. 
 
 We say that a morphism of free $\cO$-modules $M_1\to M_2$ is \select{maximal} if it does not factor through $M_2(-1)\subset M_2$. 
 
 For $a$ odd define $_{\lambda}\Char(\Pi_a)\subset {_{\lambda}\Pi_a}$ as the open subscheme given by the condition that each map $\wedge^i v$ in C2) is maximal. Then there is an open subscheme $_{\lambda,N}\Pi^0_a\subset {_{\lambda,N}\Pi_a}$ whose preimage under the projection $_{\lambda}\Pi_a\to {_{\lambda,N}\Pi_a}$ equals $_{\lambda}\Char(\Pi_a)$. One checks that any $v\in {_{\lambda}\Char(\Pi_a)}$ is generic and the corresponding lattice $M_v$ satisfies $M_v\in \Gr_{\GG_a}^{\lambda}$. Note that for $v\in {_{\lambda}\Char(\Pi_a)}$ we have an isomorphism of $\cO$-modules
$$
R_v/(M_a(b+1))\,\iso\, \cO/t^{a_1-(a+1)/2}\oplus\ldots\oplus \cO/t^{a_n-(a+1)/2}
$$
for any uniformizer $t\in\cO$. 

 Write $\IC(_{\lambda,N}\Pi^0_a)$ for the intersection cohomology sheaf of $_{\lambda,N}\Pi^0_a$. 

\begin{Pp} 
\label{Pp_action_of_Sph_GG}
Let $\lambda\in \Lambda^+_{\GG}$ with $\<\lambda, \check{\omega}_0\>=a$. \\
1) The map
$$
\pi: {_{0,r}\Pi}\ttimes \Grb_{\GG_a}^{\lambda}\to 
{_{\lambda, N}\Pi_a}
$$
is an isomorphism over the open subscheme $_{\lambda,N}\Pi^0_a$. \\ % So, $\dim {_{\lambda,N}\Pi^0_a}=?$. 
2) Assume $m>n$ then one has a canonical isomorphism 
$
\H^{\lambda}_{\GG}(I_0)\,\iso\, \IC(_{\lambda,N}\Pi^0_a)
$.
\end{Pp}
\begin{Prf}
1) The fibre of $\pi$ over $v\in {_{\lambda,N}\Pi^0_a}$ is the scheme classifying lattices $M'\in\Grb_{\GG_a}^{\lambda}$ such that $v(U_a\otimes C_a^{-1})\subset M'$. Given such a lattice $M'$ let us show that $M_v=M'$. 

 Consider first the case of $a$ odd. The inclusion $R_v\subset M'+M_a(b+1)$ must be an equality, because for $M'\in \Gr_{\GG_a}^{\mu}$ with $\mu\le\lambda$ we have 
$$
\dim(M'+M_a(b+1))/(M_a(b+1))=\epsilon(\mu)\le \epsilon(\lambda)=\dim R_v/(M_a(b+1))
$$
We have denoted here $\epsilon(\mu)=\<\mu,\check{\omega}_n\>-\frac{n}{2}(a+1)$.  So, $M_v=v(U_a\otimes C_a^{-1})+(M'\cap M_a(b))\subset M'$ is also an equality, because both $M_v$ and $M'$ have symplectic forms with values in $A_a(-a)$. 

 The case of $a$ even  is quite similar to (\cite{L2}, Lemma~13). Namely, the inclusion $R_v\subset M'+M_a(-\frac{a}{2})$ must be an equality, because for $M'\in\Gr^{\mu}_{\GG_a}$ with $\mu\le\lambda$ we get
$$
\dim(M'+M_a(-a/2))/(M_a(-a/2))=\epsilon(\mu)\le \epsilon(\lambda)=\dim R_v/(M_a(-a/2))
$$
Here for $a$ even we have set $\epsilon(\mu)=\<\mu,\check{\omega}_n\>-\frac{n}{2}a$. So, 
$$
M_v=v(U_a\otimes C_a^{-1})+(M'\cap (M_a(-a/2)))\subset M'
$$
is also an equality. The first assertion follows.

\smallskip\noindent
2) For $m\ge n$ the scheme $_{\lambda,N}\Pi^0_a$ is nonempty, so $\IC(_{\lambda,N}\Pi^0_a)$ appears in $\H^{\lambda}_{\GG}(I_0)$ with multiplicity one. So, it suffices to show that 
$$
\Hom(\H^{\lambda}_{\GG}(I_0), \H^{\lambda}_{\GG}(I_0))=\Qlb,
$$
where $\Hom$ is taken in the derived category $\D_{\GG Q\HH_a(\cO)}(\Pi_a(F))$. By adjointness, 
$$
\Hom(\H^{\lambda}_{\GG}(I_0), \H^{\lambda}_{\GG}(I_0))\,\iso\,
\Hom(\H^{-w_0^{\GG}(\lambda)}_{\GG}\H^{\lambda}_{\GG}(I_0), I_0),
$$
where $\Hom$ in the RHS is taken in $\D_{\GG Q\HH_0(\cO)}(\Pi_0(F))$. We are reduced to show that for  any $0\ne \mu\in\Lambda^+_{\GG}$ with $\<\mu,\check{\omega}_0\>=0$ one has 
$$
\Hom(\H^{\mu}_{\GG}(I_0), I_0)=0
$$ 
in $\D_{\GG Q\HH_0(\cO)}(\Pi_0(F))$. The latter assertion is true for $m>n$, it is proved in (\cite{L2}, part 2) of Lemma~13). 
\end{Prf}

\medskip
\begin{Rem} 
\label{Rem_Weil_a_depends_on_parity}
For any $a,b\in\ZZ$ let us construct an equivalence $\Weil_a\,\iso\, \Weil_{a+2b}$. Pick isomorphisms of $\cO$-modules 
\begin{equation}
\label{iso_a_and_a+2b}
L_a(b)\,\iso\, L_{a+2b}, \;\; A_a(2b)\,\iso\, A_{a+2b}, \;\; U_a\,\iso\, U_{a+2b}, 
\end{equation}
They yield isomorphisms $C_a\,\iso\, C_{a+2b}$, 
$V_a\,\iso\, V_{a+2b}$, $M_a(b)\,\iso\, M_{a+2b}$. Hence, also isomorphisms $Q(\GG_a)\,\iso\, Q(\GG_{a+2b})$, $\GG_a\,\iso\, \GG_{a+2b}$ of group schemes over $\Spec\cO$ (and similarly for $\HH$). We also get isomorphisms of group schemes over $\Spec\cO$
$$
Q\GG\HH_a\,\iso\, Q\GG\HH_{a+2b}, \;\; \GG Q\HH_a\,\iso\, \GG Q\HH_{a+2b},\;\;  \HH Q\GG_a\,\iso\, \HH Q\GG_{a+2b}
$$
The isomorphisms (\ref{iso_a_and_a+2b}) also yield $\Pi_a(b)\,\iso\, \Pi_{a+2b}$ and $\Upsilon_a(b)\,\iso\, \Upsilon_{a+2b}$. In turn, we get equivalences 
$$
\P_{\HH Q\GG_a}(\Upsilon_a(F))\,\iso\, \P_{\HH Q\GG_{a+2b}}(\Upsilon_{a+2b}(F)),\;\;
\P_{\GG Q\HH_a}(\Pi_a(F))\,\iso\, \P_{\GG Q\HH_{a+2b}}(\Pi_{a+2b}(F))
$$
which yield the desired equivalence $\Weil_a\,\iso\, \Weil_{a+2b}$. The diagram commutes
$$
\begin{array}{ccc}
_{-a}\Sph_{\GG} & \to & \Weil_a\\
\downarrow\lefteqn{\epsilon} && \downarrow\lefteqn{\scriptstyle \!\!\wr}\\
_{-a-2b}\Sph_{\GG} & \to & \Weil_{a+2b},
\end{array}
$$
where the horizontal arrows are given by $\cS\mapsto \H^{\la}_{\GG}(\cS, I_0)$, and $\epsilon$, at the level of representations of $\check{G}$, is given by $V\mapsto V\otimes V^{b\omega}$. Here $V^{\omega}$ is the one-dimensional representation of $\check{\GG}$ with highest weight $\omega$ such that $\<\omega,\check{\omega}_0\>=2$. 
So, the case of $a$ even in Proposition~\ref{Pp_action_of_Sph_GG} also follows from (\cite{L2}, Lemma~13). 
\end{Rem}

\sssec{} 
\label{Sect_finite_field_k_0}
%old section 4.8.7.2 
Let $k_0\subset k$ be a finite subfield. In this subsection we assume that all the objects of Sections~4 are defined over $k_0$. In particular, $\cO_0\subset\cO$ is a complete discrete valuation $k_0$-algebra, and $F_0$ its fraction field.

 Write $\Weil_{a,k_0}$ for the category of triples $(\cF_1,\cF_2,\beta)$ as in Definition~\ref{Def_Weil_a} of
$\Weil_a$ but with a $k_0$-structure and, as such, pure of weight zero. It is understood that the Fourier transform functors are normalized to preserve purity. Note that for any $(\cF_1,\cF_2,\beta)\in \Weil_a$ the perverse sheaf $\cF_1$ is $\Gm$-equivariant with respect to the homotheties on $\Upsilon_a(F)$, this follows from the $\GL(L_a)(\cO)$-equivariance of $\cF_1$. 
 
 Denote by $\D\Weil_{a, k_0}$ the category of bounded complexes as in the definition of $\D\Weil_a$ but, in addition, with a  $k_0$-structure and, as such, pure of weight zero. So, for an object of $\D\Weil_{a, k_0}$ its semi-simplification is a bounded complex of the form $\oplus_{i\in\ZZ} F_i[i](\frac{i}{2})$ with $F_i\in\Weil_{a, k_0}$. 

 For a totally disconnected locally compact space $Y$ write $\cS(Y)$ for the Schwarz space of locally constant $\Qlb$-valued functions on $Y$ with compact support. 
Write $\Weil_a(k_0)$ for the $\Qlb$-vector space of pairs $(\cF_1, \cF_2)$, where $\cF_1\in \cS_{\HH Q\GG_a(\cO_0)}(\Upsilon_a(F_0))$, $\cF_2\in \cS_{\GG Q\HH_a(\cO_0)}(\Pi_a(F_0))$ with $\zeta_a(\cF_1)\,\iso\,\cF_2$. 

 Write $\cP$ for the composition of functors 
$$
\D\Weil_a\,\toup{f_{\HH}} \,\D\P_{\HH Q\GG_a(\cO)}(\Upsilon_a(F))^c\toup{J^*_{P_{\HH_a}}} 
\D_{Q\GG\HH_a(\cO)}(Q\Upsilon_a(F))^c,
$$
where $f_{\HH}$ sends $(\cF_1,\cF_2,\beta)$ to $\cF_1$. Write $\bar\cP: \D\Weil_{a, k_0}\to \D_{Q\GG\HH_a(\cO_0), mixed}(Q\Upsilon_a(F_0))^c$ for the similarly defined functor over $k_0$. Here we denoted by
$$
\D_{Q\GG\HH_a(\cO_0), mixed}(Q\Upsilon_a(F_0))\subset 
\D_{Q\GG\HH_a(\cO_0)}(Q\Upsilon_a(F_0)) 
$$
the full subcategory of mixed complexes (\cite{BBD}, 5.1.5). The following is similar to (\cite{L2}, Proposition~7).
\begin{Pp}
\label{Pp_about_functor_cP}
For $i=1,2$ let $K_i\in \D\Weil_{a, k_0}$. If $\bar\cP(K_1)\,\iso\, \bar\cP(K_2)$ 
then $K_1\,\iso\, K_2$ noncanonically in $\D\Weil_a$.  
\end{Pp} 
\begin{Prf}  
Write $K_{k_0}$ (resp., $DK_{k_0}$) for the Grothendieck group of the category $\Weil_{a, k_0}$ (resp., of $\D\Weil_{a, k_0}$). Note that $DK_{k_0}\,\iso\, K_{k_0}\otimes_{\ZZ} \ZZ[t, t^{-1}]$. Denote by $\Upsilon K_{k_0}$ the Gro\-then\-dieck group of the category $\D_{Q\GG\HH_a(\cO_0), mixed}(Q\Upsilon_a(F_0))^c$.
%of pure complexes of weight zero on $Q\Upsilon_a(F_0)$, whose all perverse cohomologies lie in $\P_{Q\GG\HH_a(\cO_0)}(Q\Upsilon_a(F_0))$. 

 Recall that the hyperbolic localization (of equivariant complexes) preserves purity (\cite{Brad}, Theorem~2). The functor $J^*_{P_{\HH_a}}$ yields a homomorphism $J^*_{P_{\HH_a}}: DK_{k_0}\to \Upsilon K_{k_0}$. Let us show that it is injective. 
 Let $\cF$ be an objects in its kernel. For a finite subfield $k_0\subset k_1\subset k$ let $\cO_1\subset F_1$ be obtained from $\cO_0\subset F_0$ by the extension of scalars $k_0\to k_1$. The trace of Frobenius map $\tr_{k_1}$ over $k_1$ fits into the diagram
\begin{equation}
\label{diag_for_functor_cP}
\begin{array}{ccc}
DK_{k_0} & \toup{J^*_{P_{\HH_a}}} & \Upsilon K_{k_0}\\
\downarrow\lefteqn{\scriptstyle \tr_{k_1}} && \downarrow\lefteqn{\scriptstyle \tr_{k_1}}\\
\Weil_a(k_1) & \toup{J_{k_1}} & \cS_{Q\GG\HH_a(\cO_1)}(Q\Upsilon_a(F_1))
\end{array}
\end{equation}
By (Lemma~\ref{Lm_injectivity_appendixA}, Appendix~A), $J_{k_1}$ is injective, so $\tr_{k_1}(\cF)=0$ for any finite extension $k_0\subset k_1$. By the result of Laumon (\cite{La}, Theorem~1.1.2) this implies $\cF=0$ in $DK_{k_0}$. Finally, if $K_1=K_2$ in $DK_{k_0}$ then $K_1\,\iso\, K_2$ in $\D\Weil_a$ (cf. \cite{L2}, Remark~10).
\end{Prf}

\medskip
 
The following result will not be used in this paper, its proof is found in Appendix~A. 
 
\begin{PpA} 
\label{PpA.1}
Assume $m>n$. The map 
$K_0({_{-a}\Sph_{\GG}})\otimes\Qlb\to \Weil_a(k_0)$ given by $\cS\mapsto \tr_{k_0} \H^{\la}_{\GG}(\cS, I_0)$ is an isomorphism of $\Qlb$-vector spaces. 
\end{PpA} 

 Write $\Weil_a^{ss}\subset\Weil_a$ for the full subcategory of semi-simple objects. 

\begin{Con}
\label{Con_Weil_a_m_greater_n}
Assume $m>n$. The functor $_{-a}\Sph_{\GG}\to \Weil_a^{ss}$ given by $\cS\to \H^{\la}_{\GG}(\cS, I_0)$ is an equivalence of categories.
\end{Con} 

\ssec{Action of $\Sph_{\HH}$}
\label{Sect_Action_Sph_HH}

\sssec{}
We write $V^{\check{\lambda}}$ for the irreducible $\HH$-module with highest weight $\check{\lambda}$. Let $V_0, C_0$,  $\check{\alpha}_0$ be as in Section~\ref{Sect_3.2_The_group_HH}. 
 For $0< i<m$ let $\check{\alpha}_i$ denote the highest weight of the irreducible $\HH$-module $\wedge^i V_0$. Remind that
$$
\wedge^m V_0\,\iso\, V^{\check{\alpha}_m}\oplus V^{\check{\alpha}'_m}
$$
is a direct sum of two irreducible representations, this is our definition of $\check{\alpha}_m, \check{\alpha}'_m$. Say that a maximal isotropic subspace $\cL\subset V_0$ is \select{$\check{\alpha}_m$-oriented} (resp., \select{$\check{\alpha}'_m$-oriented}) if $\wedge^m \cL\subset V^{\check{\alpha}_m}$ (resp., $\wedge^m \cL\subset V^{\check{\alpha}'_m}$). The group $\HH$ has two orbits on the variety of maximal isotropic subspaces in $V_0$ given by the orientation.  

 Remind that $\Gr_{\HH}^b$ classifies lattices $V'\subset V_0(F)$ such that the induced form $\Sym^2 V'\to C(b)$ is regular and nondegenerate, here $C=C_0(\cO)$. 
   
 Let $\lambda\in\Lambda^+_{\HH}$, set $a=\<\lambda,\check{\alpha}_0\>$. Remind that $\cA^{\lambda}_{\HH}\in\Sph_{\HH}$ denotes the IC-sheaf of $\Grb_{\HH}^{\lambda}$, so $\cA^{\lambda}_{\HH}\in{_{-a}\Sph_{\HH}}$. By definition, the complex
$$
\H^{\lambda}_{\HH}(I_0)=
\H^{\la}_{\HH}(\cA^{\lambda}_{\HH}, I_0)\in \D_{\HH Q\GG_a}(\Upsilon_a(F))
$$
is as follows. Set $r=\<\lambda, \check{\alpha}_1\>$ and $N=\<-w_0^{\HH}(\lambda), \check{\alpha}_1\>$. Let $_{0,r}\Upsilon\ttimes \Grb^{\lambda}_{\HH_a}$ be the scheme classifying $h\in \Grb^{\lambda}_{\HH_a}$, $x\in 
L_a^*\otimes A_a\otimes ((hV_a)/V_a(-r))$. Let 
\begin{equation}
\label{map_pi_for_Hecke_HH}
\pi: {_{0,r}\Upsilon\ttimes \Grb^{\lambda}_{\HH_a}}\to {_{N,r}\Upsilon_a}
\end{equation}
be the map sending $(x, h\HH_a(\cO))$ to $x$. Then $\H^{\lambda}_{\HH}(I_0)\,\iso\, \pi_!(\Qlb\tboxtimes \cA^{\lambda}_{\HH})$ canonically, where $\Qlb\tboxtimes \cA^{\lambda}_{\HH}$ is normalized to be perverse. 

 View a point of $\Upsilon_a(F)$ as a map $L_a\otimes A_a^*\to V_a(F)$. Define a closed subscheme $_{\lambda}\Upsilon_a\subset \Upsilon_a(N)$ as follows. A point $v\in\Upsilon_a(N)$ lies in $_{\lambda}\Upsilon_a$ if the following conditions hold:
 
\begin{itemize}
\item[C1)] $v\in\Char(\Upsilon_a)$;

\item[C2)] for $1\le i<m$ the map $\wedge^i v: \wedge^i(L_a\otimes A_a^*)\to (\wedge^i V_a)(\<-w_0^{\HH}(\lambda), \check{\alpha}_i\>)$ is regular;
\item[C3)] the map 
$$
(v_m,v'_m): \wedge^m( L_a\otimes A_a^*)\to V_a^{\check{\alpha}_m}(\<-w_0^{\HH}(\lambda), \check{\alpha}_m\>)\oplus
V_a^{\check{\alpha}'_m}(\<-w_0^{\HH}(\lambda), \check{\alpha}'_m\>)
$$ 
induced by $\wedge^m v$ is regular.
\end{itemize}
 
 The scheme $_{\lambda}\Upsilon_a$ is stable under translations by $\Upsilon_a(-r)$, so there is a closed subscheme $_{\lambda,N}\Upsilon_a\subset {_{N,r}\Upsilon_a}$ such that $_{\lambda}\Upsilon_a$ is the preimage of $_{\lambda,N}\Upsilon_a$ under the projection $\Upsilon_a(N)\to {_{N,r}\Upsilon_a}$.
Clearly, the map (\ref{map_pi_for_Hecke_HH}) factors through the closed subscheme $_{\lambda,N}\Upsilon_a\subset {_{N,r}\Upsilon_a}$. 

 For each $v\in\Char(\Upsilon_a)$ define a $\cO$-lattice $V_v\subset V_a(F)$ as follows. For a $\cO$-lattice $R\subset V_a(F)$ set
$$
R^{\perp}=\{x\in V_a(F)\mid \<x, y\>\in C_a(-a)\;\;\;\mbox{for all}\;\; y\in R\}
$$ 
Consider two cases. 

\smallskip
\noindent
{\scshape CASE:} $a$ is even. For $v\in\Char(\Upsilon_a)$ set $R_v=v(L_a\otimes A_a^*)+V_a(-\frac{a}{2})$ and 
$$
V_v=v(L_a\otimes A_a^*)+R_v^{\perp}
$$ 
Then $V_v\in \Gr_{\HH}^{-a}$. In this case we get a stratification of $\Char(\Upsilon_a)$ by locally closed subschemes $_{\lambda}\Char(\Upsilon_a)$ indexed by $\{\lambda\in\Lambda^+_{\HH}\mid \<\lambda,\check{\alpha}_0\>=a\}$. Namely, $v\in \Char(\Upsilon_a)$ lies in $_{\lambda}\Char(\Upsilon_a)$ iff $V_v\in\Gr_{\HH}^{\lambda}$. 

 Clearly, $_{\lambda}\Char(\Upsilon_a)\subset {_{\lambda}\Upsilon_a}$. There is a unique open subscheme $_{\lambda,N}\Upsilon_a^0\subset {_{\lambda,N}\Upsilon_a}$ whose preimage under the projection $_{\lambda}\Upsilon_a\to {_{\lambda,N}\Upsilon_a}$ equals $_{\lambda,N}\Upsilon_a^0$. 

\medskip\noindent
{\scshape CASE:} $a$ is odd. Let $b=(-a-1)/2$. We have $(V_a(b+1))^{\perp}=V_a(b)$. Set 
$$
R_v=v(L_a\otimes A_a^*)+V_a(b+1)
$$ 
and $V_v=v(L_a\otimes A_a^*)+R_v^{\perp}$. Then the induced form $\Sym^2 V_v\to C_a(-a)$ is regular, but still can be degenerate. We call $v$ \select{generic} if the form $\Sym^2 V_v\to C_a(-a)$ is nondegenerate. In this case $V_v\in\Gr_{\HH}^{-a}$. 

 For $a$ odd define an open subscheme $_{\lambda}\Char(\Upsilon_a)\subset {_{\lambda}\Upsilon_a}$ as follows. Note that $\<w_0^{\HH}(\lambda), \check{\alpha}_m-\check{\alpha}'_m\>\ne 0$. A point $v\in {_{\lambda}\Upsilon_a}$ lies in $_{\lambda}\Char(\Upsilon_a)$ if the following conditions hold: 
 
\begin{itemize}
\item the maps in C2) are maximal; 
\item if $\<w_0^{\HH}(\lambda), \check{\alpha}_m-\check{\alpha}'_m\> <0$ then $v_m$ in C3) is maximal, otherwise $v'_m$ in C3) is maximal.  
\end{itemize}
There is a unique open subscheme $_{\lambda,N}\Upsilon_a^0\subset {_{\lambda,N}\Upsilon_a}$ whose preimage under the projection $_{\lambda}\Upsilon_a\to {_{\lambda,N}\Upsilon_a}$ equals $_{\lambda}\Char(\Upsilon_a)$. 

 Write $\IC(_{\lambda,N}\Upsilon_a^0)$ for the intersection cohomology sheaf of $_{\lambda,N}\Upsilon_a^0$. 
 
\begin{Pp} Let $\lambda\in\Lambda_{\HH}^+$ with $\<\lambda, \check{\alpha}_0\>=a$. \\
1) The map
$$
\pi: {_{0,r}\Upsilon\ttimes \Grb^{\lambda}_{\HH_a}}\to {_{\lambda,N}\Upsilon_a}
$$
is an isomorphism over the open subscheme $_{\lambda,N}\Upsilon_a^0$. \\
2) Assume $m\le n$ then one has a canonical isomorphism $\H^{\lambda}_{\HH}(I_0)\,\iso\, \IC(_{\lambda,N}\Upsilon_a^0)$.
\end{Pp}
\begin{Prf}
1) Let $v\in {_{\lambda,N}\Upsilon_a^0}$. 
The fibre of $\pi$ over $v$ is the scheme classifying lattices $V'\in \Grb^{\lambda}_{\HH_a}$ such that $v(L_a\otimes A_a^*)\subset V'$. Given such a lattice $V'$ let us show that $V_v=V'$. 

 In view of Remark~\ref{Rem_Weil_a_depends_on_parity} the case of $a$ even is reduced to the case $a=0$, and the latter is done in (\cite{L2}, Lemma~9). 
 
 Consider the case of $a$ odd. The inclusion $R_v\subset V'+V_a(b+1)$ must be an equality, because for $V'\in \Gr^{\mu}_{\HH_a}$ with $\mu\le\lambda$ we have
$$
\dim(V'+V_a(b+1))/(V_a(b+1))=\epsilon(\mu)\le \epsilon(\lambda)=\dim R_v/(V_a(b+1)) 
$$
We have denoted here 
$$
\epsilon(\mu)=-m(b+1)+\max\{\<-w_0^{\HH}(\mu), \check{\alpha}_m\>, \<-w_0^{\HH}(\mu), \check{\alpha}'_m\>\}
$$ 

 It follows that $V_v=v(L_a\otimes A_a^{-1})+(V'\cap V_a(b))\subset V'$. To prove that $V'=V_v$, it suffices to show that $v$ is generic. This follows from the fact that 
$(v(L_a\otimes A_a^{-1})+R_v^{\perp})/R_v^{\perp}$ is a maximal isotropic subspace in $R_v/R_v^{\perp}$. 

\smallskip\noindent
2) For $m\le n$ the scheme $_{\lambda,N}\Upsilon_a^0$ is nonempty, so $\IC(_{\lambda,N}\Upsilon_a^0)$ appears in $\H^{\lambda}_{\HH}(I_0)$ with multiplicity one. Now it remains to show that
$$
\Hom(\H^{\lambda}_{\HH}(I_0), \H^{\lambda}_{\HH}(I_0))=\Qlb,
$$
where $\Hom$ is taken in the derived category $\D_{\HH Q\GG_a(\cO)}(\Upsilon_a(F))$. By adjointness, 
$$
\Hom(\H^{\lambda}_{\HH}(I_0), \H^{\lambda}_{\HH}(I_0))\,\iso\, \Hom(\H^{-w_0^{\HH}(\lambda)}_{\HH} \H^{\lambda}_{\HH}(I_0), I_0),
$$
where $\Hom$ in the RHS is taken in $\D_{\HH Q\GG_0(\cO)}(\Upsilon_0(F))$. We are reduced to show that for any $0\ne \mu\in\Lambda^+_{\HH}$ with $\<\mu, \check{\alpha}_0\>=0$ one has
$$
\Hom(\H^{\mu}_{\HH}(I_0), I_0)=0
$$
in $\D_{\HH Q\GG_0(\cO)}(\Upsilon_0(F))$. For $m\le n$ this is proved in (\cite{L2}, part 2) of Lemma~9). 
\end{Prf}
 
\sssec{} As in the case $m>n$, assume for a moment that $k_0\subset k$ is a finite subfield, and all the objects introduced in Section~4 have a $k_0$-structure. The following result can be proved exactly as Proposition~A.1 (it is not used in the paper, the proof is omitted). 
 
\begin{PpAA} 
\label{PpA.2}
Assume $m\le n$. Then the map 
$K_0({_{-a}\Sph_{\HH}})\otimes\Qlb \to \Weil_a(k_0)$ given by $\cS\mapsto \tr_{k_0}\H^{\la}_{\HH}(\cS, I_0)$ is an isomorphism of $\Qlb$-vector spaces. 
\end{PpAA} 

\begin{Con}
\label{Con_Weil_a_m_less_n}
Assume $m\le n$. The functor $_a\Sph_{\HH}\to\Weil_a^{ss}$ given by $\cS\mapsto \H^{\la}_{\HH}(\cS, I_0)$ is an equivalence of categories. 
\end{Con}
 
\ssec{Proof of Theorem~\ref{Th2_local}}
\label{Sect_Proof of Th2_local}

\sssec{} For a homomorphism of groups $h: H_1\times\Gm\to H_2$, recall that, by abuse of notations, we also write $h: H_1\times\Gm\to H_2\times\Gm$ for the map $(h, \pr)$, where $\pr: H_1\times\Gm\to\Gm$ is the projection.

\sssec{}
\label{Sect_4.12.2}
Use the notations of Sections~\ref{Sect_Action_of_Sph_GG} and \ref{Sect_Action_Sph_HH}. Assume that $U_0$ is $\check{\alpha}_m$-oriented, so $Q(\HH)$ acts on $\det U_0$ by the character
$\check{\alpha}_m=\check{\nu}_{\HH}$, the notation $\check{\nu}_{\HH}$ is that of Section~\ref{Sect_4.9.2}. As in Section~\ref{Sect_Root data}, we identify $\check{\omega}_0: \Gm\,\iso\,\check{\GL}(A_0)$ and $\check{\alpha}_0: \Gm\,\iso\, \check{\GL}(C_0)$. Recall the notations $\bar U_0, \bar L_0$ of Section~\ref{Sect_Root data}. We use the identifications $\check{\GL}(L_0)\,\iso\, \GL(\bar L_0)$, $\check{\GL}(U_0)\,\iso\, \GL(\bar U_0)$ of Section~\ref{Sect_Root data}. 
 
 For $m>n$ consider the decomposition $\bar U_0=\bar L_0
 \oplus {_1{\bar U}}$, where $\bar L_0$ is generated by $e_1,\ldots, e_n$, and $_1{\bar U}$ is generated by $e_{n+1},\ldots, e_m$. Let 
$$
\kappa_0: \check{\GL}(L_0)\times\Gm\to \check{\GL}(U_0)
$$ 
be the composition
$$
\GL(\bar L_0)\times\Gm
\toup{\tau\times\id}\GL(\bar L_0)\times\Gm
\toup{\id\times 2\check{\rho}_{\GL(_1{\bar U})}} \GL(\bar L_0)\times \GL(_1{\bar U})\toup{\Levi} \GL(\bar U_0)
$$
where $\tau$ is an automorphism $g\mapsto (^tg)^{-1}$
of $\check{\GL}(L_0)$. 
 
Let $\kappa_Q: \check{Q}(\GG)\times\Gm\to \check{Q}(\HH)$ be the map
$$
\check{\GL}(L_0)\times \check{\GL}(A_0)\times\Gm\to 
\check{\GL}(U_0)\times\check{\GL}(C_0)
$$
given by $(x,y,z)\mapsto (\kappa_0(x,z), y\omega_n(x))$,  here $\omega_n: \check{\GL}(L_0)\to\Gm$ is defined in Section~\ref{Sect_3.3.3_Levi_Q(GG)}. 

 The restriction of $\kappa_Q: \check{Q}(\GG)\times\Gm\to \check{Q}(\HH)$ to $\check{Q}(\GG)$ is the composition 
$$
\check{Q}(\GG)\toup{\tau_{\GG}}\check{Q}(\GG)\toup{i_Q}\check{Q}(\HH),
$$ 
where $\tau_{\GG}, i_Q$ are those of Section~\ref{Sect_3.5.1_case_m>n}. Indeed, in the notations of Section~\ref{Sect_Root data}, the map $\check{\Lambda}_{\GG}\to \check{\Lambda}_{\HH}$ given by $i_Q\tau_{\GG}$ sends $(a,b)$ to $(a, -b)$. So, $i_Q\tau_{\GG}$ sends $(a,b)\in \check{\Lambda}_1\subset \check{\Lambda}_{\GG}$ to the sum of $(-a, -b)\in\check{\Lambda}_1$ and of $(2a, 0)\in \check{\Lambda}_2$. This explains the appearance of $\omega_n$ in the above formula.

% Write $\kappa_{Q,ex}: \check{Q}(\GG)\times\Gm\to \check{Q}(\HH)\times\Gm$ for the map $(\kappa_Q,\pr)$, where $\pr: \check{Q}(\GG)\times\Gm\to\Gm$ is the projection.
   
 For $m\le n$ consider the decomposition $\bar L_0=\bar U_0\oplus {_1{\bar L}}$, where $\bar U_0$ is generated by $e_1,\ldots, e_m$, and $_1{\bar L}$ is generated by $e_{m+1},\ldots, e_n$. Let $\kappa_0: \check{\GL}(U_0)\times\Gm\to \check{\GL}(L_0)$ be the composition
$$
\GL(\bar U_0)\times\Gm\toup{\id\times 2\check{\rho}_{\GL({_1{\bar L}})}}\GL(\bar U_0)\times\GL({_1{\bar L}})\toup{\Levi} \GL(\bar L_0)\toup{\tau}\GL(\bar L_0),
$$ 
here $\tau(g)=(^tg)^{-1}$ for $g\in \GL(\bar L_0)$.

Denote by $\kappa_Q: \check{Q}(\HH)\times\Gm\to \check{Q}(\GG)$ the map 
$$
\check{\GL}(U_0)\times \check{\GL}(C_0)\times\Gm\to \check{\GL}(L_0)\times\check{\GL}(A_0)
$$ 
given by $(x,y,z)\mapsto (\kappa_0(x,z), y \alpha_m(x))$, here $\alpha_m: \check{\GL}(U_0)\to\Gm$ is defined in Section~\ref{Sect_3.2.3_Levi_Q(HH)}. The restriction of $\kappa_Q$ to $\check{Q}(\HH)$ equals the composition 
$$
\check{Q}(\HH)\toup{\tau_{\HH}}\check{Q}(\HH)\toup{i_Q}\check{Q}(\GG),
$$ 
here the notations $\tau_{\HH}, i_{\kappa}$ are those of Section~\ref{Sect_3.4_m_le_n}. 

 Recall the Hecke functors on $Q\Upsilon_a(F)$ for $Q(\GG)$ defined in Section~\ref{Sect_4.8.4}. 

\begin{Pp} 
\label{Pp_geom_theta_lifting_GL_n_GL_m}
1) For $m>n$ the two functors $_{-a}\Sph_{Q(\HH)}\to \D_{Q\GG\HH_a(\cO)}(Q\Upsilon_a(F))$ given by
$$
\cT\mapsto \H^{\la}_{Q(\HH)}(\cT, I_0)\;\;\;\;\mbox{and}\;\;\;\; \cT\mapsto \H^{\la}_{Q(\GG)}(\gRes^{\kappa_Q}(\cT), I_0)
$$
are canonically isomorphic.\\
2) For $m\le n$ the two functors $_{-a}\Sph_{Q(\GG)}\to \D_{Q\GG\HH_a(\cO)}(Q\Upsilon_a(F))$ given by
$$
\cT\mapsto \H^{\la}_{Q(\GG)}(\cT, I_0)\;\;\;\;\mbox{and}\;\;\;\; \cT\mapsto \H^{\la}_{Q(\HH)}(\gRes^{\kappa_Q}(\cT), I_0)
$$
are canonically isomorphic.  
\end{Pp}
\begin{proof}
2) One has $\kappa_Q\comp \check{\alpha}_0=\check{\omega}_0$. So, if $\cT\in {_{-a}\Sph_{Q(\GG)}}$ then $\gRes^{\kappa_Q}(\cT)$ vanishes off $\Gr_{Q(\HH)}^{-a}$. For a dominant weight $\lambda$ of $\check{Q}(\GG)$ write temporary $V^{\lambda}$ for the irreducible representation of $\check{Q}(\GG)$ with highest weight $\lambda$. Recall the notations $\bar\omega_{\HH}, \alpha_m$ from Section~\ref{Sect_3.2.3_Levi_Q(HH)}. One gets $\Loc(V^{\bar\omega_{\HH}})\in {_{-1}\Sph_{Q(\HH)}}$ and $\H^{\la}_{Q(\HH)}(V^{\bar\omega_{\HH}}, I_0)\,\iso\, I_1$, where $I_1$ is the constant perverse sheaf on $Q\Upsilon_1$. Besides,  
$\H^{\la}_{Q(\HH)}(V^{\alpha_m}, I_1)$ is the constant perverse sheaf on $Q\Upsilon_1(-1)$. 

 For a dominant weight $\lambda$ of $\check{Q}(\GG)$
write temporary $W^{\lambda}$ for the irreducible representation of $\check{Q}(\GG)$ with highest weight $\lambda$. Recall the notation $\bar\omega_{\GG}$ from Section~\ref{Sect_3.3.3_Levi_Q(GG)}. One gets $\Loc(W^{\bar\omega_{\GG}})\in {_{-1}\Sph_{Q(\GG)}}$ and 
$\H^{\la}_{Q(\GG)}(W^{\bar\omega_{\GG}}, I_0)$ is the constant perverse sheaf on $Q\Upsilon_1(-1)$. Thus,
$$
\H^{\la}_{Q(\GG)}(W^{\bar\omega_{\GG}}, I_0)\,\iso\, \H^{\la}_{Q(\HH)}(V^{\alpha_m+\bar\omega_{\HH}}, I_0)
$$
% This explains the map $\kappa_Q$ composed with the projection to $\check{\GL}(A_0)$. 

 We are reduced to the following claim. For $\cT\in {_0\Sph_{Q(\GG)}}$ the two functors ${_0\Sph_{Q(\GG)}}\to \D_{Q\GG\HH_0(\cO)}(Q\Upsilon_0(F))$ given by
$$
\cT\mapsto \H^{\la}_{Q(\GG)}(\cT, I_0)\;\;\;\;\mbox{and}\;\;\;\; \cT\mapsto \H^{\la}_{Q(\HH)}(\gRes^{\kappa_0}(\cT), I_0)
$$
are isomorphic. This follows from (\cite{L2}, Proposition~4 and Corollary~5).

\medskip\noindent
1) similar to 2). 
\end{proof}

\sssec{} As in (\cite{L2}, Theorem~7), for each $a\in\ZZ$ the diagram of functors is canonically commutative
$$
\begin{array}{ccccc}
&& \D\!\Weil_a\\
& \swarrow\lefteqn{\scriptstyle f_{\HH}} && \searrow\lefteqn{\scriptstyle f_{\GG}}\\
\DP_{\HH Q\GG_a(\cO)}(\Upsilon_a(F)) &&&&
\DP_{\GG Q\HH_a(\cO)}(\Pi_a(F))\\
\downarrow\lefteqn{\scriptstyle J^*_{P_{\HH_a}}} &&&&
\downarrow\lefteqn{\scriptstyle J^*_{P_{\GG_a}}}\\
\D_{Q\GG\HH_a(\cO)}(Q\Upsilon_a(F)) && \toup{\Four_{\psi}} && \D_{Q\GG\HH_a(\cO)}(Q\Pi_a(F)),
\end{array}
$$
where $f_{\HH}$ (resp., $f_{\GG}$) sends $(\cF_1,\cF_2,\beta)$ to $\cF_1$ (resp., to $\cF_2$). 
 
 Recall the maps $\kappa_{\HH}: \check{Q}(\HH)\times\Gm\to\check{H}$ and $\kappa_{\GG}: \check{Q}(\GG)\times\Gm\to\check{\GG}$ from Sections~\ref{Sect_4.9.2}, 
\ref{Sect_4.9.4}.  The restriction of $\kappa_{\HH}$ and of $\kappa_{\GG}$ to $\Gm$ equals 
$$
2(\check{\rho}_{\HH}-\check{\rho}_{Q(\HH)})+nm\check{\alpha}_0-n\check{\alpha}_m
$$ 
and $2(\check{\rho}_{\GG}-\check{\rho}_{Q(\GG)})+mn\check{\omega}_0-m\check{\omega}_n$ respectively. From definitions one gets
$$
\left\{
\begin{array}{l}
2(\check{\rho}_{\HH}-\check{\rho}_{Q(\HH)})=(m-1)\check{\alpha}_m-\frac{m(m-1)}{2}\check{\alpha}_0\\ \\
2(\check{\rho}_{\GG}-\check{\rho}_{Q(\GG)})=(n+1)\check{\omega}_n-\frac{n(n+1)}{2}\check{\omega}_0
\end{array}
\right.
$$ 
% Write $\kappa_{\HH,ex}:  \check{Q}(\HH)\times\Gm\to\check{H}\times\Gm$ for the map, whose first component is $\kappa_{\HH}$ and the second $\check{Q}(\HH)\times\Gm\to\Gm$ is the projection, and similarly for $\kappa_{\GG,ex}$.  

  By Corollary~\ref{Cor_Braden_rules}, for $\cT\in {_{-a}\Sph_{\HH}}$ and $\cS\in {_{-a}\Sph_{\GG}}$ we get canonical isomorphisms
\begin{equation}
\label{iso_for_Step_1_th_2_first}
\bar\cP (\H^{\la}_{\HH}(\cT, I_0))\,\iso\, \H^{\la}_{Q(\HH)}(\gRes^{\kappa_{\HH}}(\cT), I_0)
\end{equation}
and
\begin{equation}
\label{iso_for_Step_1_th_2_second}
\bar\cP (\H^{\la}_{\GG}(\cS, I_0))\,\iso\, \H^{\la}_{Q(\GG)}(\gRes^{\kappa_{\GG}}(\cS), I_0)
\end{equation}
in $\D_{Q\GG\HH_a(\cO_0), mixed}(Q\Upsilon_a(F_0))$. % The Hecke functors in the RHS of (\ref{iso_for_Step_1_th_2_first}) and (\ref{iso_for_Step_1_th_2_second}) are from $\D_{Q\GG\HH_0(\cO_0)}(Q\Upsilon_0(F_0))$ to the category $\D_{Q\GG\HH_a(\cO_0)}(Q\Upsilon_a(F_0))$.

\sssec{CASE  $m>n$} 
\label{Sect_Proof_of_Th_local_case_m>n}
Proposition~\ref{Pp_geom_theta_lifting_GL_n_GL_m} together with (\ref{iso_for_Step_1_th_2_first}) yield a canonical isomorphism
$$
\bar \cP (\H^{\la}_{\HH}(\cT, I_0))\,\iso\, \H^{\la}_{Q(\GG)}(\gRes^{\kappa_{Q,ex}}\gRes^{\kappa_{\HH}}(\cT), I_0)
$$
Recall the automorphism $\tau_{\GG}$ of $\check{\GG}$ from Section~\ref{Sect_Case m>n_3.5}. Recall that $i_{\kappa}: \check{\GG}\to\check{\HH}$ denotes the restriction of $\kappa$ to $\check{\GG}$. We will check below that the diagram commutes
\begin{equation}
\label{diag_groups_Rallis_m_greater_n}
\begin{array}{ccc}
\check{\GG}\times\Gm & \toup{\tau_{\GG}\times\id} \check{\GG}\times\Gm\toup{\kappa} & \check{\HH}\\
\uparrow\lefteqn{\scriptstyle \kappa_{\GG}} && \uparrow\lefteqn{\scriptstyle 
\kappa_{\HH}}\\
\check{Q}(\GG)\times\Gm & \toup{\kappa_{Q}} & \check{Q}(\HH)\times\Gm 
\end{array}
\end{equation}
The above diagram together with (\ref{iso_for_Step_1_th_2_first}),  (\ref{iso_for_Step_1_th_2_second}) yield canonical isomorphisms 
\begin{multline*}
\bar\cP(\H^{\la}_{\GG}(\ast\gRes^{\kappa}(\cT), I_0)\,\iso\, \H^{\la}_{Q(\GG)}(\gRes^{\kappa_{\GG, ex}}(\ast\gRes^{\kappa}(\cT)), I_0)\,\iso
\\  \H^{\la}_{Q(\GG)}(\gRes^{\kappa_{Q,ex}}\gRes^{\kappa_{\HH}}(\cT), I_0)\,\iso\, \H^{\la}_{Q(\HH)}(\gRes^{\kappa_{\HH}}(\cT), I_0)\,\iso\, \bar\cP(\H^{\la}_{\HH}(\cT, I_0))
\end{multline*}
Thus, we get a canonical isomorphism
$$
\bar\cP(\H^{\la}_{\GG}(\ast\gRes^{\kappa}(\cT), I_0))\,\iso\,
\bar\cP (\H^{\la}_{\HH}(\cT, I_0))
$$
Proposition~\ref{Pp_about_functor_cP} allows to lifts it to the desired isomorphism in $\D\Weil^{ss}_a$
$$
\H^{\la}_{\GG}(\ast\gRes^{\kappa}(\cT), I_0)\,\iso\, \H^{\la}_{\HH}(\cT, I_0)
$$

 First, removing the $\Gm$-factors the diagram (\ref{diag_groups_Rallis_m_greater_n}) becomes 
$$
\begin{array}{ccccc}
\check{\GG} & \toup{\tau_{\GG}} & \check{\GG} & \toup{i_{\kappa}} & \check{\HH}\\
\uparrow &&&& \uparrow\\
\check{Q}(\GG) & \toup{\tau_{\GG}} & \check{Q}(\GG) &\toup{i_Q} & \check{Q}(\HH)
\end{array}
$$
It commutes according to Section~\ref{Sect_Case m>n_3.5}. So, there is a unique $\delta_{\kappa}:\Gm\to \check{\TT}_{\HH}$ such that for $\kappa=(i_{\kappa},\delta_{\kappa})$ the diagram (\ref{diag_groups_Rallis_m_greater_n}) commutes. One gets
$$
\delta_{\kappa}+i_{\kappa}(m\check\omega_n-mn\check{\omega}_0
-2(\check{\rho}_{\GG}-\check{\rho}_{Q(\GG)}))=
2\check{\rho}_{\GL(_1{\bar U})}+ 2(\check{\rho}_{\HH}-\check{\rho}_{Q(\HH)})+nm\check{\alpha}_0-n\check{\alpha}_m
$$ 
One checks that $\kappa$ coincides with the map defined in Section~\ref{Sect_2.4.3_m_bigger_n}. If $m=n+1$ then $2\check{\rho}_{\GL(_1{\bar U})}=0$
and $\delta_{\kappa}$ is trivial. 

 The compatibility with the tensor structures on $\Rep(\check{\GG}), \Rep(\check{\HH})$ is proved exactly as in (\cite{L2}, Theorem~7).

\sssec{CASE  $m\le n$} 
\label{Sect_5.12.6_m_le_n}
Proposition~\ref{Pp_geom_theta_lifting_GL_n_GL_m} together with (\ref{iso_for_Step_1_th_2_second}) yield a canonical isomorphism
$$
\bar\cP (\H^{\la}_{\GG}(\cS, I_0))\,\iso\, \H^{\la}_{Q(\HH)}(\gRes^{\kappa_{Q,ex}}\gRes^{\kappa_{\GG}}(\cS), I_0)
$$
for $\cS\in {_{-a}\Sph_{\GG}}$. Recall the automorphism $\tau_{\HH}$ of $\check{\HH}$ from Section~\ref{Sect_3.4_m_le_n}. Recall that in Section~\ref{Sect_2.4.2_case_m_le_n} we defined the map $\kappa$ and its restriction $i_{\kappa}: \check{\HH}\to \check{\GG}$ to $\check{\HH}$. We check below that the following diagram commutes
\begin{equation}
\label{diag_groups_Rallis_m_less_n}
\begin{array}{ccc}
\check{\HH}\times\Gm & \!\!\toup{\tau_{\HH}\times\id} \; \check{\HH}\times\Gm \;\, \toup{\kappa}  \!\! & \check{\GG}\\
\uparrow\lefteqn{\scriptstyle \kappa_{\HH}} && \uparrow\lefteqn{\scriptstyle \kappa_{\GG}}\\
\check{Q}(\HH)\times\Gm & \toup{\kappa_{Q}} & \check{Q}(\GG)\times\Gm
\end{array}
\end{equation}
The above diagram together with (\ref{iso_for_Step_1_th_2_first}),  (\ref{iso_for_Step_1_th_2_second}) yield canonical isomorphisms
\begin{multline*}
\bar\cP (\H^{\la}_{\HH}(\ast\gRes^{\kappa}(\cS), I_0))\,\iso\, \H^{\la}_{Q(\HH)}(\gRes^{\kappa_{\HH}}(\ast\gRes^{\kappa}(\cS)), I_0)\,\iso\\  \H^{\la}_{Q(\HH)}(\gRes^{\kappa_{Q,ex}}\gRes^{\kappa_{\GG}}(\cS), I_0)\,\iso\,\H^{\la}_{Q(\GG)}(\gRes^{\kappa_{\GG}}(\cS), I_0)\,\iso\, \bar\cP (\H^{\la}_{\GG}(\cS, I_0))
\end{multline*}
Thus, we get a canonical isomorphism
$$
\bar\cP (\H^{\la}_{\HH}(\ast\gRes^{\kappa}(\cS), I_0))\,\iso\, \bar\cP (\H^{\la}_{\GG}(\cS, I_0))
$$
Proposition~\ref{Pp_about_functor_cP} allows to lift it to the desired isomorphism in $\D\Weil^{ss}_a$
$$
\H^{\la}_{\HH}(\ast\gRes^{\kappa}(\cS), I_0)\,\iso\, \H^{\la}_{\GG}(\cS, I_0)
$$
 
Removing the $\Gm$-factor, the diagram (\ref{diag_groups_Rallis_m_less_n}) becomes
$$
\begin{array}{ccccc}
\check{\HH} & \toup{\tau_{\HH}} & \check{\HH} & \toup{i_{\kappa}} & \check{\GG}\\
\uparrow &&&& \uparrow\\
\check{Q}(\HH) & \toup{\tau_{\HH}} & \check{Q}(\HH) & \toup{i_Q} & \check{Q}(\GG)
\end{array}
$$
It commutes according to Section~\ref{Sect_3.4_m_le_n}. So, there is a unique 
$\delta_{\kappa}: \Gm\to \check{\TT}_{\GG}$ such that for $\kappa=(i_{\kappa},\delta_{\kappa})$ the diagram (\ref{diag_groups_Rallis_m_less_n}) commutes. Our $\delta_{\kappa}$ is determined by
$$
\delta_{\kappa}-i_{\kappa}(2(\check{\rho}_{\HH}-\check{\rho}_{Q(\HH)})+nm\check{\alpha}_0-n\check{\alpha}_m)=
2(\check{\rho}_{\GG}-\check{\rho}_{Q(\GG)})+mn\check{\omega}_0-m\check{\omega}_n-2\check{\rho}_{\GL(_1{\bar L})}
$$
One checks that $(i_{\kappa},\delta_{\kappa})$ coincides with the map $\kappa$ defined in Section~\ref{Sect_2.4.2_case_m_le_n}. 
If $m=n$ then $2\check{\rho}_{\GL(_1{\bar L})}=0$ and $\delta_{\kappa}$ is trivial. 

The compatibility with the tensor structures on $\Rep(\check{\GG}), \Rep(\check{\HH})$ is proved exactly as in (\cite{L2}, Theorem~7).
Theorem~\ref{Th2_local} is proved.

\section{Global theory}
\label{Sect_Global_theory}

%\noindent
In Section~\ref{Sect_6.1_proof_Th} we derive Theorem~\ref{Th_Hecke_property_Aut_GH} from Theorem~\ref{Th_local_main}. In Section~\ref{Sect_6.2_proof_Th_main} we derive Theorem~\ref{Th_theta-lifting_functors-general} from Theorem~\ref{Th_Hecke_property_Aut_GH}. 

\ssec{Proof of Theorem~\ref{Th_Hecke_property_Aut_GH}}
\label{Sect_6.1_proof_Th}

\sssec{} To simplify notations, fix a closed point $\tilde x\in\tilde X$. Let $^{a\tilde x}\Bun_{G, H}$ be obtained from $^a\Bun_{G, H}$ by the base change $\Spec k\toup{\tilde x}\tilde X$.
We will establish isomorphisms (\ref{iso_Th_Hecke_Aut_GH_m_less_n}) and (\ref{iso_Th_Hecke_Aut_GH_m_greater_n}) over $^{a\tilde x}\Bun_{G,\tilde H}$. The fact that these isomorphisms depend on $\tilde x$ as expected is left to the reader. Set $x=\pi(\tilde x)$. 

  Recall the line bundle $\cE$ from Section~\ref{Sect_2.3.2_def_tilde_pi}, we have $\pi^*\cE\,\iso\, \cO_{\tilde X}$ canonically. So, the above choice of $\tilde x$ yields a trivialization $\cE\,\iso\,\cO\mid_{D_x}$ over $D_x=\Spec\cO_x$. The corresponding trivialization for $\sigma\tilde x$ is the previous one multiplied by $-1$. In Section~\ref{Sect_Dual_pair_GSp_GO} we worked with a complete discrete valuation algebra $\cO$. We will apply Theorem~\ref{Th_local_main} for $\cO$ replaced by $\cO_x$.

\sssec{} Recall the stack $^a\cXL$ and a line bundle $^a\cA_{\cXL}$ on it introduced in Section~\ref{Sect_5.2_stack_a_cXL}. A point of $^{a\tilde x}\Bun_{G, H}$ is given by a collection: $(M,\cA)\in\Bun_G$, $(V,\cC)\in\Bun_H$, and an isomorphism $\cA\otimes\cC\,\iso\,\Omega(ax)$. Let $^a\xi: {^{a\tilde x}\Bun_{G,H}}\to {^a\cXL}$ be the map sending $(M,\cA, V,\cC)$ to $(M,\cA, V,\cC)\mid_{D_x}$ together with the discrete lagrangian subspace $L=\H^0(X-x, M\otimes V)\subset M\otimes V(F_x)$. 
 
\begin{Lm} 
\label{Lm_line_bundle_on_rx_Bun_GH}
For a point $(\cM,\cA,\cV,\cC)$ of $^{a\tilde x}\Bun_{G,H}$ there is a canonical $\ZZ/2\ZZ$-graded isomorphism
$$
\det\RG(X, \cM\otimes\cV)\otimes\cC_x^{-anm}\,\iso\, 
\frac{\det\RG(X,\cM)^{2m}\otimes\det\RG(X,\cV)^{2n}}{\det\RG(X,\cC)^{2nm}\otimes\det\RG(X,\cO)^{2nm}}
$$
Here $\cC_x$ is of parity zero as $\ZZ/2\ZZ$-graded. 
\end{Lm}
\begin{Prf}
By (\cite{L}, Lemma~1), we get a canonical $\ZZ/2\ZZ$-graded isomorphism
$$
\det\RG(X, \cM\otimes\cV)\,\iso\, \frac{\det\RG(X,\cM)^{2m}\otimes\det\RG(X,\cV)^{2n}}{\det\RG(X, \cA^n)\otimes\det\RG(X,\det \cV)}\otimes
\frac{\det\RG(X, \cA^n\otimes\det\cV)}{\det\RG(X, \cO)^{4nm-1}}
$$
Applying this to $\cM=\cO^n\oplus\cA^n$ with natural symplectic form $\wedge^2\cM\to\cA$, we get
$$
\frac{\det\RG(X, \cV\otimes\cA)^n}{\det\RG(X, \cV)^n}\,\iso\, \frac{\det\RG(X, \cA)^{2nm}\otimes\det\RG(X, \cA^n\otimes\det \cV)}
{\det\RG(X, \cA^n)\otimes \det\RG(X, \det \cV)\otimes\det\RG(X, \cO)^{2nm-1}}
$$
Since $\cA\otimes\cC\,\iso\, \Omega(ax)$ and $\cV\,\iso\, \cV^*\otimes\cC$, the LHS of the above formula idetifies with 
$$
\det\RG(X, V/V(-ax))^{-n}\,\iso\, (\det V_x)^{-an}\otimes\det(\cO/\cO(-ax))^{-2mn}
$$ 
We have used a canonical $\ZZ/2\ZZ$-graded isomorphism 
$$
\det(V/V(-ax))\,\iso\, (\det V_x)^a\otimes (\det(\cO/\cO(-ax))^{2m}
$$ 
Since $\det V_x\,\iso\, \cC_x^m$, we get
$$
\det\RG(X,\cM\otimes\cV)\,\iso\, \frac{\det\RG(X,\cM)^{2m}\otimes\det\RG(X,\cV)^{2n}\otimes\cC_x^{-anm}}
{\det\RG(X, \cA)^{2nm}\otimes\det\RG(X, \cO)^{2nm}\otimes\det\RG(X, \cO/\cO(-ax))^{2nm}}
$$
To simplify the above expression, note that $\det\RG(X,\cA)\,\iso\, \det\RG(X, \cC(-ax))$ and 
$$
\det\RG(X,\cC)\,\iso\, \det\RG(X, \cC(-ax))\otimes \cC_x^a\otimes\det\RG(X, \cO/\cO(-ax))
$$ 
Our assertion follows.
\end{Prf}  
 
\sssec{} Let $^a\cA$ be the line bundle on $^{a\tilde x}\Bun_{G, H}$ with fibre $\det\RG(X, M\otimes V)\otimes C_x^{-anm}$ at $(M,\cA,V,\cC)$. We have canonically $(^a\xi)^*(^a\cA_{\cXL})\,\iso\, {^a\cA}$. Extend $^a\xi$ to a morphism $^a\tilde\xi: {^{a\tilde x}\Bun_{G, H}}\to {^a\wt\cXL}$ sending $(M,\cA,V,\cC)$ to its image under $^a\xi$ together with the one-dimensional space
$$
\cB=\frac{\det\RG(X,\cM)^m\otimes\det\RG(X,\cV)^n}{\det\RG(X,\cC)^{nm}\otimes\det\RG(X, \cO)^{nm}}
$$
equipped with the isomorphism $\cB^2\,\iso\, \det\RG(X, \cM\otimes\cV)\otimes\cC_x^{-anm}$ of Lemma~\ref{Lm_line_bundle_on_rx_Bun_GH}.  
 
\sssec{}  Let $^{a\tilde x}\cH_{G, H}$ be the stack classifying collections: a point $(M,\cA, M',\cA',\beta)\in {_x\cH_G}$ of the Hecke stack such that the isomorphism $\beta$ of the $G$-torsors $(M,\cA)$ and $(M',\cA')$ over $X-x$ induces an isomorphism $\cA(-ax)\,\iso\, \cA'$; a $H$-torsor $(V,\cC)\in\Bun_H$, and an isomorphism $\cA'\otimes\cC\,\iso\,\Omega$. We have the diagram
$$
^{a\tilde x}\Bun_{G,H} \getsup{h^{\la}} {^{a\tilde x}\cH_{G, H}} \toup{h^{\ra}} \Bun_{G, H},
$$   
where $h^{\ra}$ (resp., $h^{\la}$) sends the above point of $^{a\tilde x}\cH_{G, H}$ to $(M',\cA', V,\cC)\in\Bun_{G,H}$ (resp., to $(M,\cA, V,\cC)\in {^{a\tilde x}\Bun_{G, H}}$).  

 Restriction to $D_x$ gives rise to the diagram
\begin{equation}
\label{diag_for_Hecke_global_local_GG_HH}
\begin{array}{ccccc} 
^{a\tilde x}\Bun_{G,H} & \getsup{h^{\la}} &{^{a\tilde x}\cH_{G,H}} &\toup{h^{\ra}} &\Bun_{G,H} \\
\downarrow\lefteqn{\scriptstyle {^a\xi}} && \downarrow\lefteqn{\scriptstyle {^a\xi_G}}&& \downarrow\lefteqn{\scriptstyle ^0\xi}\\
^a\cXL & \getsup{h^{\la}} & {^{a,0}\cH_{\GG,\cXL}} & \toup{h^{\ra}}  & {^0\cXL},
\end{array}
\end{equation}
where the low row is the diagram (\ref{diag_Hecke_cXL_for_G}) for $a'=0$. Now Lemma~\ref{Lm_line_bundle_on_rx_Bun_GH} allows to extend (\ref{diag_for_Hecke_global_local_GG_HH}) to the following diagram, where both squares are cartesian
$$
\begin{array}{ccccc} 
^{a\tilde x}\Bun_{G, H} & \getsup{h^{\la}} &{^{a\tilde x}\cH_{G,H}} &\toup{h^{\ra}} &\Bun_{G, H} \\
\downarrow\lefteqn{\scriptstyle {^a\tilde\xi}} && \downarrow\lefteqn{\scriptstyle {^a\tilde\xi_G}}&& \downarrow\lefteqn{\scriptstyle ^0\tilde\xi}\\
^a\wt\cXL & \getsup{\tilde h^{\la}} & {^{a,0}\wt\cH_{\GG,\cXL}} & \toup{\tilde h^{\ra}}  & {^0\wt\cXL},
\end{array}
$$
and the low row is the diagram (\ref{diag_tilde_Hecke_cXL_for_G}) for $a'=0$ from Section~4.3.1. This provides an isomorphism 
$$
\H^{\la}_G(\cS, (^0\tilde\xi)^*K)\,\iso\, (^a\tilde\xi)^*\H^{\la}_{\GG}(\cS, K)
$$
on $^{a\tilde x}\Bun_{G, H}$ functorial in $\cS\in{_{-a}\Sph_{\GG}}$ and $K\in\D_{\cT_0}(\wt\cL_d(W_0(F_x)))$. Here the functors 
$$
(^a\tilde\xi)^*: \D_{\cT_a}(\wt\cL_d(W_a(F_x)))\to \D(^{a\tilde x}\Bun_{G, H})
$$ 
are defined as in (\cite{LL}, Section~7.2). 

\sssec{} Let $^{a\tilde x}\cH_{H, G}$ be the stack classifying collections: a point of the Hecke stack $(V,\cC, V',\cC', \beta)\in {_x\cH_H}$ such that the isomorphism $\beta$ of $H$-torsors $(V,\cC)$ and $(V',\cC')$ over $X-x$ induces an isomorphism $\cC(-ax)\,\iso\, cC'$; a $G$-torsor $(M,\cA)$ on $X$ and an ismorphism $\cA\otimes\cC'\,\iso\,\Omega$. As above, we get a diagram
$$
^{a\tilde x}\Bun_{G, H} \getsup{h^{\la}} {^{a\tilde x}\cH_{H, G}} \toup{h^{\ra}} \Bun_{G, H},
$$   
where $h^{\ra}$ (resp., $h^{\la}$) sends the above point of $^{a\tilde x}\cH_{H, G}$ to $(M,\cA, V',\cC')$ (resp., to $(M,\cA, V,\cC)$). 

 As in the case of the Hecke functor for $G$, we get the diagram, where both squares are cartesian
$$
\begin{array}{ccccc}
^{a\tilde x}\Bun_{G, H} &\getsup{h^{\la}} &{^{a\tilde x}\cH_{H, G}} &\toup{h^{\ra}} &\Bun_{G, H}\\
\downarrow\lefteqn{\scriptstyle {^a\tilde\xi}} && \downarrow\lefteqn{\scriptstyle {^a\tilde\xi_H}}&& \downarrow\lefteqn{\scriptstyle ^0\tilde\xi}\\
^a\wt\cXL & \getsup{\tilde h^{\la}} & {^{a,0}\wt\cH_{\HH,\cXL}} & \toup{\tilde h^{\ra}}  & {^0\wt\cXL},
\end{array}
$$
and the low row is the diagram (\ref{diag_Hecke_wt_HH_cXL}) for $a'=0$ from Section~4.3.2. This provides an isomorphism
$$
\H^{\la}_H(\cS, (^0\tilde\xi)^*K)\,\iso\, (^a\tilde\xi)^*\H^{\la}_{\HH}(\cS, K)
$$
on $^{a\tilde x}\Bun_{G, H}$ functorial in $\cS\in {_{-a}\Sph_{\HH}}$ and $K\in \D_{\cT_0}(\wt\cL_d(W_0(F_x)))$. By (\cite{LL}, Theorem~3), we have $(^0\tilde\xi)^* S_{W_0(F)}\,\iso\, \Aut_{G, H}$ naturally.
Now Theorem~\ref{Th_Hecke_property_Aut_GH} follows from Theorem~\ref{Th_local_main} by applying the functor $(^a\tilde\xi)^*$. Theorem~\ref{Th_Hecke_property_Aut_GH} is proved.
 
\ssec{Proof of Theorem~\ref{Th_theta-lifting_functors-general}}
\label{Sect_6.2_proof_Th_main}

\sssec{} We derive Theorem~\ref{Th_theta-lifting_functors-general} from Theorem~\ref{Th_Hecke_property_Aut_GH}. The proof is similar to (\cite{L2}, Theorem~3). We give the argument only for $m\le n$ (the case $m>n$ is similar). 
   
 Let $a\in\ZZ$. Let us establish the isomorphism (\ref{iso_theta_lifting_Hecke_Th1}) for any $\cS\in{_{-a}\Sph_{\GG}}$. By base change theorem, for $K\in \D(\Bun_H)$ we get
$$
\H^{\la}_G(\cS, F_G(K))\,\iso\, (^a\gp)_!(^a\gq^*K\otimes \H^{\la}_G(\cS, \Aut_{G, H}))[-\dim\Bun_H],
$$ 
where $^a\gq: {^a\Bun_{G,H}}\to \Bun_H$ and $^a\gp: {^a\Bun_{G, H}}\to \tilde X\times \Bun_G$ send a collection $(\tilde x\in\tilde X, M,\cA,V,\cC)\in {^a\Bun_{G, H}}$ to $(V,\cC)$ and $(\tilde x, M,\cA)$ respectively. 
 
 By Theorem~\ref{Th_Hecke_property_Aut_GH}, the latter complex identifies with
\begin{equation}
\label{complex_prf_Th1}
(^a\gp)_!(^a\gq^*K\otimes \H^{\ra}_H(\gRes^{\kappa}(\cS), \Aut_{G, H}))[-\dim\Bun_H]  
\end{equation}
Consider the diagram
$$
\begin{array}{ccccc}
\tilde X\times \Bun_H & \getsup{\supp\times h^{\la}} & {^a\cH_H} & \toup{h^{\ra}} & \Bun_H\\
\uparrow\lefteqn{\scriptstyle \id\times\gq} && \uparrow && \uparrow\lefteqn{\scriptstyle {^a\gq}}\\
\tilde X\times \Bun_{G, H} & \getsup{\supp\times h^{\la}} & {^a_{\tilde X}\cH_{H, G}} & \toup{\supp\times h^{\ra}} & {^a\Bun_{G,H}}\\
 & \searrow\lefteqn{\scriptstyle \id\times\gp} && \swarrow\lefteqn{\scriptstyle {^a\gp}}\\
 && \tilde X\times\Bun_G,
\end{array}
$$
where $^a\cH_H$ is the stack classifying $\tilde x\in\tilde X$, $H$-torsors
$(V,\cC)$ and $(V',\cC')$ on $X$ identified via an isomorphism $\beta$ over $X-\pi(\tilde x)$ so that $\beta$ yields $\cC'\,\iso\, C(a\pi(\tilde x))$. The map $\supp\times h^{\la}$ (resp., $h^{\ra}$) in the top row sends this point to $(\tilde x, V,\cC)$ (resp., to $(V',\cC')$). 

 The stack $^a_{\tilde X}\cH_{H, G}$ is the above diagram classifies collections: $(\tilde x, V,\cC, V',\cC', \beta)\in {^a\cH_H}$, a $G$-torsor $(M,\cA)$ on $X$, and an isomorphism $\cA\otimes\cC\,\iso\,\Omega$. The map $\supp\times h^{\la}$ (resp., $\supp\times h^{\ra}$) is the middle row sends this collection to $(\tilde x, M,\cA, V,\cC)$ (resp., to $(\tilde x, M,\cA, V',\cC')$).
 
 By the projection formulas, now (\ref{complex_prf_Th1}) identifies with
$$
(\id\times\gp)_!(\Aut_{G,H}\otimes (\id\times\gq)^* \H^{\la}_H(\gRes^{\kappa}(\cS), K))[-\dim\Bun_H]
$$ 
%The verification of the compatibilities with the $\Sigma$-actions is left to a reader. 
Theorem~\ref{Th_theta-lifting_functors-general} is proved. 

\appendix

\section{Invariants in the classical setting}
\label{Sect_appendix_A}

\ssec{} In this appendix we assume that $k_0\subset k$ is a finite subfield, and all the objects introduced in Section~4 are defined over $k_0$. In particular, $\cO_0\subset\cO$ is a complete discrete valuation $k_0$-algebra, and $F_0$ its fraction field. Our purpose is to prove Proposition~\ref{PpA.1} formulated in Section~\ref{Pp_about_functor_cP} and Lemma~\ref{Lm_injectivity_appendixA}. 

\begin{Lm} 
\label{Lm_Jean-Francois}
Let $G$ be a reductive group scheme over $\Spec\cO_0$, $P\subset G$ be a parabolic and $U\subset P$ its unipotent radical. Let $V$ be a smooth $\Qlb$-representation of $G(F_0)$. Then the natural map $V^{G(\cO_0)}\to V_{U(F_0)}$ is injective, here $V_{U(F_0)}$ denotes the corresponding Jacquet module.
\end{Lm}
\begin{Prf} Pick a Borel subgroup $B\subset P$, write $I\subset G(\cO_0)$ for the corresponding Iwahori subgroup. It suffices to show that $V^I\to V_{U(F_0)}$ is injective. 
 
 Let $v\in V^I$ vanish in $V_{U(F_0)}$. Then one may find a semisimple $t\in B(F_0)$ such that the characteristic function $\phi$ of $ItI$ annihilates $v$ (it suffices that the action of $t$ on $U(F_0)$ be sufficiently contracting). However, $\phi$ is invertible in the Iwahori-Hecke algebra of $(G(F_0),I)$ by (\cite{AB}, Lemma~8), so $v=0$.
\end{Prf}

\medskip

\begin{Lm} 
\label{Lm_injectivity_appendixA}
The maps 
$J^*_{P_{\HH_a}}: \Weil_a(k_0)\to \cS_{Q\GG\HH_a(\cO_0)}(Q\Upsilon_a(F_0))
$ and 
\begin{equation}
\label{map_for_proof_of_PpA.1}
J^*_{P_{\GG_a}}: \Weil_a(k_0)\to \cS_{Q\GG\HH_a(\cO_0)}(Q\Pi_a(F_0))
\end{equation}
are injective.
\end{Lm}
\begin{Prf}
Both claims being similar, we prove only the second one.  
Apply Lemma~\ref{Lm_Jean-Francois} for the parabolic $P_{\HH_a}\subset \HH_a$ and the representation $\cS(\Pi_a(F))$ of $\cT_a(F)$. Remind that $\cT_a=\{(g_1,g_2)\in \GG_a\times\HH_a\mid (g_1,g_2)\; \mbox{acts trivially on}\; A_a\otimes C_a\}$, and $U_{\HH_a}\subset P_{\HH_a}$ is the unipotent radical. 
  
  For $v\in\Pi_a(F_0)$ let $s_{\Pi}(v): C_a^*\otimes \wedge^2 U_a(F_0)\to\Omega(F_0)$ be the map introduced in Section~\ref{Sect_Parabolic subgroups}. Write $\Cr(\Pi_a)$ for the space of $v\in\Pi_a(F_0)$ such that $s_{\Pi}(v)=0$. 
By (\cite{MVW}, page 72), the Jacquet module $\cS(\Pi_a(F_0))_{U_{\HH_a}(F_0)}$ identifies with the Schwarz space $\cS(\Cr(\Pi_a))$, and the projection 
$$
\cS(\Pi_a(F_0))\to \cS(\Pi_a(F_0))_{U_{\HH_a}(F_0)}
$$ 
identifies with the restriction map $\cS(\Pi_a(F_0))\to \cS(\Cr(\Pi_a))$. So, the restriction map $\Weil_a(k_0)\to \cS_{\GG Q\HH_a(\cO_0)}(\Cr(\Pi_a))$ is injective. Note that $Q\Pi_a(F_0)\subset \Cr(\Pi_a)$. 

 We claim that the restriction $\cS_{\GG Q\HH_a(\cO_0)}(\Cr(\Pi_a))\to \cS_{Q\GG\HH_a(\cO_0)}(Q\Pi_a(F_0))$ is injective. Indeed, let $h$ be in the kernel of this map and $v\in \Cr(\Pi_a)$. We must show that $h(v)=0$. The image of $v: U_a\otimes C_a^*\to M_a(F_0)$ is contained in some lagrangian $F_0$-subspace $\cL\subset M_a(F_0)$. Since the group $\GG_a(\cO_0)$ acts transitively on $\GG_a(F_0)/P_{\GG_a}(F_0)$, and $h$ is $\GG Q\HH_a(\cO_0)$-invariant, we may assume $\cL=L_a(F_0)$, hence we may assume $v\in \Pi_a(F_0)$ and $h(v)=0$.

Thus, (\ref{map_for_proof_of_PpA.1}) is also injective.
\end{Prf}

\bigskip

\begin{Prf}\select{of Proposition~A.1}

\medskip\noindent
For $b\in\ZZ$ set $_b\cH_{\GG}=K_0(_b\Sph_\GG)\otimes\Qlb$ and $_b\cH_{Q(\GG)}=K_0(_b\Sph_{Q(\GG)})\otimes\Qlb$. So, 
$$
\cH_{\GG}=\mathop{\oplus}\limits_{b\in\ZZ} \; {_b\cH_{\GG}},\;\;\;\; \cH_{Q(\GG)}=\mathop{\oplus}\limits_{b\in\ZZ} \; {_b\cH_{Q(\GG)}}
$$
are the Hecke algebras for $\GG$ and $Q(\GG)$ respectively. From Proposition~\ref{Pp_description_modules_for_Levi}, we learn that the map 
$$
_{-a}\cH_{Q(\GG)}\to \cS_{Q\GG\HH_a(\cO)}(Q\Pi_a(F_0))
$$ 
given by $\cS\mapsto \tr_{k_0}\H^{\la}_{Q(\GG)}(\cS, I_0)$ is an isomorphism of $\Qlb$-vector spaces. Write $_{-a}W\subset {_{-a}\cH_{Q(\GG)}}$ for the image of the map (\ref{map_for_proof_of_PpA.1}). We get a $\ZZ$-graded subspace 
$$
\cW:=\mathop{\oplus}\limits_{a\in\ZZ} \;{_aW}\subset \cH_{Q(\GG)}$$

 For $a,a'\in\ZZ$ we have the Hecke operators
$$
\H^{\la}_{\GG}: {_{a'-a}\cH_{\GG}}\times \cS_{\GG Q\HH_{a'}(\cO_0)}(\Pi_{a'}(F_0))\to \cS_{\GG Q\HH_a(\cO_0)}(\Pi_a(F_0))
$$
defined as in Section~\ref{Sect_4.7.1}. We claim that for $\cS\in {_{a'-a}\cH_{\GG}}$ the operator $\H^{\la}_{\GG}(\cS,\cdot)$ sends $\Weil_{a'}(k_0)$ to the subspace  $\Weil_a(k_0)\subset \cS_{\GG Q\HH_a(\cO)}(\Pi_a(F_0))$. This follows from the fact the actions of the groupoids $\GG Q\HH$ and $\HH Q\GG$ on the spaces $\cS_{\GG Q\HH_a(\cO_0)}(\Pi_a(F_0))$ commute with each other. 

 More precisely, for $a,b\in\ZZ$ given $g=(g_1,g_2)\in \cT_{b,a}$ such that $g_2: V_a\,\iso\, V_b$ is an isomorphism of $Q(\HH)$-torsors over $\Spec\cO$, let $h=(h_1,h_2)\in\cT_b$ be any element such that $h_1: M_b\,\iso\, M_b$ is a scalar automorphism of the $\GG$-torsor $M_b$ over $\Spec\cO$. Here $h_2$ is an automorphism of the $\HH$-torsor $V_b$ over $\Spec\cO$. Set $h'_2=g_2^{-1} h_2 g_2$, so $h'_2$ is an automorphism of the $\HH$-torsor $V_a$ over $\Spec\cO$. Set $h'_1=h_1$ then $h'=(h_1,h_2)\in \cT_a$. The equality $g h'=hg$ in $\cT$ shows that $g: \cS(\Pi_a(F_0))\to \cS(\Pi_b(F_0))$ sends $\HH_a(\cO_0)$-equivariant objects to $\HH_b(\cO_0)$-equivariant objects. We have used the action of the groupoid $\cT$ on the spaces $\cS(\Pi_a(F_0))$ obtained as in Remark~\ref{Rem_groupoid_action_two_Shrodinger_models}. 

 Thus, $\cW$ is a $\ZZ$-graded module over the $\ZZ$-graded ring $\cH_{\GG}$. We also know from (\cite{L2}, Proposition~2) that $_0W={_0\cH_{\GG}}$. 
Our statement is reduced to Lemma~\ref{Lm_graded_modules_appendix} below.  
\end{Prf}

\medskip

 Remind that we have picked a maximal torus $T_{\GG}\subset Q(\GG)$. Write $W$ (resp., $W_Q$) for the Weyl group of $(\GG,T_{\GG})$ (resp., of $(Q(\GG),T_{\GG})$). Then 
$$
\cH_{Q(\GG)}\,\iso\, \Qlb[\check{T}_{\GG}]^{W_Q},\;\;\;\;\;\;\;
\cH_{\GG}\,\iso\, \Qlb[\check{T}_{\GG}]^W
$$ 
Recall the map $\kappa_{\GG}: \check{Q}(\GG)\times\Gm\to \check{\GG}$ from Section~\ref{Sect_4.9.4}. The homomorphism $\Res^{\kappa_{\GG}}:\cH_{\GG}\to\cH_{Q(\GG)}$ comes from the map
$f^{\kappa_{\GG}}: \check{T}_{\GG}^{W_Q}\to \check{T}_{\GG}^W$ obtained by taking the Weil group invariants of the map
$\check{T}_{\GG}\to\check{T}_{\GG}$, $t\mapsto t \nu(q^{1/2})$, where $\nu$ is some coweight of the center $Z(\check{Q}(\GG))$, and $q$ is the number of elements of $k_0$. 

\begin{Lm} 
\label{Lm_graded_modules_appendix}
View $\cH_{Q(\GG)}$ as a $\ZZ$-graded $\cH(\GG)$-module via $\Res^{\kappa_{\GG}}:\cH_{\GG}\to\cH_{Q(\GG)}$. 
Let 
$$
\cW=\mathop{\oplus}\limits_{a\in\ZZ} {_aW}\subset \cH_{Q(\GG)}=\mathop{\oplus}\limits_{a\in\ZZ} {_a\cH_{Q(\GG)}}
$$ 
be a $\ZZ$-graded submodule over the $\ZZ$-graded ring $\cH_{\GG}$. Assume that $_0W={_0\cH_{\GG}}$. Then $\cW=\cH_{\GG}$.
\end{Lm} 
\begin{Prf}
Given $x\in {_aW}$, pick a nonzero $h\in {_{-a}\cH_{\GG}}$ then $hx\in {_0\cH_{\GG}}$. So, $x$ is a rational function on $\check{T}_{\GG}^W$ which becomes everywhere regular after restriction under $f^{\kappa_{\GG}}: \check{T}_{\GG}^{W_Q}\to \check{T}_{\GG}^W$. Since $\check{T}_{\GG}^W$ is normal by Remark~\ref{Rem_normality} below, and 
$x$ is entire over $\Qlb[\check{T}_{\GG}]^W$, it follows that $x\in \Qlb[\check{T}_{\GG}]^W$.
\end{Prf}

\medskip

\begin{Rem} 
\label{Rem_normality}
Let $A$ be an entire normal ring, $W$ be a finite group acting on $A$. Assuming that $A$ is finite over $A^{W}$, one checks that $A^W$ is normal. 
\end{Rem}


\begin{thebibliography}{99}
\bibitem{AGKRRV} D. Arinkin, D. Gaitsgory, D. Kazhdan, S. Raskin, N. Rozenblyum, Y. Varshavsky, The stack of local systems with restricted variation and geometric Langlands theory with nilpotent singular support, arXiv:2010.01906
\bibitem{AB} S. Arkhipov, R. Bezrukavnikov, Perverse sheaves on affine flags and langlands dual group, Israel J. Math. 170, 135 (2009), 135 - 183
\bibitem{BBD} A. Beilinson, J. Bernstein, P. Deligne, Faisceaux pervers, Asterisque 100 (1982)
\bibitem{BD} A.~Beilinson, V.~Drinfeld, Quantization of the Hitchin integrable system and Hecke eigensheaves, preprint downloadable from\\ {http://www.math.utexas.edu/$\sim$benzvi/Langlands.html}
\bibitem{B} A. Borel, Automoprhic $L$-functions, Proc. of Simposia in Pure Math., 33 (1979), part 2, 27 - 61
\bibitem{Brad} T. Braden, Hyperbolic localization of intersection cohomology, Transform. Groups 8 (2003), 209 - 216
\bibitem{BG} A. Braverman, D. Gaitsgory, Geometric Eisenstein series, Inv. Math. 150 (2002), 287 - 384
\bibitem{FGV} E. Frenkel, D. Gaitsgory, K. Vilonen, On the geometric Langlands conjecture, J. Amer. Math. Soc. 15 (2002), 367 - 417
\bibitem{G5} D. Gaitsgory, The local and global versions of the Whittaker category, arXiv:1811.02468 version 6
\bibitem{GL} D. Gaitsgory, S. Lysenko, Parameters and duality for the metaplectic geometric Langlands theory, Selecta Math., (2018) Vol. 24, Issue 1, 227 - 301. Corrected version: of August 30, 2020 on my webpage and arxiv.
\bibitem{GaRo} D. Gaitsgory, N. Rozenblyum, A study in derived algebraic geometry, Math. surveys and monographs, vol. 221 (2017), AMS Providence, Rhode Island, corrected version available at {https://people.math.harvard.edu/$\sim$gaitsgde/GL/}
\bibitem{GP} Ph. Gille, Patrick Polo,  R\'e\'edition de SGA3, available at \\ https://webusers.imj-prg.fr/$\sim$patrick.polo/SGA3/
\bibitem{HNY} J. Heinloth, B.-Ch. Ng\^o, Zh. Yun, Kloosterman sheaves for reductive groups, Ann. Math., 177 (2013), 241 - 310
\bibitem{La} G. Laumon, Transformation de Fourier, constantes d'\'equations fonctionnelles et conjectures de Weil, Publ. IHES, t. 65 (1987), p. 131 - 210
\bibitem{L} S. Lysenko, Geometric Waldspurger periods, Compos. Math. 144 (2008), no. 2, p. 377 - 438 
\bibitem{L1} S. Lysenko, Moduli of metaplectic bundles on curves and Theta-sheaves, Ann. Scient. \'Ec. Norm. Sup. 4 s\'erie, t.39 (2006), 415 - 466
\bibitem{L2} S. Lysenko, Geometric theta-lifting for the dual pair $\SO_{2m}$, $\Sp_{2n}$, math.RT/0701170 version 7 (this is a corrected version of S. Lysenko, Geometric theta-lifting for the dual pair $\SO_{2m}$, $\Sp_{2n}$, Annales ENS, 4e serie, t. 44 (2011), p. 427 - 493)
\bibitem{L3} S. Lysenko, On the automorphic sheaves for $\GSp_4$, arXiv:1901.04447
\bibitem{LL} V. Lafforgue, S. Lysenko, Geometric Weil representation: local field case, Compos. Math. 145 (2009), no. 1, 56 - 88 
\bibitem{LL1} V. Lafforgue, S. Lysenko,  Compatibility of the Theta correspondence with the Whittaker functors,  Bull. Soc. math. France, 139 (1), (2011), 75 - 88 
\bibitem{HTT} J. Lurie, Higher topos theory, Ann. of math. studies, no. 170, Princeton University Press (2009)
\bibitem{HA} J. Lurie, Higher algebra, version September 18, 2017
\bibitem{MV} I. Mirkovi\'c, K. Vilonen, Geometric Langlands duality and representations of algebraic groups over commutative rings, Ann. of Math., 166 (2007), 95 - 143
\bibitem{MVW} C. Moeglin, M.-F. Vign\'eras, J.-L. Waldspurger, Correspondance de Howe sur un corps p-adique, Lecture Notes in Math. 1291 (1987)
\bibitem{Ral} S. Rallis, Langlands functoriality and the Weil representation, Amer. J. Math., vol.104, no. 3, p. 469 - 515 (1982)
\bibitem{R} B. Roberts, The theta correspondence for similitudes, Israel J. Math. 94 (1996), 285 - 317
\bibitem{Zh} X. Zhu, The geometric Satake correspondence for ramified groups, Annales Sci. ENS, 48 (2), (2015), p. 409 - 451
\end{thebibliography}
\end{document}